\documentclass[11pt]{article}

\usepackage{grffile}
\usepackage{ctable}

\usepackage{amssymb}
\usepackage{amsmath}
\usepackage[dvips]{epsfig}
\usepackage[small]{caption}
\usepackage{graphicx}
\usepackage[all]{xy}

\newtheorem{Lemma}{Lemma}[section]

\newtheorem{Proposition}[Lemma]{Proposition}

\newtheorem{Remark}[Lemma]{Remark}
\newtheorem{Definition}[Lemma]{Definition}
\newtheorem{Hypothesis}[Lemma]{Hypothesis}

\newenvironment{Proof}%
 {\begin{trivlist} \item[]{\bf Proof. }}%
 {\hspace*{\fill}$\rule{.4\baselineskip}{.4\baselineskip}$\end{trivlist}}

 {\begin{trivlist}\item[]\textbf{Acknowledgments.}}{\end{trivlist}}

\makeatletter\@addtoreset{figure}{section}\makeatother

\makeatletter \@addtoreset{equation}{section} \makeatother

\newcommand{\R}{\mathbb{R}}
\newcommand{\C}{\mathbb{C}}

\newcommand{\Z}{\mathbb{Z}}

\def\Re{\mathop{\mathrm{Re}}}
\def\Im{\mathop{\mathrm{Im}}}

\newcommand{\rmO}{\mathrm{O}}

\newcommand{\rmd}{\mathrm{d}}
\newcommand{\rme}{\mathrm{e}}
\newcommand{\rmi}{\mathrm{i}}

\newcommand{\Rg}{\mathrm{Rg}}

\renewcommand{\leq}{\leqslant}
\renewcommand{\geq}{\geqslant}

\def\XXint#1#2#3{{\setbox0=\hbox{$#1{#2#3}{\int}$}
     \vcenter{\hbox{$#2#3$}}\kern-.5\wd0}}

\newfam\bifam
\font\tenbi=cmmib10 scaled \magstep1 \font\sevenbi=cmmib10 at 11pt
\font\fivebi=cmmib10 at 6pt \textfont\bifam = \tenbi
\scriptfont\bifam = \sevenbi \scriptscriptfont\bifam= \fivebi

\ifx\pdfoutput\undefined
   \pdffalse
\else
   \pdfoutput=1
   \pdftrue
\fi
\ifpdf
   \usepackage{graphicx}
   \usepackage{epstopdf}
\else
   \usepackage{graphicx}
\fi

\begin{document}
\vspace*{0.4in}
\begin{center}
{\fontsize{16}{16}\fontfamily{cmr}\fontseries{b}\selectfont{Characterizing the effect of boundary conditions on striped phases}}\\[0.2in]

David Morrissey and  Arnd Scheel\footnote{Research partially supported by the National Science Foundation through  grants NSF- DMS-0806614 and NSF-DMS-1311740.}\\
\textit{\footnotesize University of Minnesota, School of Mathematics,   206 Church St. S.E., Minneapolis, MN 55455, USA}

\date{\small \today} 
\end{center}

\begin{abstract}
\noindent 
We study the influence of boundary conditions on stationary, periodic patterns in one-dimensional systems. We show how a conceptual understanding of the structure of equilibria in large domains can be based on the characterization of boundary layers through displacement-strain curves. Most prominently, we distinguish wavenumber-selecting and phase-selecting boundary conditions and show how they impact the set of equilibria as the domain size tends to infinity. We illustrate the abstract concepts in the phase-diffusion and the Ginzburg-Landau approximation. We also show how to compute displacement-strain curves in more general systems such as the Swift-Hohenberg equation using continuation methods. 
\end{abstract}

%

\vspace*{0.2in}

{\small
{\bf Running head:} {Boundary layers for striped phases}

}
\vspace*{0.2in}

%
%
%
%

\section{Introduction}

Stripe patterns are arguably the simplest non-trivial patterns observed in nature. Stripes form in numerous contexts, ranging from classical fluid experiments such as Taylor vortices in the Taylor-Couette flow and convection rolls in Rayleigh-B\'enard convection, over bilayers in diblock copolymers, to granular media, and reaction-diffusion systems. A well-studied context of stripe formation is an instability of a spatially homogeneous state in a spatially extended medium. A linearized analysis can then predict wavenumbers as linearly fastest growing modes. 
While known for more than a century in fluid dynamics, the occurrence of patterned states as fastest growing modes was conjectured by Turing in 1952 \cite{turing} for reaction-diffusion systems and experimentally realized in 1991 \cite{dekepper}. 

The arguably simplest example for the formation of stripes is the celebrated Swift-Hohenberg equation
\[
u_t=-(\Delta+1)^2u + \mu u - u^3.
\]
Considered on $x\in\R$, for instance, the linearization at the trivial state $u\equiv 0$,
\[
u_t=-(\Delta+1)^2u + \mu u,
\]
can be readily analyzed using Fourier transform,
\[
\frac{\rmd}{\rmd t}\hat{u}(k)=(-(1-k^2)^2+\mu)\hat{u}(k).
\]
For $\mu>0$, wavenumbers with $|k^2-1|<\sqrt{\mu}$ are unstable. The fastest-growing wavenumber is $|k|=1$. In fact, for all $\mu<1$, \emph{only} wavenumbers with $k\neq 0$ are unstable, so that for small amplitude only spatially patterned perturbations are amplified. 

A nonlinear analysis readily reveals that together with the linear instability, nonlinear patterns bifurcate from the trivial state. For $\mu>0$ small, one finds $2\pi/k$-periodic ``striped'' solutions $u_\mathrm{st}(kx;k,\mu)$, $u_\mathrm{st}(\xi;k,\mu)=u_\mathrm{st}(\xi+2\pi;k,\mu)$, for all unstable $k$-vectors, $|k^2-1|<\sqrt{\mu}$. Many of those periodic solutions are compatible with the boundary conditions in large but finite systems. With, say, periodic boundary conditions and system size $L$, we need to require $k\in (2\pi/L)\Z$, which yields $\rmO(L\sqrt{\mu})$ different $k$-values. 

For boundary conditions other than Neumann, Dirichlet, or periodic, equilibria are difficult to characterize completely. In fact, bifurcations diagrams tend to be very complex, and of little help when trying to ``describe'' typical dynamics. Interestingly, different types of boundary conditions can provoke quite different types of dynamics. For instance, heated walls in Rayleigh-B\'enard convection are known to favor parallel alignment of roll solutions with the boundary, as opposed to the otherwise typical perpendicular alignments. Slow drift of rolls has also been observed in contexts of Rayleigh-B\'enard convection. 

\begin{figure}
\centering{
\includegraphics[width=0.49\textwidth]{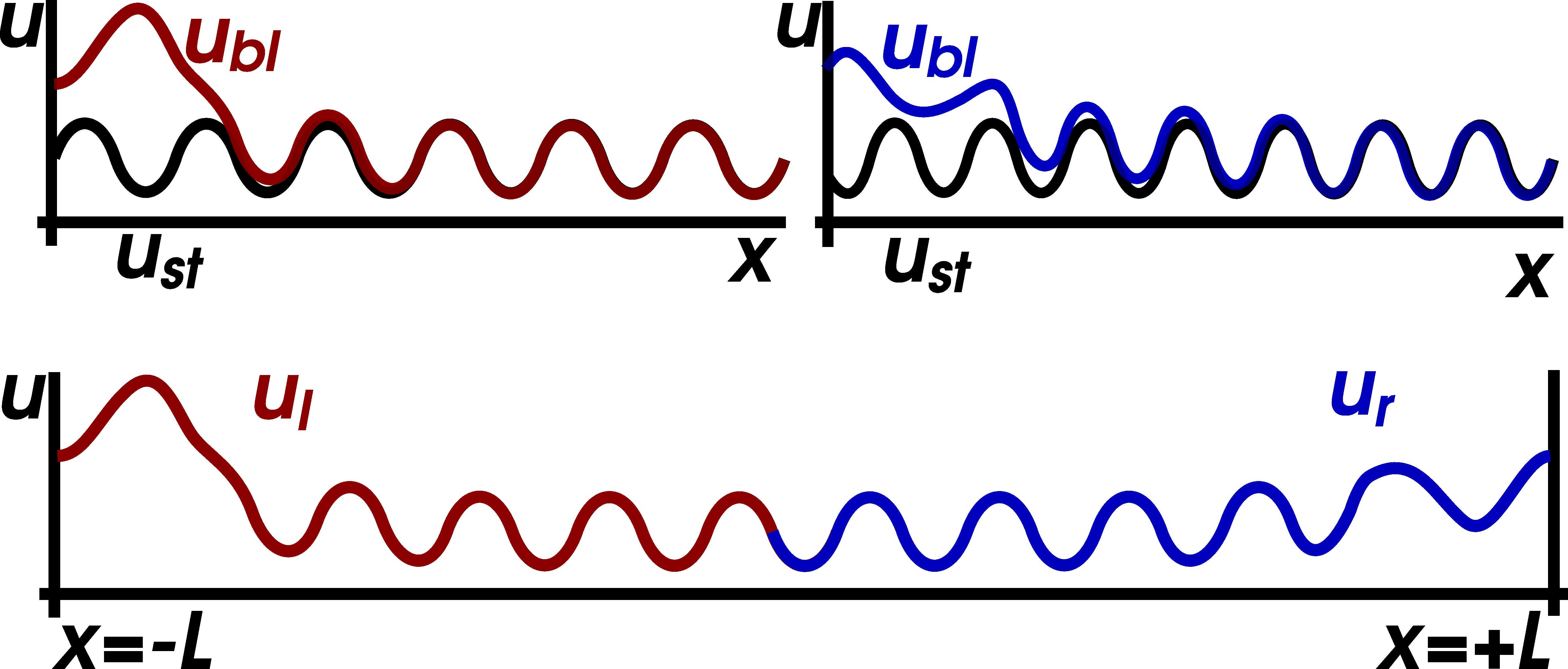}\hfill\includegraphics[width=0.45\textwidth]{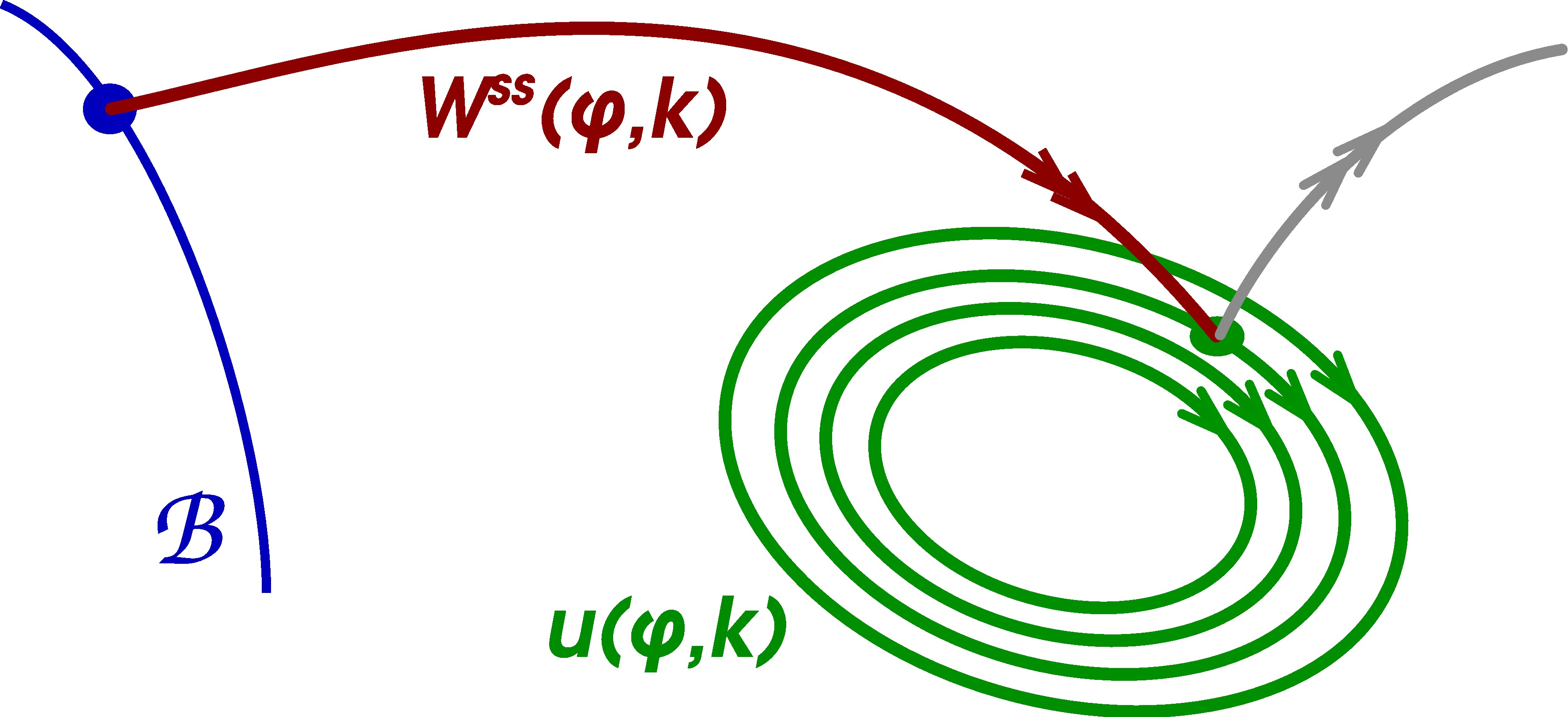}
}
\caption{Boundary layers selecting different wavenumbers and phases of the periodic pattern $u_\mathrm{st}$ (top left); boundary layers for left and right boundary, $u_\mathrm{l/r}$, matched at $x=0$ (bottom left). Schematic plot of boundary manifold $\mathcal{B}$, strong stable fiber $W^\mathrm{ss}$ intersecting $\mathcal{B}$ and thus yielding a boundary layer, the family of periodic orbits, and the unstable fibers (right).}\label{f:sch}
\end{figure}

The purpose of this work is to study the role of boundaries in such pattern forming systems in a systematic fashion. The key idea is to decompose the problem of describing equilibria in large finite domains into two steps (see Figure \ref{f:sch} for an illustration):
\begin{enumerate}
\item describe equilibria in semi-infinite domains $x\in\R_+$ and $x\in\R_-$;
\item patch equilibria from semi-infinite domains at $x=0$ to obtain equilibria in bounded domains $x\in(-L,L)$. 
\end{enumerate}
The boundary conditions in step (ii) at $x=L$ should be studied as boundary conditions at $x=0$ in step (i) for $x\in\R_-$, boundary conditions at $x=-L$ give boundary conditions for $x\in\R_+$.

\paragraph{Outline.} We study a general class of equations in Section \ref{s:2}. In particular, we characterize stripes, introduce spatial dynamics, define boundary layers, and discuss displacement-strain relations for boundary layers. We also present our main general matching result for bounded domains. Section \ref{s:2a} illustrates the concepts at the hands of the (somewhat trivial) example of the phase-diffusion equation. In Section \ref{s:3}, we study the specific (integrable) example of the Ginzburg-Landau equation, which arises as a modulation equation for small-amplitude striped patterns. Section \ref{s:4} outlines an effective numerical procedure for the computation of displacement-strain relations and illustrates results for the Swift-Hohenberg equation. We conclude with a discussion in Section \ref{s:5}. 

\section{Boundary layers for striped phases: wavenumber and phase selection}\label{s:2}
We define families of stable stripes in Section \ref{s:2.1}, and  introduce spatial dynamics in Section \ref{s:2.2}. We introduce boundary layers in Section \ref{s:2.3}. We briefly explore stability in Section \ref{s:2.4}. Section \ref{s:2.5} is concerned with gluing boundary layers to obtain equilibrium configurations in large but finite domains. Section \ref{s:2.6} formalizes the gluing process as a subtraction of curves followed by quantization. Section \ref{s:2.7} contains a list of conceptual examples. 

\subsection{Families of stripes}\label{s:2.1}
We consider a general semilinear parabolic system on the real line
\begin{equation}\label{e:par}
u_t=-(\rmi\partial_{x})^{2m}u+f(u,\partial_x u,\ldots,\partial_x^{2m-1}u),\qquad u\in\R^N,x\in\R
\end{equation}
where $f$ is smooth and reflection symmetric, 
\[
f(u,u_1,u_2,\ldots,u_{2m-1})=f(u,-u_1,u_2,\ldots,-u_{2m-1}).
\]
Such systems may possess \emph{even, periodic solutions} $u_\mathrm{st}$, which we refer to as stripes. The linearization at such stripes,
\[
\mathcal{L}_\mathrm{st}v=-(\rmi\partial_{x})^{2m}v+\sum_{j=1}^{2m}\partial_jF(u,\partial_x u_\mathrm{st},\ldots,\partial_x^{2m-1}u_\mathrm{st})\partial_x^{j-1}v,
\]
defines an elliptic operator on $L^2(\R,\R^N)$. The spectrum of $\mathcal{L}$ is the union of the point spectra of the associated Bloch operators,
\[
\mathcal{L}_{\mathrm{st},\ell} v=-(\rmi(k\partial_{\xi}+\rmi\ell))^{2m}v+\sum_{j=1}^{2m}\partial_jf(u,\partial_x u_\mathrm{st},\ldots,\partial_x^{2m-1}u_\mathrm{st})(k\partial_\xi+\rmi\ell)^{j-1}v,
\]
where the operators $\mathcal{L}_{\mathrm{st},\ell}$ are defined on $L^2((0,2\pi),\R^N)$ with periodic boundary conditions. Since $u_\mathrm{st}'(kx;k)$ contributes to the kernel of $\mathcal{L}_0$, one always finds a band of eigenvalues $\lambda(\ell)$, $\lambda(0)=0$, to $\mathcal{L}_\ell$. Symmetry shows that $\lambda(\ell)=\lambda(-\ell)\in\R$. We therefore say that a periodic pattern is stable if the spectrum of $\mathcal{L}$ is contained in $\Re\lambda<0$, except for the critical branch of eigenvalues $\lambda(\ell)$, for which we require $\lambda''(0)<0$. Such spectral stability has been shown to imply nonlinear stability with respect to localized perturbations and we will henceforth refer to this criterion simply as ``stability''; see \cite{schneidderguido,johnsonzumbrun,scheelwu}. Instabilities can occur due to other bands of spectrum crossing the imaginary axis, or due to an Eckhaus instability, when $\lambda''(0)$ changes sign \cite{turingjens}.

\begin{Hypothesis}[Families of stripe solutions]\label{h:1}
We assume that there exists a bounded family of stable stripes $u_\mathrm{st}(kx;k)$,for wavenumbers $k\in J_k=(k_\mathrm{min},k_\mathrm{max})$, $0<k_\mathrm{min}<k_\mathrm{max}<\infty$,
\[
u_\mathrm{st}(\xi;k)=u_\mathrm{st}(\xi+2\pi;k), \quad u_\mathrm{st}(\xi;k)=u_\mathrm{st}(-\xi;k)\quad u_\mathrm{st}(\cdot;k)\in C^2(J_k,C^{2m}_\mathrm{per}(0,2\pi)).
\]
\end{Hypothesis}
We emphasize that the interval $J_k$ is not assumed to be small. In the examples that we consider below, the boundaries $k_\mathrm{min/max}$ are determined by Eckhaus instabilities. Most of our analysis is insensitive to oscillatory instabilities, $\lambda\in\rmi\R\setminus\{0\}$ so that in this case the effective range of wavenumbers could be further extended. 

\begin{Lemma}[Robustness]\label{l:1}
Given a stable stripe, there exists a family of stripes nearby. Moreover, families of stripes depend smoothly on system parameters.
\end{Lemma}
\begin{Proof}
The kernel of the Bloch-wave operator $\mathcal{L}_0$ is trivial in $L^2_\mathrm{even}$, so that one can continue a stable periodic pattern smoothly in the period. Since the condition $\lambda''(0)<0$ is open, all members in the local family will be stable. One similarly establishes persistence of global families with respect to parameters.
\end{Proof}
\begin{Remark}[Minimal period]\label{r:min}
We may assume throughout that the period $2\pi$ is the minimal period. We may choose the period to be minimal for a fixed member in the family. Changes of the minimal period always imply a bifurcation, hence would imply additional spectrum at $\lambda=0$ and contradict our assumption on stability. As a consequence, the stripes have precisely two reflection symmetries, $x\to -x$ and $x\to 2\pi-x$. 
\end{Remark}

\subsection{Spatial dynamics}\label{s:2.2}

The steady-state equation 
\[
0=-(\rmi\partial_{x})^{2m}u+f(u,\partial_x u,\ldots,\partial_x^{2m-1}u),
\]
can be rewritten as a first-order system,
\begin{equation}\label{e:ode}
U_x=F(U),\qquad U=(u,\partial_xu,\ldots,\partial_x^{2m-1}u)\in\R^{2mN},
\end{equation}
which possesses a reversibility symmetry 
\[
F(RU)=-RF(U),\quad  R(u_0,u_1,\ldots,u_{2m-1})=(u_0,-u_1,\ldots,-u_{2m-1}).
\]
Stripes correspond to reversible periodic solutions $RU(kx;k)=U(-kx;k)$. Purely imaginary Floquet exponents of the linearization at such a periodic solution give rise to kernels of the operators $\mathcal{L}_\ell$ in a straight-forward fashion. Moreover, the algebraic multiplicity of the simple zero-exponent is given by the order order of tangency of the curve of critical spectrum $\lambda(\ell)$; see for instance \cite[Lemma 2.1]{radial}. The family of stripes therefore forms a two-dimensional normally hyperbolic manifold $\mathcal{S}$ in $\R^{2mN}$. Reversibility implies that both stable and unstable manifold $W^\mathrm{cs/cu}$ are $mN+1$-dimensional, and mapped into each other by $R$. Both are smoothly fibered by $mN-1$-dimensional strong (un)stable fibers $W^\mathrm{ss/uu}$. 

\subsection{Boundary layers}\label{s:2.3}

We consider (\ref{e:par}) on $x\in\R_+$, supplemented with $mN$ boundary conditions, which we write in the form $U(x=0)\in\mathcal{B}$, where $\mathcal{B}$ is an $mN$-dimensional smooth manifold in $\R^{2mN}$. Natural examples are 
\begin{itemize}
\item Dirichlet, $\partial_x^{2j}u=a_j, 0\leq j<mN$;
\item Neumann, $\partial_x^{2j+1}u=a_j, 0\leq j<mN$;
\item clamped,  $\partial_x^{j}u=a_j, 0\leq j<mN$;
\item free,  $\partial_x^{j}u=a_j, mN\leq j<2mN$.
\end{itemize}
The terms clamped and free allude to the plate equation. Of course, many mixed and nonlinear variations are possible, and not all guarantee well-posedness of the parabolic equation (\ref{e:par}).

Equilibria of (\ref{e:par}) in $BC^{2m}(\R_+,\R^N)$ that satisfy the boundary conditions $\mathcal{B}$ are bounded solutions $U(x)$ to the ODE (\ref{e:ode}), with $U(0)\in\mathcal{B}$, where $\mathcal{B}$ is now considered as a natural $mN$-dimensional submanifold of the phase space of $\R^{2mN}$ of (\ref{e:ode}). In other words, equilibria are naturally identified via a multiple shooting problem; see Figure \ref{f:sch}.  Equilibria of interest to us converge to stripes as $x\to\infty$. They therefore lie in the intersection $\mathcal{B}\cap W^\mathrm{cs}$. Simple dimension counting suggests that this intersection is one-dimensional. Thom transversality indeed guarantees just this for generic manifolds $\mathcal{B}$. Moreover, the intersection projects smoothly along the fibers onto a curve in $\mathcal{S}$. Parameterizing $\mathcal{S}$ naturally by $k\in J_k$ and an angle $\varphi$, we can now characterize boundary layers through curves in the cylinder $\varphi,k\in S^1\times J_k$. 

\begin{Definition}[Boundary layers and displacement-strain relations]\label{d:bl}
We say a boundary layer $U_\mathrm{bl}(\xi)$ is transverse if 
$
T_{U_\mathrm{bl}(0)}\mathcal{B}\pitchfork T_{U_\mathrm{bl}(0)}W^\mathrm{cs} 
$
is one-dimensional, not contained in $ T_{U_\mathrm{bl}(0)}W^\mathrm{ss}$. We refer to the base point of the  fiber $W^\mathrm{ss}$ that contains $\mathcal{U_\mathrm{bl}}(0)$ as the asymptotic phase $\varphi$ and wavenumber $k$.
If all boundary layers are transverse, we refer to the collection of base points as \emph{displacement-strain ($k,\varphi$) curves} $\gamma$. Unless otherwise noted, we parameterize displacement-strain curves by arc length $\gamma(s)$. We refer to a relation $d(k,\varphi)$ which vanishes precisely at $\gamma$ as displacement-strain relation; see (\ref{e:bl}) for an explicit formula.
\end{Definition}
\begin{Remark}\label{r:ds}
\begin{itemize}
\item The terminology displacement-strain refers to the fact that the asymptotic stripes are ``displaced'' by $\varphi$ relative to a fixed, even, stripe pattern on $x\in\R$. Strain refers to the intuitive compression and expansion of the stripe relative to a fixed wavenumber.
\item 
The definition of $\varphi$ in $\mathcal{S}$ is unique with the convention $\mathcal{S}=\{U_\mathrm{st}(kx-\varphi;k),\ k\in J_k\}$, when requiring $U_\mathrm{st}(\cdot,k)$ to be even and after fixing $U_\mathrm{st}(\cdot;k_0)$ (which rules out a shift by $\pi$). We then have
\begin{equation}\label{e:bl}
|U(x)-U_\mathrm{st}(kx-\varphi;k)|\leq C\rme^{-\eta x},
\end{equation}
for some positive constants $C,\eta$. 
\item In a completely analogous fashion, one can define boundary layers for right-bounded domains, $x\in\R_-$, $U(0)\in \mathcal{B}_+$, and 
\begin{equation}\label{e:blr}
|U(x)-U_\mathrm{st}(kx-\varphi;k)|\leq C\rme^{\eta x}.
\end{equation}
Note that boundary layers in left-bounded domains give boundary layers in right-bounded domains via reflection, but for \emph{reflected boundary conditions} $R\mathcal{B}_+=\mathcal{B}_-$, replacing $\partial_x$ by $-\partial_x$. 
\item 
We emphasize that the curve $\gamma$ may have self-intersections, stemming from the projection of a smooth non-intersecting curve along the strong stable fibration. 
\end{itemize}
\end{Remark}

Given the (generic) transversality in Definition \ref{d:bl}, we can 
envision families of boundary layers, globally in $k,\varphi$. Note however that even when all boundary layers lie in a bounded set in phase space $\R^{2mN}$, the resulting displacement-strain curves may still terminate, for instance when the $\varphi-k$-curve approaches the boundary of $J_k$,or when the ``length'' of the heteroclinic diverges; see \cite{homsan} for some context on homoclinic and heteroclinic continuation. We will in fact encounter both possibilities in the example of the Ginzburg-Landau equation, Section \ref{s:3}. 

\begin{Definition}[Phase- and wavenumber selection]\label{d:sel}
A displacement-strain curve $\gamma\subset S^1\times J_k$ encodes 
\begin{itemize}
\item \emph{wavenumber selection}, if $P_k\gamma\neq J_k$;
\item \emph{phase selection}, if $P_\varphi\gamma\neq S^1$.
\end{itemize}
Here, $P_{\varphi/k}$ denote projections onto the $\varphi$- or $k$-component, respectively. If $\gamma$ is a closed curve, we write $i(\gamma)$ for the winding number of the curve $P_{\varphi}\gamma$ in $S^1$.
\end{Definition}
We will commonly encounter winding numbers 0 and 1, but also an example if winding number 2 in the Ginzburg-Landau example. 

We will explore a definition of boundary layers not based on spatial dynamics in Section \ref{s:4}. This definition will prove useful when constructing boundary layers numerically. 

\subsection{Stability of boundary layers}\label{s:2.4}
Without any further assumptions, little can be said about the stability of boundary layers. Given the stability of the asymptotic stripe pattern, we can say that the essential spectrum is contained in the negative half plane with the exception of a branch of essential spectrum touching the imaginary axis; see \cite{fiesch}. The mere existence does, of course, not give information on possible oscillatory instabilities, but we do have information on possible zero eigenvalues. We therefore count the parity of unstable eigenvalues of the linearization at a boundary layer $U_\mathrm{bl}=(u_\mathrm{bl},\ldots,\partial_x^{2m-1}u_\mathrm{bl})$,
\[
\mathcal{L}_\mathrm{bl}v=-(\rmi\partial_{x})^{2m}v+\sum_{j=1}^{2m}\partial_jF(u,\partial_x u_\mathrm{bl},\ldots,\partial_x^{2m-1}u_\mathrm{bl})\partial_x^{j-1}v,\qquad x>0,
\]
equipped with boundary conditions $(v,\partial_x v,\ldots,\partial_x^{2m-1}v)\in T\mathcal{B}$. In other words, $p=1$ refers to an even (or zero) number of unstable eigenvalues of $\mathcal{L}_\mathrm{bl}$ and $p=-1$ to an odd number. 

The following lemma relates the infinitesimal displacement-strain relations to stability.

\begin{Lemma}[Parity and displacement-strain]\label{l:par}
For any family of boundary layers and corresponding displacement-strain curve $\gamma(s)$, there exists a parity $e\in\{\pm 1\}$ so that 
\[
p=e\cdot \mathrm{sign}\,P_k\left(\frac{\rmd \gamma}{\rmd s}\right).
\]
\end{Lemma}
\begin{Proof}
We need to show that $\lambda=0$ belongs to the extended point spectrum precisely when $k'=0$ along a family of boundary layers, and that we have a strict crossing whenever $k'$ changes sign. 

Under our assumption on symmetry and stability of periodic patterns, the essential spectrum in the neighborhood of the origin is given by the negative real axis $\lambda\leq0$. We consider the linearized equation at a boundary layer, with spectral parameter $\lambda=\gamma^2$. Since coefficients converge exponentially, there is a complex linear change of variables, analytic in $\gamma$ and smoothly depending on $x$, so that the strong stable ($Nm-1$-dimensional), strong unstable ($Nm-1$-dimensional), and center subspaces $2$-dimensional are independent of $x$ and $\gamma$. in these coordinates, the (linearized) boundary conditions at $x=0$ define an $Nm$-dimensional subspace that depends analytically on $\gamma$. We list vectors in these coordinates as $(U^\mathrm{c},U^\mathrm{ss},U^\mathrm{uu})^T$. Since $\gamma=0$ is a simple branch point of the dispersion relation, we can choose coordinates in the center subspace such that the stable subspace is spanned by the $Nm-1$ column vectors in $(0,\mathrm{id},
0)^T$ and the vector $((1,-\gamma),0,0)^T+\rmO(\gamma^2)$.  The boundary conditions are spanned by the $Nm-1$ column vectors in $(0,0,\mathrm{id})^T$ and $(Q
\frac{\rmd \gamma}{\rmd s},0,0)^T$, where $Q$ is an invertible linear transformation. We can form the determinant $\mathcal{E}(\gamma)$ from basis vectors of stable eigenspace and boundary tangent space to build the Evans function which then tracks eigenvalues near the origin $\gamma=0$. Eigenvalues in $\gamma>0$ correspond to actual unstable eigenvalues, whereas $\gamma<0$ corresponds to resonance poles. Since the stable subspace limits on the eigenspace at $\gamma=0$, which is spanned by the derivative of the periodic pattern, we find that $Q\gamma'=(q_{11}\varphi'+q_{12}k',q_{22}k')^T$, with $q_{11},q_{22}\neq 0$ from invertibility. We therefore find that 
\[
\mathcal{E}(\gamma)=\left(q_{22}k'-q_{11}\varphi'\gamma+\rmO(\gamma^2)\right)\mathcal{E}_0(\gamma),
\]
where $\mathcal{E_0}$ is analytic and nonzero near the origin. We conclude that $\gamma=0$ is a root precisely when $k'=0$, and the sign of $\gamma$ changes precisely when the sign of $k'$ changes. 
\end{Proof}

\begin{Remark}[Nonlinear stability]
Given the results on stability of stripes on $\R$, one expects spectrally stable boundary layers to be stable. In fact, one would expect a slightly faster pointwise decay, $\sim t^{-3/2}$, given that the Evans function does not vanish at the origin for boundary layers, but does vanish for stripes. This is in analogy to the heat equation with, say, Dirichlet boundary conditions on the half line.
\end{Remark}

\begin{Remark}[Orientation]
We can orient displacement curves so that $p=\mathrm{sign}\,\left(P_k\left(\frac{\rmd \gamma}{\rmd s}\right)\right)$. Stable boundary layers within a connected family then necessarily either compress or expand along $\gamma$. The somewhat trivial examples of phase diffusion, below, illustrate that there is no obvious way to avoid an ambiguity in orientation (or, equivalently, choice of $e$ in Lemma \ref{l:par}).
\end{Remark}

\begin{Remark}[Displacement-strain dependence]
The results on stability suggest that, typically, $k$ could be thought of as a parameter so that turning points correspond to exchange of stability as in a traditional saddle-node bifurcation. It is mainly this motivation that leads us to displaying displacement-strain relations with $k$ as the independent, horizontal variable in the remainder of this paper. 
\end{Remark}

\subsection{Finite domains --- gluing boundary layers}\label{s:2.5}

In this section, we describe equilibria in large but finite domains, based on left and right boundary layers. The key conclusion is that displacement strain curves for boundary layers associated with boundary conditions on the left and right side of the domain allow for a complete description of the set of equilibria for large enough domains $(-L,L)$, with boundary conditions $\mathcal{B}_\mathrm{l/r}$ at $\pm L$. . 

We employ a two-sided shooting approach. We start with transverse boundary layers $U_\mathrm{l}(x;k,\varphi), x>0$ and $U_\mathrm{r}(x;k,\varphi),x<0$, which satisfy the boundary conditions at $\mathcal{B}_\mathrm{l/r}$ at $x=0$, and 
\begin{align*}
U_\mathrm{l}(x;k,\varphi)-U_\mathrm{st}(kx-\varphi;k)\to & 0,\qquad x\to+\infty,\\
U_\mathrm{r}(x;k,\varphi)-U_\mathrm{st}(kx-\varphi;k)\to & 0,\qquad x\to-\infty,
\end{align*}
with associated displacement-strain relations
\[
d_\mathrm{l/r}(k,\varphi)=0.
\]
Shifting $U_\mathrm{l/r}$ to the boundaries, we look for solutions $U$ that are close to $U_\mathrm{l}(x+L;k_\mathrm{l},\varphi_\mathrm{l})$ on $(-L,0)$ and 
$U_\mathrm{r}(x-L;k_\mathrm{r},\varphi_\mathrm{r})$ on $(0,L)$. 
Continuity of the solution at $x=0$ up to exponentially small terms gives \emph{phase and wavenumber matching},
\begin{equation}\label{e:gl}
k_\mathrm{l}L-\varphi_\mathrm{l}\equiv -
k_\mathrm{r}L-\varphi_\mathrm{r} \,
\mathrm{mod}\,2\pi,\qquad k_\mathrm{l}=k_\mathrm{r}.
\end{equation}
The following proposition shows that this formal phase and wavenumber matching procedure actually determines solutions.

\begin{Proposition}\label{p:phasematching}
Assume the existence of families of transverse boundary layers  $U_\mathrm{l}(x;k,\varphi), x>0$ and $U_\mathrm{r}(x;k,\varphi),x<0$, to boundary conditions $\mathcal{B}_\mathrm{l/r}$ at $x=0$, with displacement-strain relations $d_\mathrm{l/r}(k,\varphi)=0$. Then 
there exist solutions on $(-L,L)$ with boundary conditions $\mathcal{B}_\mathrm{l/r}$ at $x=\pm L$ which are  $L^\infty$-close to $U_\mathrm{r}(x-L;k_\mathrm{l},\varphi_\mathrm{l})$ on $x<0$ and  to  $U_\mathrm{l}(x+L;k_\mathrm{r},\varphi_\mathrm{r})$ on $x>0$ if exponentially corrected phase and wavenumber matching conditions hold,
\begin{align}
k_\mathrm{l}&=k_\mathrm{r}+\mathcal{R}_k\label{e:pm1}\\
(k_\mathrm{l}+k_\mathrm{r})L&\equiv\varphi_\mathrm{l} - \varphi_\mathrm{r}+\mathcal{R}_\varphi\mathrm{mod}\,2\pi,\label{e:pm6}
\end{align}
where $\mathcal{R}_j(k_\mathrm{l},k_\mathrm{r},\varphi_\mathrm{l},\varphi_\mathrm{r},L)$, $j=\varphi,k$, are smooth and $
\mathcal{R}_j=\rmO(\rme^{-\eta L}),
$
together with its derivatives. 
\end{Proposition}

\begin{Proof}
We write $U=U_\mathrm{l}(x+L)\chi_-(x)+U_\mathrm{r}(x-L)\chi_+(x)+W(x)$, where $\chi_\pm$ is a smooth partition of unity, $\chi_-+\chi_+=1$, $\chi_\pm(x)=1$ for $\pm x>1$. Exploiting the fact that the boundary layers solve the ODE, we find the equation
\begin{equation}\label{e:main}
\dot{W}=F(U_\mathrm{l}(x+L)\chi_-(x)+U_\mathrm{r}(x-L)\chi_+(x)+W(x))-
F(U_\mathrm{l}(x+L)\chi_-(x)+U_\mathrm{r}(x-L)\chi_+(x))+R(x),
\end{equation}
with boundary conditions
\begin{equation}\label{e:bc1}
W(\pm L)\in \mathcal{B}_\mathrm{l/r}-U_\mathrm{l/r}(0),
\end{equation}
where the remainder is 
\begin{align*}
R(x)=&F(U_\mathrm{l}(x+L)\chi_-(x)+U_\mathrm{r}(x-L)\chi_+(x))-F(U_\mathrm{l}(x+L))\chi_-(x)-F(U_\mathrm{r}(x-L))\chi_+(x)\\
&-U_\mathrm{l}(x+L)\chi_-'(x)-U_\mathrm{r}(x-L)\chi_+'(x).
\end{align*}
We view the system (\ref{e:main}) as an equation with variables $W$ and $\tau_\pm$, where $\tau_\pm$ enter as a parameterization of the boundary layers 
\[
U_\mathrm{l/r}(x)=U_\mathrm{l/r}(x;k_\mathrm{l/r}(\tau_\pm),\varphi_\mathrm{l/r}(\tau_\pm)),\quad 
d_\mathrm{l/r}(k_\mathrm{l/r}(\tau_\pm),\varphi_\mathrm{l/r}(\tau\varphi_\mathrm{l/r}(\tau_+)))=0,\quad |k_\mathrm{l/r}'(\tau_\pm)|^2+|\varphi_\mathrm{l/r}'(\tau_\pm)|^2=1.
\]

Note that the remainder is small when the system (\ref{e:pm1})--(\ref{e:pm6}) is satisfied. We therefore need to understand the linearized operator
\[
\dot{W}-F'(U_\mathrm{l}(x+L)\chi_-(x)+U_\mathrm{r}(x-L)\chi_+(x))W,\quad W(\pm L)\in T\mathcal{B}_\mathrm{l/r}|_{U_\mathrm{l/r}(0)}.
\]
Note that this operator is close to the linearization at individual boundary layers on $(-L,0)$ and $(0,L)$, respectively. This linearization at individual boundary layers is not invertible in the limit $L=\infty$. It is however Fredholm in exponentially weighted spaces, with index -1. We therefore study relaxed boundary conditions
\begin{equation}\label{e:bc2}
W(\pm L)+\alpha_\mathrm{l/r} \psi_\mathrm{l/r}(0)\in \mathcal{B}_\mathrm{l/r}-U_\mathrm{l/r}(0),
\end{equation}
replacing (\ref{e:bc1}). 
Here, $\alpha_\mathrm{l/r}$ are independent scalar variables and $\psi_\mathrm{l/r}(x)$ are solutions to 
\[
\dot{\psi}_\mathrm{l/r}-F'(U_\mathrm{l/r}(x))\psi_\mathrm{l/r}=0,
\]
such that $\psi_\mathrm{l/r}\not\in T\mathcal{B}_\mathrm{l/r}$, with at most linear growth at $\pm\infty$, respectively. 

We will show that the associated linearized operator, $\mathcal{L}_{[-L,L]}$,
\[
\dot{W}-F'(U_\mathrm{l}(x+L)\chi_-(x)+U_\mathrm{r}(x-L)\chi_+(x))W,\quad W(\pm L)+\alpha_\mathrm{l/r} \psi_\mathrm{l/r}(0)\in T\mathcal{B}_\mathrm{l/r}|_{U_\mathrm{l/r}(0)},
\]
is invertible with $L$-uniform bounds in exponentially weighted spaces,
\begin{equation}\label{e:ew}
\|W\|_{L^2_\eta}=\|\min\{\rme^{\eta(x+L)},\rme^{\eta(-x+L)}\} W(x)\|_{L^2},
\end{equation}
with domain $H^1_\eta$, where $W,W'\in L^2_\eta$. 
%

For this, we solve $\mathcal{L}_{[-L,L]}W=H$ by decomposing $H=H_-+H_+$, $H_\pm=\chi_\pm H$,
\begin{align*}
\dot{W_-}-F'(U_\mathrm{l})W_- - 
[F'(U_\mathrm{l}\chi_-+U_\mathrm{r}\chi_+)-F'(U_\mathrm{l})]W_+=H_-,\\
\dot{W_+}-F'(U_\mathrm{r})W_+ - 
[F'(U_\mathrm{l}\chi_-+U_\mathrm{r}\chi_+)-F'(U_\mathrm{r})]W_-=H_+,\\
W_-(-L)+W_+(-L)+\alpha_-\psi_\mathrm{l}\in T\mathcal{B}_\mathrm{l},\\
W_+(L)+W_-(L)+\alpha_+\psi_\mathrm{r}\in T\mathcal{B}_\mathrm{r}.
\end{align*}
In the exponentially weighted norms (\ref{e:ew}), the boundary  terms $W_+(-L)$ and $W_-(L)$, as well as the coupling terms involving $W_+$ in the first equation and $W_-$ in the second equation are exponentially small, so that this system is a small perturbation of 
\begin{align*}
\dot{W_-}-F'(U_\mathrm{l})W_- =H_-,\\
\dot{W_+}-F'(U_\mathrm{r})W_+ =H_+,\\
W_-(-L)+\alpha_-\psi_\mathrm{l}\in T\mathcal{B}_\mathrm{l},\\
W_+(L)+\alpha_+\psi_\mathrm{r}\in T\mathcal{B}_\mathrm{r}.
\end{align*}
We can extend $H_\pm$ to half lines and thereby find $W_\pm$ by inverting the linearization at the boundary layers. 

Summarizing, we have shown uniform bounded invertibility of the linearization of (\ref{e:main}), with extended boundary conditions (\ref{e:bc2}), variables $W,\alpha_\pm$, and parameters $\tau_\pm,L$. The residual is small, $R=\rmO(\rme^{-\eta L})$ when the phase matching conditions,
\[
k_\mathrm{l}(\tau_-)=k_\mathrm{r}(\tau_+),\qquad  2k_\mathrm{r}L=\varphi_\mathrm{l}(\tau_-)-\varphi_\mathrm{r}(\tau_+),
\]
are satisfied. We therefore obtain families of solutions $(W,\alpha_\pm)=(W^*,\alpha_\pm^*)(\tau_\pm,L)$. It remains to solve 
\[
\alpha_+^*(\tau_+,\tau_-,L)=0,\qquad 
\alpha_-^*(\tau_+,\tau_-,L)=0.
\]
A straightforward but tedious expansion of this system shows that it is equivalent to wavenumber and phase matching, up to exponentially small corrections,
\[
k_\mathrm{l}=k_\mathrm{r}+\rmO(\rme^{-\delta L}),\quad k_\mathrm{l}L-\varphi_\mathrm{l}=k_\mathrm{r}L-\varphi_\mathrm{r}+\rmO(\rme^{-\delta L}) \mod 2\pi,
\]
where 
\[
(k_\mathrm{l},\varphi_\mathrm{l})=(k_\mathrm{l}(\tau_-),\varphi_\mathrm{l}(\tau_-))+\alpha_-(k_\mathrm{l}'(\tau_-),\varphi_\mathrm{l}'(\tau_-))^\perp.
\]
This establishes the desired phase and wavenumber matching conditions. 

%
%
%
%
%
%
%
%
%
%

\end{Proof}

\begin{Remark}
A more geometric proof would be based on flowing the boundary conditions from $x=\pm L$ to $x=0$ and exploiting inclination lemmas (see for instance 
\cite[Thm 3.1]{exchange}) to see that the intersection is transverse in a complement of the tangent space to the two-dimensional family of periodic orbits, so that intersections can be computed from matching within this manifold up to exponentially small terms.
Our approach here is closer in spirit to \cite{rade}.
\end{Remark}

\subsection{Subtraction and quantization --- from displacement-strain curves to equilibria}\label{s:2.6}

In the following, we discuss conceptually the solutions to  (\ref{e:pm1})--(\ref{e:pm6}). We will neglect exponentially small terms, which can be easily accommodated in the analysis. 
As in the proof of Proposition \ref{p:phasematching}, we parameterize solutions to the displacement-strain relation as curves via 
\[
\Gamma_\mathrm{l/r}=\{(k_\mathrm{l/r}(\tau_\pm),\varphi_\mathrm{l/r}(\tau_\pm),
\]
with nonvanishing tangent vectors $|(\varphi',k')|^2=1$ and study wavenumber and phase-matching, (\ref{e:pm1})--(\ref{e:pm6}) as two equations in three variables  $\tau_\pm,L$. It is conceptually helpful to add an intermediate step, solving 
\begin{align*}
0&=k_\mathrm{l}(\tau_-)-k,\\
0&=k_\mathrm{r}(\tau_+)-k,\\
0&=\varphi_\mathrm{l}(\tau_-)-\varphi_\mathrm{r}(\tau_+)-\psi.
\end{align*}
Note that with the definition of $k,\psi$ from this system, phase matching reduces to $\psi=2kL\mathrm{mod}\, 2\pi$. 

The linearization of this auxiliary system is readily seen to be onto as long as  $k'_\mathrm{l}$ and $k'_\mathrm{r}$ do not vanish simultaneously. Since in this case the linearization with respect to $k,\varphi$, and either $\tau_+$ or $\tau_-$ is invertible, we can parameterize the solution set as curves $(k,\psi)(\tau)$. Critical points of $k(\tau)$ are given by the union of critical points of $k_\mathrm{l/r}(\tau_\pm)$. Critical points of $\psi$ correspond to values of $k$ where tangent vectors of left and right displacement-strain relations are collinear, 
\[
\psi'=0 \quad \Longleftrightarrow \left|\begin{array}{cc}
k_\mathrm{l}'&k_\mathrm{r}'\\
\psi_\mathrm{l}'&\psi_\mathrm{r}'
\end{array}
\right|=0.
\]
We can formalize these operations as follows. 
\begin{Definition}[Differential displacement-strain]
We refer to the curves $(k,\psi)(\tau)$ as \emph{differential displacement-strain curves}, obtained by taking the difference between left and right displacement-strain curve, pointwise in $k$. 
\end{Definition}

From differential displacement-strain curves, we obtain equilibria in the bounded domain for large $L$ by intersecting with the straight line $\psi=2kL\mathrm{mod}\, 2\pi$, effectively quantizing the curves. 

In summary, the operations that lead to a description of equilibria in bounded domain are:

\begin{itemize}
\item \emph{vertical subtraction}, which gives differential displacement-strain curves;
\item \emph{quantization} of differential displacement-strain, by intersecting with $\psi=2kL$.
\end{itemize}

\begin{figure}
\centering{\includegraphics[width=0.95\textwidth]{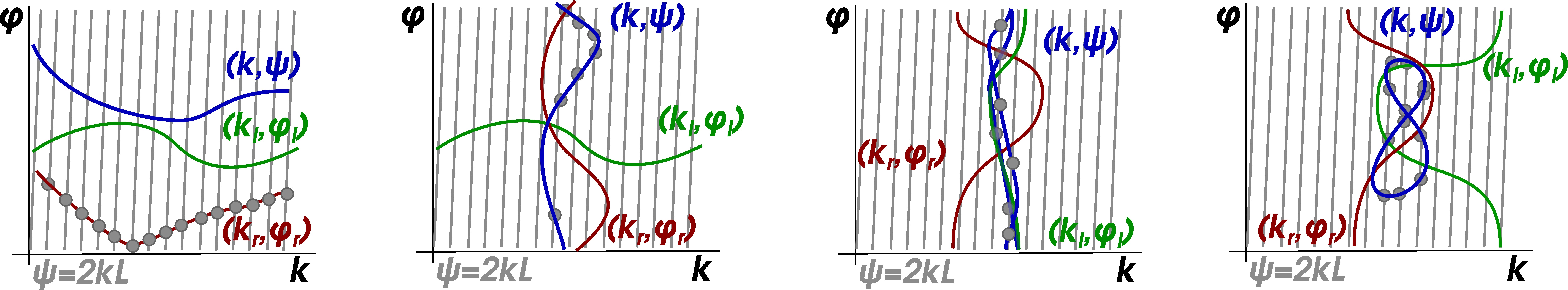}}
\caption{The figures show schematic plots of left and right displacement-strain relations, their differential displacement-strain curves, and the quantized intersections with $\varphi=2kL$. From left to right: two-sided phase selection, phase and wavenumber selection, two-sided wavenumber selection with $J_\mathrm{l}\subset J_\mathrm{r}$ (note the two branches of the differential displacement-strain curve), and two-sided wavenumber selection with marginal overlap. 
}\label{f:dds}
\end{figure}

\begin{figure}
\centering{\includegraphics[width=0.95\textwidth]{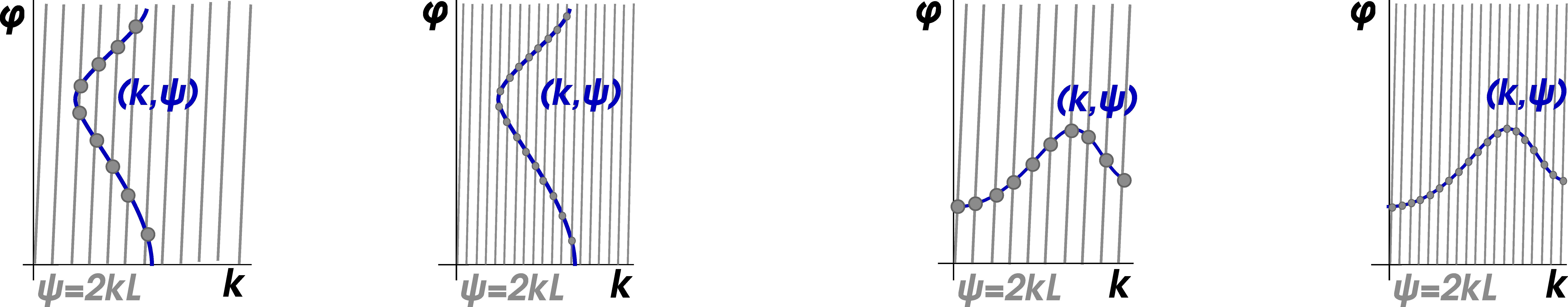}}
\caption{The figure illustrates the effect of increasing the domain size on equilibria. Wavenumber selection with nonzero winding number (left) leads to a continuous (snaking) family of equilibria as $L$ is increased. Phase selection (right) leads to chains of creation and annihilation events of equilibria near the Eckhaus boundary. }\label{f:dds2}
\end{figure}

\subsection{From displacement-strain to equilibria: conceptual examples}\label{s:2.7}
In the following, we analyze several representative cases in some more detail; see Figure \ref{f:dds} for illustrations of differential displacement-strain curves and Figure \ref{f:dds2} for the effect of varying $L$ in the case of wavenumber and phase selection by the differential displacement-strain curves, respectively. 

\paragraph{Two-sided phase selection.}
In this case, both displacement-strain curves are defined for $k_\mathrm{r/l}\in J$, and so is the differential displacement-strain curve. In the simplest case, $\psi=\psi(k)$ is a graph, and intersections with the straight lines $2kL$ are monotone in $L$, increasing from $k_-$ to $k_+$ as $L$ increases. In particular, specific equilibria appear and disappear as $L$ is increased by an order-one amount because of instabilities at the boundary of $J$. 

\paragraph{Phase and wavenumber selection.}

Suppose that $k_\mathrm{r}\in J$ and $k_\mathrm{l}\in J_\mathrm{l}$. The differential displacement-strain curve is now defined on $J_\mathrm{l}$, as a closed curve, with the same winding number of the $\psi$ component as $\varphi_\mathrm{r}$. Intersections with the straight lines $2kL$ give continuous families of equilibria in the case of nonzero winding number (the right boundary layer is only selecting wavenumbers, allowing for arbitrary phases). 

\paragraph{Two-sided wavenumber selection.}

When $J_\mathrm{r}\subset J_\mathrm{l}$, the situation is equivalent to left-sided wavenumber selection. New phenomena appear when $J_\mathrm{r}$ and $J_\mathrm{l}$ are not contained one in the other. Of course, no equilibria exist when these intervals are mutually disjoint. In the simplest case, a mutual intersection gives figure-eight type differential displacement-strain curves. 
Again, one expects phase slips as the size of the domain is gradually increased. 

%
%
%
%
%
%

\paragraph{Summary.}

Envisioning slowly increasing the domain size $L$, we expect continuous dependence of the solution on $L$  only when one of the boundary conditions is wavenumber selecting with nonzero winding number in the phase, and when all wavenumbers of this boundary condition are compatible with the other boundary condition. Other combinations will typically lead to phase slips. 

Incompatible boundary conditions (selecting disjoint wavenumber intervals) will lead to phase drift as we shall see in the simple example of the phase-diffusion problem, below.


%
%
%
%
%

\section{Boundary layers in phase-diffusion}\label{s:2a}

Using a multiple-scale expansion, one can derive simplified dynamics near roll solutions, known as the Cross-Newell phase diffusion equation. One expands an Ansatz $u_\mathrm{st}(\theta(t,x);\theta_x(t,x))$ assuming long-wavelength modulations of the phase $\theta$, and finds the nonlinear diffusion equation
\[
\theta_t=(b(\theta_x))_x,
\]
where $b'(k)=\lambda''(0)$. In particular, the long-wavelength problem is well-posed when $b'>0$. Boundary conditions are given  by smooth curves $\mathcal{B}$ in the $(\theta,\theta_x)$-plane. Each point $(\theta,\theta_x)=(\varphi(s),k(s))\in\mathcal{B}$, gives rise to a unique equilibrium in $x>0$, $\varphi+kx$. Linearizing the equation and boundary conditions gives
\[
\theta_t=b'(k)\theta_{xx}, \quad k'\theta-\varphi'\theta_x=0.
\]
Patterns with $k'\varphi'<0$ are unstable with eigenfunction $\rme^{-\gamma x}$, $\gamma=k'/\varphi'$, and eigenvalue $\lambda=\gamma^2$.  Stable patterns possess $k'\varphi'>0$, with resonance pole $\lambda=\gamma^2$.  Note that when $\varphi'\to 0$, an eigenvalue (or resonance pole) disappears at infinity, not unsurprisingly due to the singular perturbation in the boundary conditions when passing from Dirichlet to Robin (or mixed).  
The allowed interval $J$ is given through $b'(J)>0$. Clearly, Dirichlet boundary conditions $\theta=\varphi_0$ select a phase and Neumann boundary conditions $\theta_x=k_0$ select a wavenumber. In reaction-diffusion systems or the Swift-Hohenberg equation, phase selection is accomplished by Neumann boundary conditions. Since perturbations from Neumann to mixed are non-singular, we can exclude eigenvalues near infinity and changes of stability when $\varphi'=0$ for such more general systems. 

Within the class of phase-diffusion problems, one would require $\varphi'$ to be nonzero to guarantee well-posedness. As a consequence, one finds $k=k(\varphi)$ and isolas are excluded, that is, wavenumber selection implies a non-zero winding number of the displacement-strain curve. 

Imposing different wavenumber selection conditions $\theta_x=k_\pm$ at $x=\pm L$, one finds drifting patterns,
\[
\theta(t,x)=\int_0^x b^{-1}\left(\frac{k_-+k_+}{2}-\frac{k_--k_+}{2L}y\right)\rmd y-
\frac{k_--k_+}{2L}t,
\]
which are effectively time-periodic solutions to the equation since $\theta\in\R/2\pi\Z$. Roughly speaking, $b(\theta_x)$ interpolates the selected wavenumbers linearly, leading to a nonlinear phase profile and constant phase drift. Note that the expression is well-defined as long as the interval of wavenumbers defined by the boundary conditions lies within the stable regime, $b'>0$. 

On the other hand, our analysis validates the phase-diffusion approximation in semi-bounded (or even bounded) domains for general striped phases in a very concise fashion: boundary conditions should be replaced by displacement-strain relations, 
\[
d(\theta,\theta_x)=0,  
\]
at $x=0$, respectively, where $\theta$ is understood modulo $2\pi$, and the equation in the interior of the domain should be simply $\theta_{xx}=0$. The fact that displacement-strain relations may not lead to well-posed boundary conditions points to the fact that some defects (here, boundary conditions) can simply not be incorporated into a phase-modulation description. On the other hand, one can expand the effective boundary conditions near $\theta_x=k_*$, and solutions near $\theta=k_*x+\varphi$, to obtain locally valid effective Robin boundary conditions.

%
%

\section{Boundary layers in the Ginzburg-Landau equation}\label{s:3}

The Ginzburg-Landau equation
\[
A_t=A_{xx}+A-A|A|^2. \quad A\in\C,
\]
arises as a modulation equation at the onset of a Turing instability and as such provides a universal model for small-amplitude Turing patterns. We will study the effect of boundary conditions within this equation but will also comment on the effect of scalings used to derive Ginzburg-Landau on the boundary conditions. 
Stationary solutions satisfy 
\begin{equation}\label{e:stgl}
A_x=B,\quad B_x=-A+A|A|^2. 
\end{equation}
Exact solutions include the periodic patterns $A(x)=\sqrt{1-k^2}\rme^{\rmi k x}$, and the defect solutions, 
\begin{equation}\label{e:tanh}
A_\mathrm{d}(x;k)=\left(\sqrt{2}k+\rmi \sqrt{1-3k^2}\tanh(\sqrt{1-3k^2}x/\sqrt{2})\right)\rme^{\rmi k x}.
\end{equation}
One readily calculates that the argument of $A(x;k)\rme^{-\rmi k x}$ increases by $2\arctan(\sqrt{1-3k^2}/(\sqrt{2}k))$ when $k>0$ and decreases by the same amount when $k<0$. For $k=0$, the argument is not defined at $x=0$ and the final phase shift is $\pi$. The defects are invariant under the reflection $A(x)\mapsto \bar{A}(-x)$. Defects are known to be unstable for energetic reasons. 

As a Hamiltonian system with two degrees of freedom and a phase rotation symmetry, (\ref{e:stgl}) is integrable and all solutions are explicit, in particular the stable and unstable manifolds of periodic orbits. 

Indeed, the bounded part of stable and unstable manifolds is given by the defect, and the unbounded parts are given by
\begin{equation}\label{e:coth}
A_\mathrm{}(x;k)=\left(\sqrt{2}k+\rmi \sqrt{1-3k^2}\coth(\sqrt{1-3k^2}x/\sqrt{2})\right)\rme^{\rmi k x},
\end{equation}
with $x<0$ or $x>0$, respectively.

The Ginzburg-Landau equation arises universally as an amplitude equation for small striped patterns, and therefore is a prototype for our study. On the other hand, the additional symmetries in this equation sometimes distort phenomena. In fact, one can further ``reduce'' the Ginzburg-Landau equation by assuming slowly varying phase near periodic patterns, and then derive the phase-diffusion equation, discussed in the preceding section. Our main findings here do replicate the phenomena for the phase-diffusion equation. However, we find several new phenomena, in particular:
\begin{itemize}
\item displacement-strain curves with winding number 2 (inhomogeneous Neumann);
\item displacement-strain curves where $\varphi$ is \emph{not} monotone (inhomogeneous Dirichlet);
\item defect pinch-off at terminal points of displacement-strain curves (inhomogeneous Neumann and Dirichlet);
\item reconnection crises near the Eckhaus boundaries (inhomogeneous Neumann and Dirichlet).
\end{itemize}

\subsection{Symmetries and invariant coordinates}\label{s:3.1}
While one can find boundary layers using the explicit representation of stable manifolds in (\ref{e:tanh}) and (\ref{e:coth}), we found it significantly more manageable (albeit, presumably, equivalent) to use different coordinates and reduce to polynomial equations. 

The stationary Ginzburg-Landau equation (\ref{e:stgl}) is invariant under the gauge and reversal symmetries
\[
\mathcal{R}_\psi (A,B)=\rme^{\rmi\psi|}(A,B),\qquad \mathcal{S}(A,B)=(\bar{A},-\bar{B}),\qquad \mathcal{T}(A,B)=(A,-B).
\]
Typically, $\mathcal{S}$ corresponds to the reflection symmetry in the original equation, when Ginzburg-Landau is derived as an amplitude equation, since, say in the Swift-Hohenberg case, $u\sim A\rme^{\rmi x}+c.c.$. Complex conjugation symmetry $\mathcal{ST}$ is a ``normal form artifact. 
The system is also Hamiltonian, with $\mathcal{H}=|B|^2+|A|^2-\frac{1}{2}|A|^4$, so that 
\[
A_x=\partial_{\bar{B}}\mathcal{H},\qquad B_x=-\partial_{\bar{A}}\mathcal{H}.
\]
One can factor the gauge symmetry and consider the equation in new variables for $A,B\in\C^2$, the Hilbert invariants of the circle action $\mathcal{R}_\psi$,
\[
\alpha=|A|^2\geq 0,\quad \beta=|B|^2\geq 0, \quad \mathcal{M}=\frac{1}{2}(m+\rmi M)=A\bar{B}\in\C,\quad \mbox{with relation }\  4\alpha\beta=\mathcal{M}\bar{\mathcal{M}}.
\]
The inverse transformation gives $A$ and $B$ up to a relative phase, 
\begin{equation}\label{e:inv}
A=\sqrt{\alpha}\rme^{\rmi\varphi_1},\quad B=\sqrt{\beta}\rme^{\rmi\varphi_2}, \quad \rme^{\rmi(\varphi_1-\varphi_2)}=\frac{M}{\sqrt{4\alpha\beta}}.
\end{equation}
The equations read 
\begin{align}
\alpha_x&=m\nonumber\\
\beta_x&=m(\alpha-1)\nonumber\\
m_x&=2(\beta+\alpha^2-\alpha)\nonumber\\
M_x&=0.
\end{align}
Since $\beta$ is given implicitly through the relation and $M$ is conserved, we end up with a two-dimensional system for $(\alpha,m)$, in which we can eliminate $\beta$, most conveniently using the Hamiltonian, which gives
\begin{align}
\alpha_x&=m\nonumber\\
m_x&=2\left(\mathcal{H}-2\alpha+\frac{3}{2}\alpha^2\right).\label{e:2d}
\end{align}
From solutions to this planar ODE, for any fixed $\mathcal{H}$, we find solutions to the original system by first reconstructing $\beta$ from $\mathcal{H}$, reconstructing $M^2=\alpha\beta - m^2$ from the relation, and using the inverse transformation (\ref{e:inv}). Note however that each solution of (\ref{e:2d}) yields two solutions corresponding to the different signs of $M$.

\subsection{Stable manifolds}\label{s:3.2}
Periodic patterns are of the form $A=\sqrt{1-k^2}\rme^{\rmi k x}$, which in invariant coordinates gives
\[
\alpha=1-k^2,\quad \beta=k^2(1-k^2),\quad m=0,\quad M=-2k(1-k^2),\quad \mathcal{H}=\frac{1}{2}(1-k^2)(1+3k^2).
\]
Note that $\mathcal{H}$ is a monotone function of $k^2$ within the Eckhaus boundary $k^2\in[0,\frac{1}{3})$. We can also use the expressions for $\mathcal{H}$ and $M$ to express (\ref{e:2d}) using the parameter $k$,
\begin{align}
\alpha_x&=m\nonumber\\
m_x&=2\left(\mathcal{H}(k)-2\alpha+\frac{3}{2}\alpha^2\right), \qquad\mathcal{H}(k)= \frac{1}{2}(1-k^2)(1+3k^2).\label{e:2dk}
\end{align}
Equilibria are $\alpha_\pm=\frac{2}{3}(1\pm\sqrt{1-\frac{3}{2}\mathcal{H}})$, or, 
\[
\alpha_+=1-k^2,\quad \alpha_-=\frac{1}{3}(1+3k^2).
\]
with $\alpha_+$ corresponding to the Eckhaus-stable equilibrium for $k^2<1/3$.
We normalize the Hamiltonian of (\ref{e:2dk}),
\[
h(\alpha,m;k)=\frac{1}{2}m^2-V(\alpha;k),\quad V(\alpha;k)=\int_{1-k^2}^\alpha\left((1-k^2)(1+3k^2)-4a+3a^2\right)\rmd a=(\alpha-\alpha_*)(\alpha-\alpha_+)^2,
\]
where $\alpha_*=2k^2$. The stable manifold $W^\mathrm{s}$ of $\alpha_+$ consists of an unbounded branch with 
\[
m=-\sqrt{2(2\mathcal{H}\alpha-2\alpha^2+\alpha^3+h_+)}, \qquad \alpha>\alpha_+,
\]
and a bounded branch, 
\[
m=\pm\sqrt{2(2\mathcal{H}\alpha-2\alpha^2+\alpha^3+h_+)}, \qquad \alpha_*<\alpha<\alpha_+.
\]
A somewhat more compact representation is achieved in coordinates $\tau=\frac{m}{\alpha-\alpha_+}$, 
where $W^\mathrm{s}$ is given through
\begin{equation}\label{e:ws}
\alpha=2k^2+\frac{1}{2}\tau^2,\quad m=\tau(\frac{1}{2}\tau^2-(1-3k^2)),\quad \tau\leq \tau_*=\sqrt{2(1-3k^2)}.
\end{equation}
In summary, given a wavenumber $|k|<1/\sqrt{3}$, we can parameterize its stable manifold by $\tau$ via (\ref{e:ws}).

\begin{figure}
\includegraphics[width=0.24\textwidth]{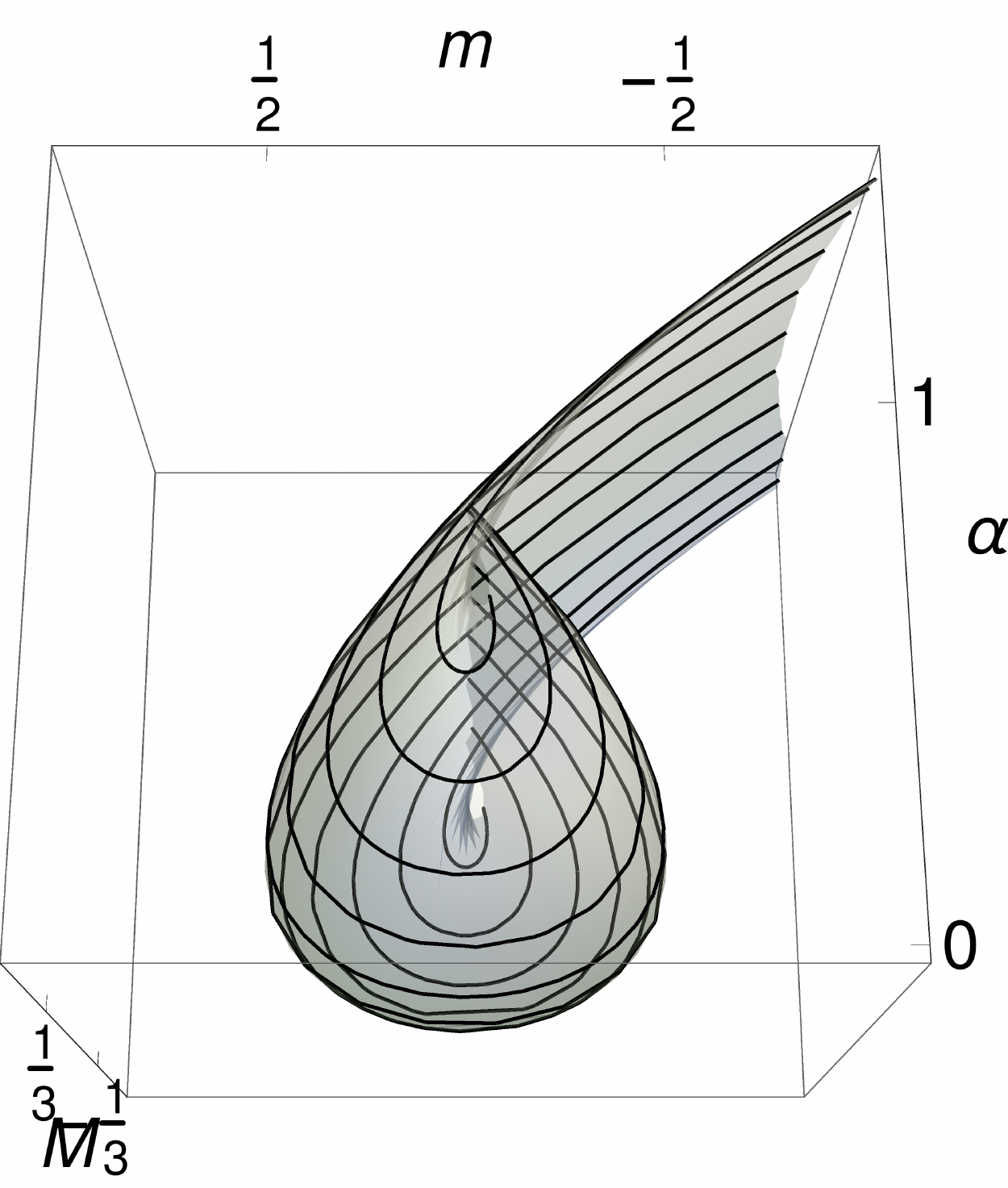}\hfill 
\includegraphics[width=0.24\textwidth]{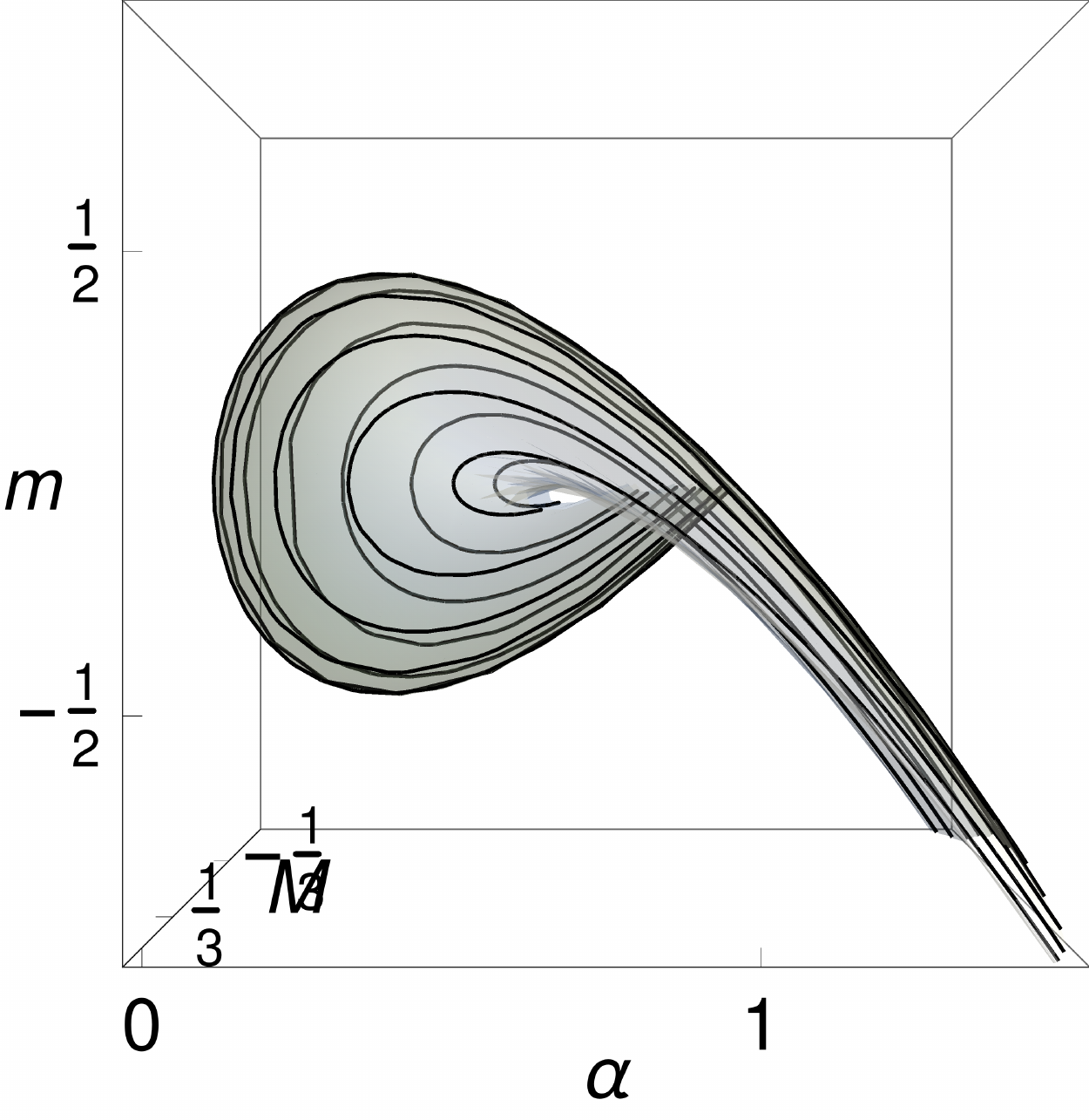}\hfill 
\includegraphics[width=0.24\textwidth]{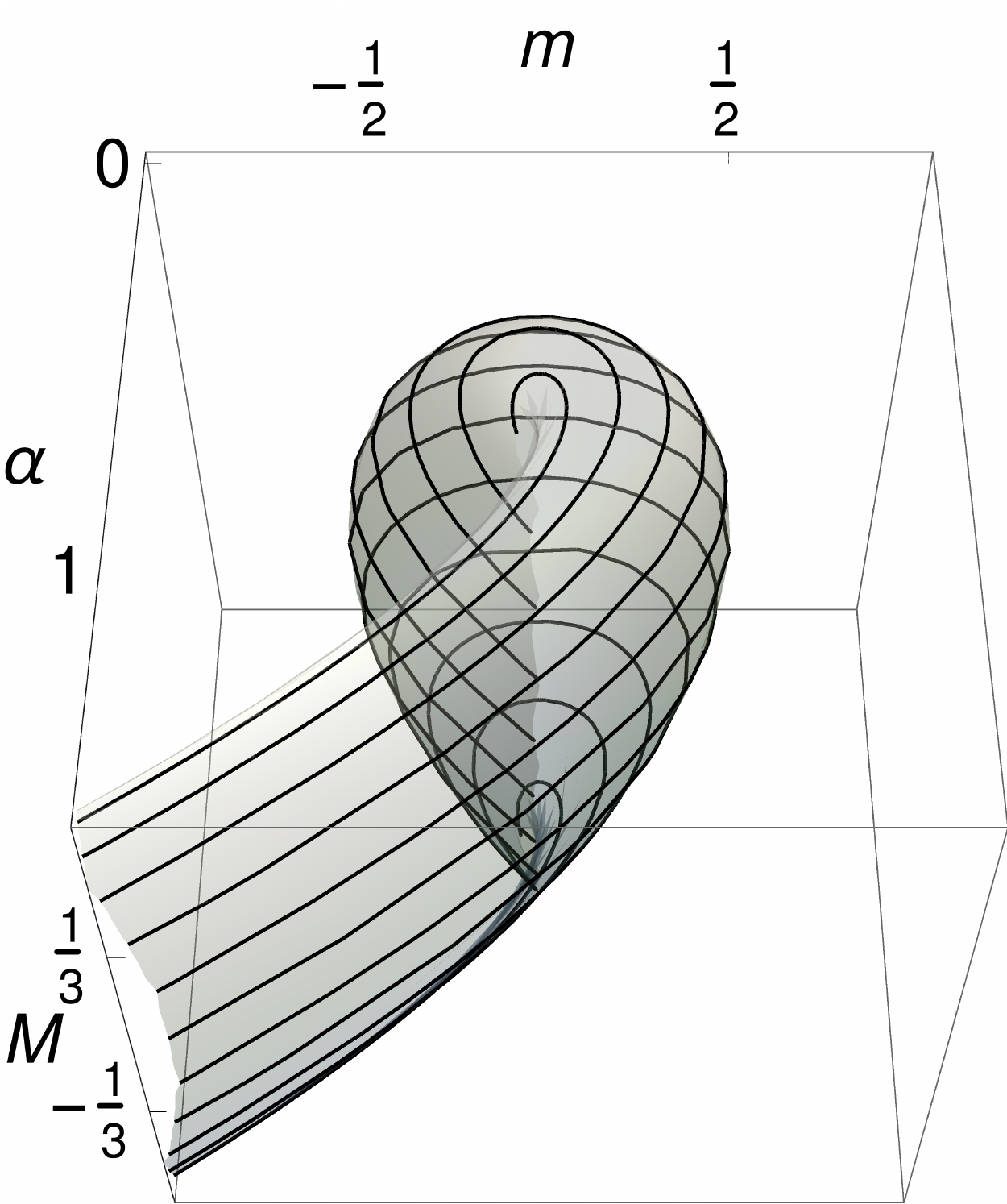}\hfill 
\includegraphics[width=0.24\textwidth]{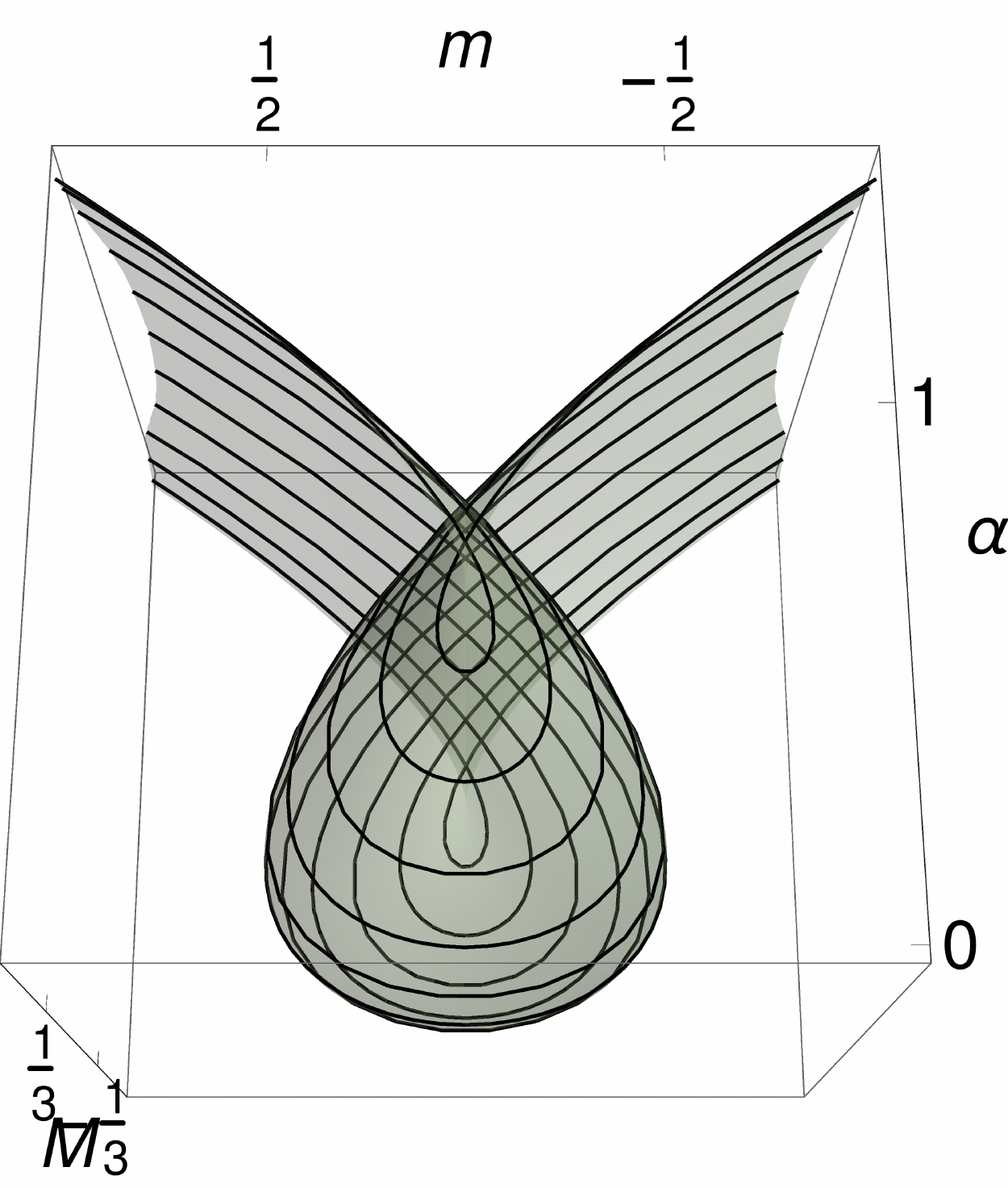}
\caption{The stable manifold of wave trains, plotted as $(\alpha,m,M)(k,\tau)$, with $|k|\leq 1/\sqrt{3}$ (left 3 figures). Stable \emph{and} unstable manifold  (right).}\label{f:fish}
\end{figure}

Figure \ref{f:fish} shows the stable manifold of the wave trains in $(\alpha,m,M)$-space.

\subsection{Phase shifts}\label{s:3.2a}

Phase shifts $\varphi$ are defined so that for a boundary layer
\[
\lim_{x\to\infty}\Im\left(\log(A(x))-(kx-\varphi)\right)=0,
\]
which gives, using the expression for the derivative of $A$, 
\begin{equation}
\varphi=-\lim_{x\to\infty}\left(\Im(\log A(x))-kx\right)=-\Im\log(A(0))-\Delta\varphi,\qquad \Delta\varphi=-\int_0^\infty\left(k+\frac{M}{2\alpha(x)}\right)\rmd x.\label{e:phaseshift0}
\end{equation}
Using (\ref{e:ws}) gives $\tau\tau'=m$ and $\tau'=\frac{1}{2}(\tau-\tau_*)(\tau+\tau_*)$ with $\tau_*=\sqrt{2(1-3k^2)}$. We can use $\tau$ as integration variable with $\tau(\infty)=-\tau_*$, to obtain 
\[
\Delta\varphi=-\int_{\tau(0)}^{-\tau_*}\frac{2k}{4k^2+\tau^2}\rmd\tau=\arctan\left(\frac{\sqrt{2(1-3k^2)}}{2k}\right)+\arctan\left(\frac{\tau(0)}{2k}\right),
\]
and, altogether,
\begin{equation}\label{e:phaseshift}
\varphi=\arctan\left(\frac{\sqrt{2(1-3k^2)}}{2k}\right)-\Im\log(A(0))+\arctan\left(\frac{\tau(0)}{2k}\right).
\end{equation}

\subsection{Computing displacement-strain curves for specific boundary conditions}\label{s:3.3}

General boundary conditions will not be invariant under the gauge symmetry. Consider a general affine boundary condition of the form 
\begin{equation}\label{e:bcpl}
B=k_1A+k_2\bar{A}+\mu,\quad k_1,k_2,\mu\in\C;
\end{equation}
note that this comprises all affine boundary conditions which are not strictly Dirichlet in either component. We will discuss Dirichlet boundary conditions at the end of this section. 

The group orbit of the two-dimensional plane defined by (\ref{e:bcpl}) under the gauge symmetry is
\[
\{(A,B);B=\rme^{\rmi\phi}\left(k_1\xi
+k_2\bar{\xi}
+ \mu\right),A=\rme^{\rmi\phi}\xi,\xi\in\C,\phi\in[0,2\pi)\},
\]
which gives the parameterized representation  in invariants 
\begin{equation}\label{e:bcinv}
\alpha=|\xi|^2, \quad \mathcal{M}=2\xi(k_1\xi
+k_2\bar{\xi}
+ \mu).
\end{equation}
One now proceeds to find boundary layers by solving (\ref{e:bcinv}) together with (\ref{e:ws}) and the relation for $M$,
\begin{equation}\label{e:bcws}
\alpha=2k^2+\frac{1}{2}\tau^2,\quad m=\tau(\frac{1}{2}\tau^2-(1-3k^2)),\quad M=-2k(1-k^2),\quad \tau\leq \sqrt{2(1-3k^2)}.
\end{equation}
Equations (\ref{e:bcinv}) and (\ref{e:bcws}) can be viewed as a polynomial system of 6 real equations in 7 real variables $\Re\xi,\Im\xi, \alpha, m,M,k,\tau$, together with an inequality constraint. We will next show how to solve this system in several simple cases. 

\paragraph{Neumann boundary conditions}

For $k_1=k_2=0$, $\mu> 0$, we can use the second equation in (\ref{e:bcinv}) and the second and third equation in (\ref{e:bcws}) to express $\xi$ in terms of $\tau,k$, which we can then substitute into the first equations of  (\ref{e:bcinv}) and (\ref{e:bcws}) to obtain an equation in $\tau,k$, only. After eliminating a factor $4k^2+\tau^2$, which only gives the trivial solution $A\equiv 0$, we find
\begin{equation}\label{e:ktn}
4k^4-8\mu^2+(\tau^2-2)^2+8k^2(\tau^2-1)=0,
\end{equation}
combined of course with the inequality $\tau\leq \sqrt{2(1-3k^2)}$. This equation is a simple quadratic in $\kappa=k^2$ and $\sigma=\tau^2$, 
\[
\left(\sigma+(4-2\sqrt{3})\kappa-2+\frac{2}{\sqrt{3}}\right)\left(\sigma+(4+2\sqrt{3})\kappa-2-\frac{2}{\sqrt{3}}\right)=8\mu^2-\frac{4}{3}.
\]
with solutions given by hyperbolas in $\sigma,\kappa$, with a double point at $|\mu|=1/\sqrt{6}\sim 0.408$.

\begin{figure}[h!]
\includegraphics[width=0.27\textwidth]{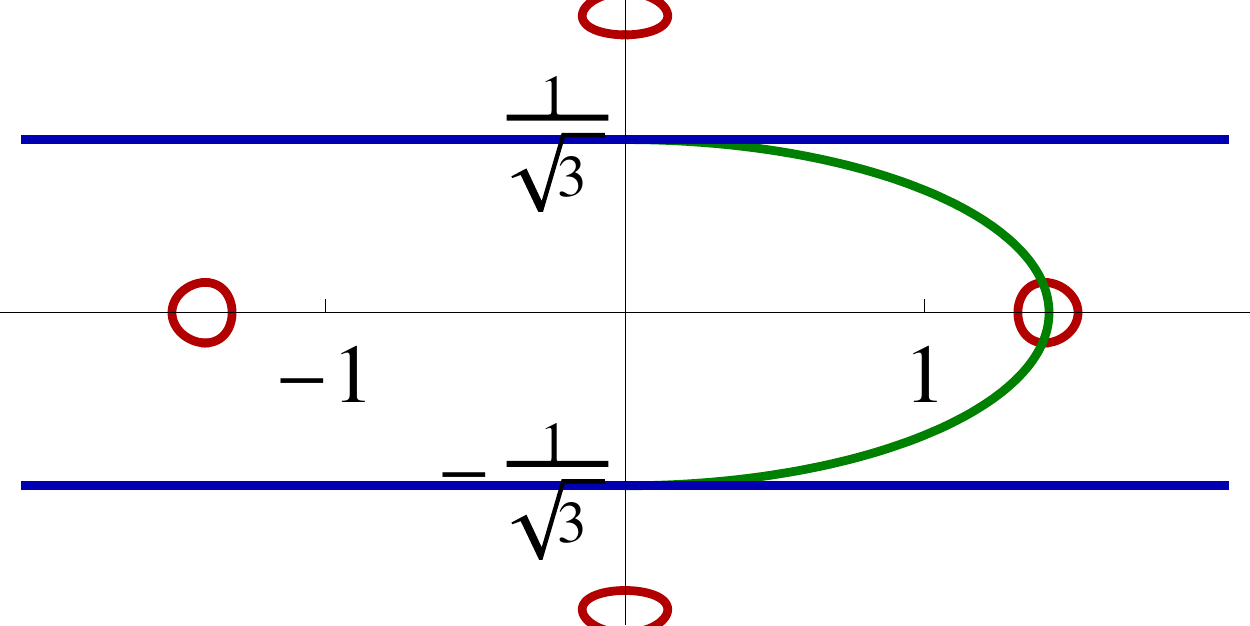}\hfill
\includegraphics[width=0.27\textwidth]{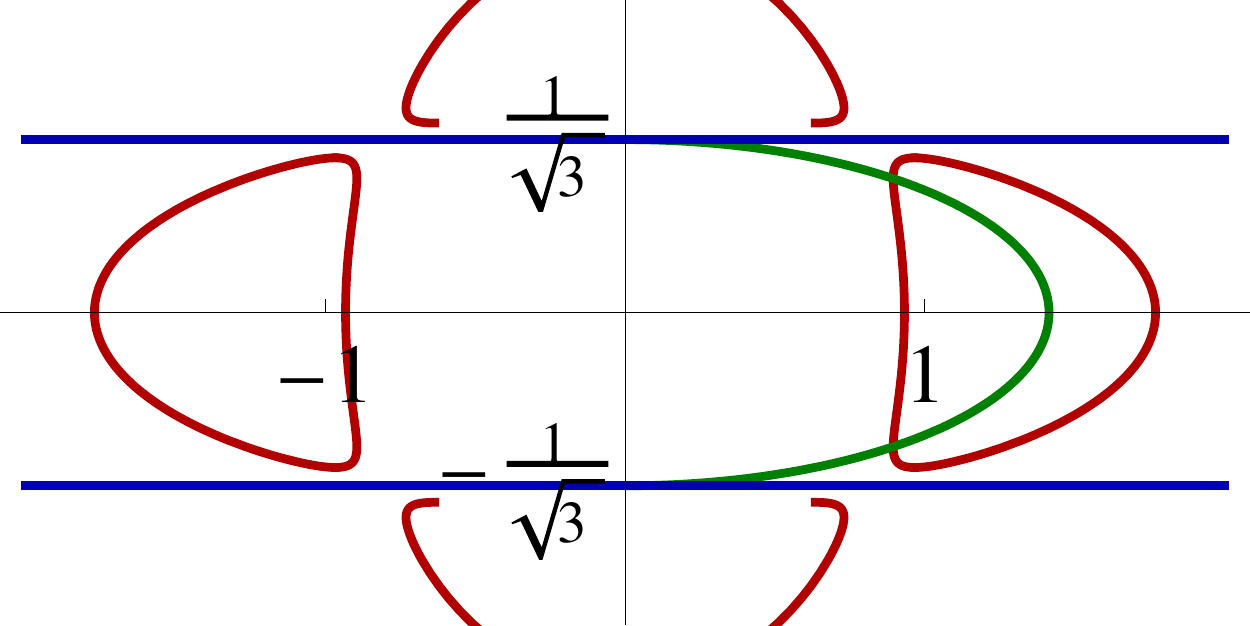}\hfill
\includegraphics[width=0.27\textwidth]{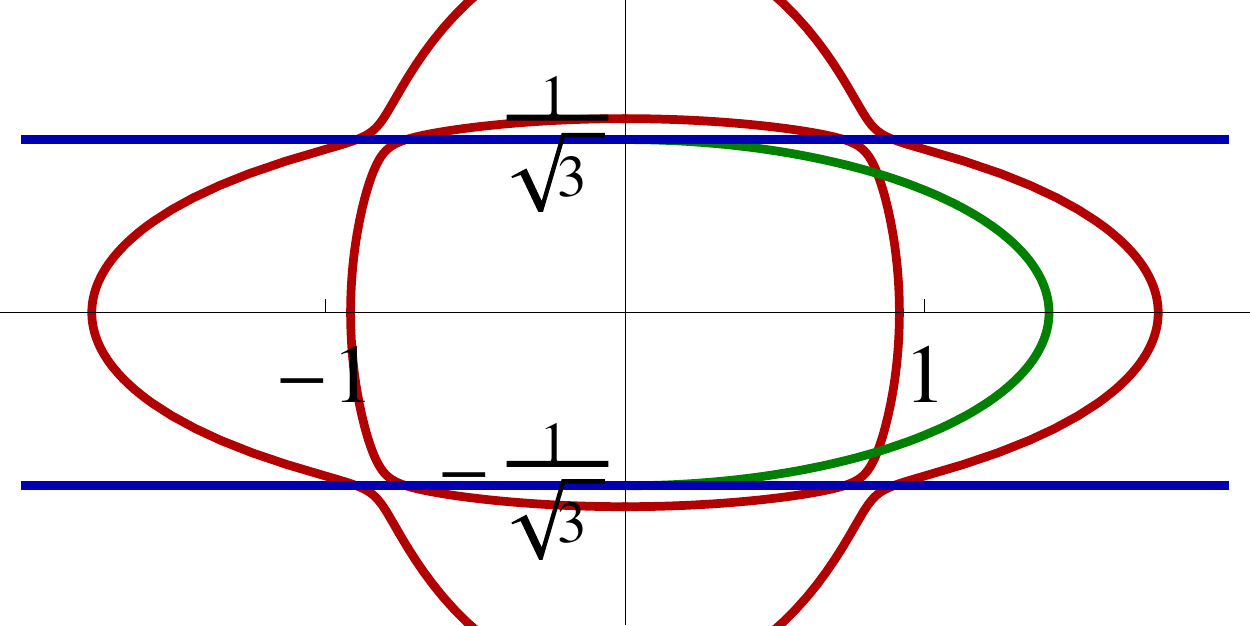}\\[0.2in]
\includegraphics[width=0.27\textwidth]{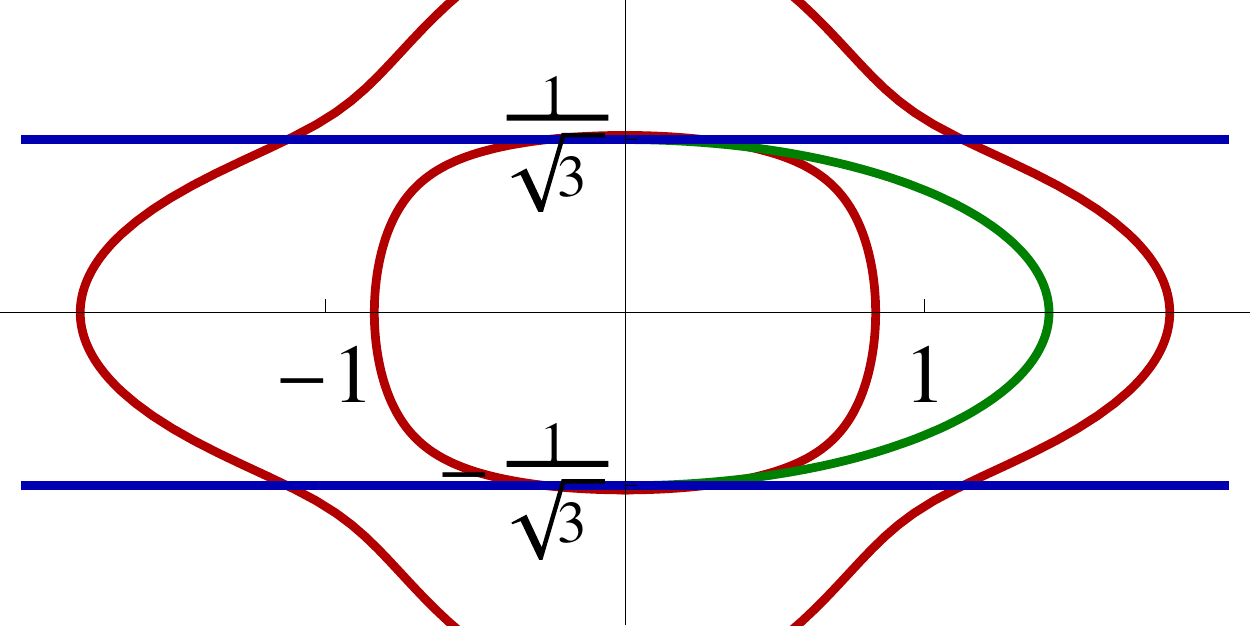}\hfill
\includegraphics[width=0.27\textwidth]{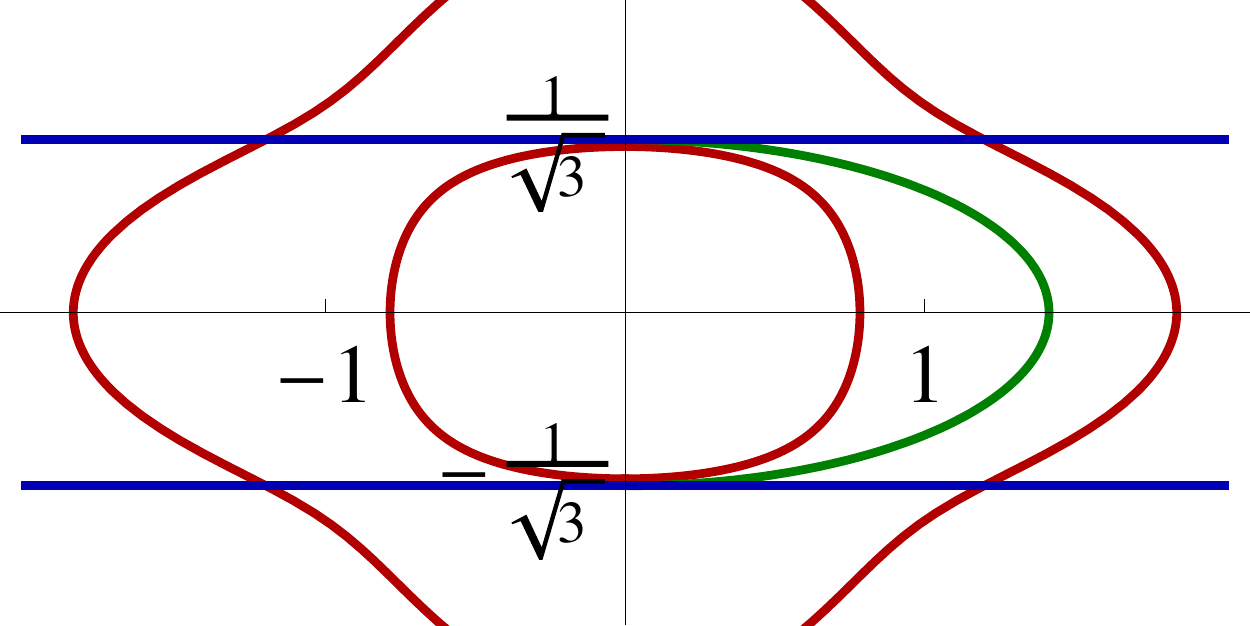}\hfill
\includegraphics[width=0.27\textwidth]{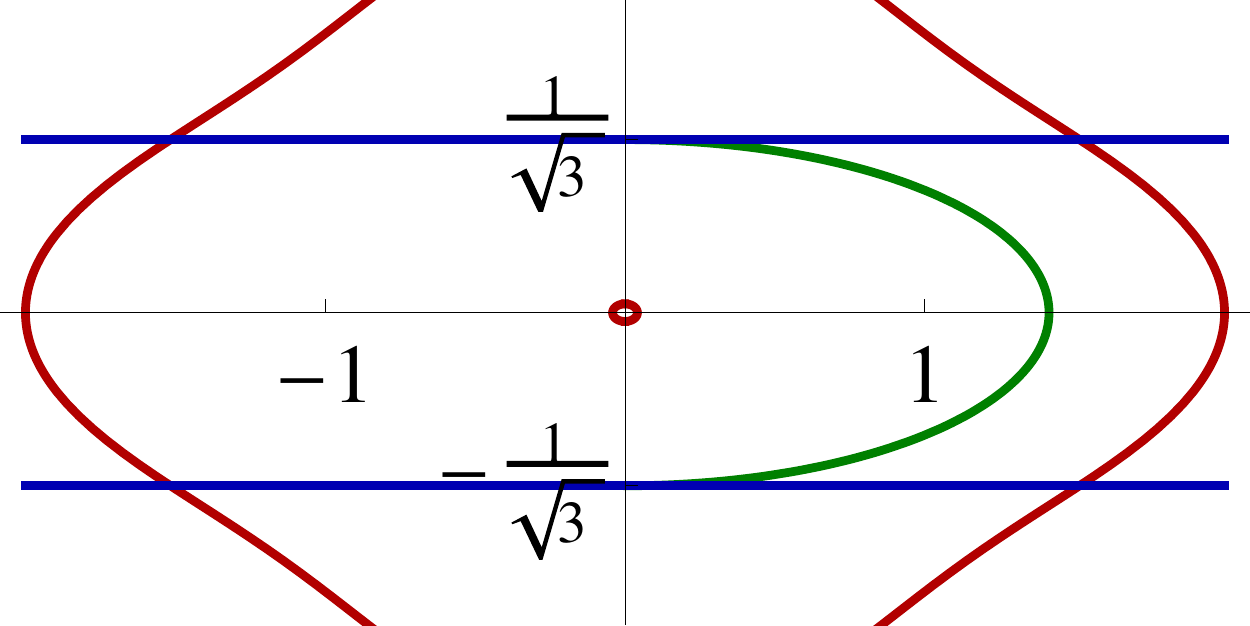}
\caption{Plots of solutions (\ref{e:kpmneu}) in the $(\tau,k)$-plane, with  reconnection crises. Also shown the restriction $\tau\leq \tau_*$ as a half-ark in the right half-plane.
Parameter values are $\mu=0.1,0.4,0.41,0.46,0.49,0.7065$, left to right, top to bottom.}\label{f:k-tau}
\end{figure}

Solving and substituting back gives the four solution branches which are plotted in the $(\tau,k)$-plane in Figure \ref{f:k-tau}. 
\begin{align}
k_{1,\pm}(\tau)&=\pm\sqrt{1-\tau^2+\sqrt{2\mu^2-\tau^2+\frac{3}{4}\tau^4}}\nonumber\\
k_{1,\pm}(\tau)&=\pm\sqrt{1-\tau^2-\sqrt{2\mu^2-\tau^2+\frac{3}{4}\tau^4}}.\label{e:kpmneu}
\end{align}
\begin{figure}[h!]
\includegraphics[width=0.27\textwidth]{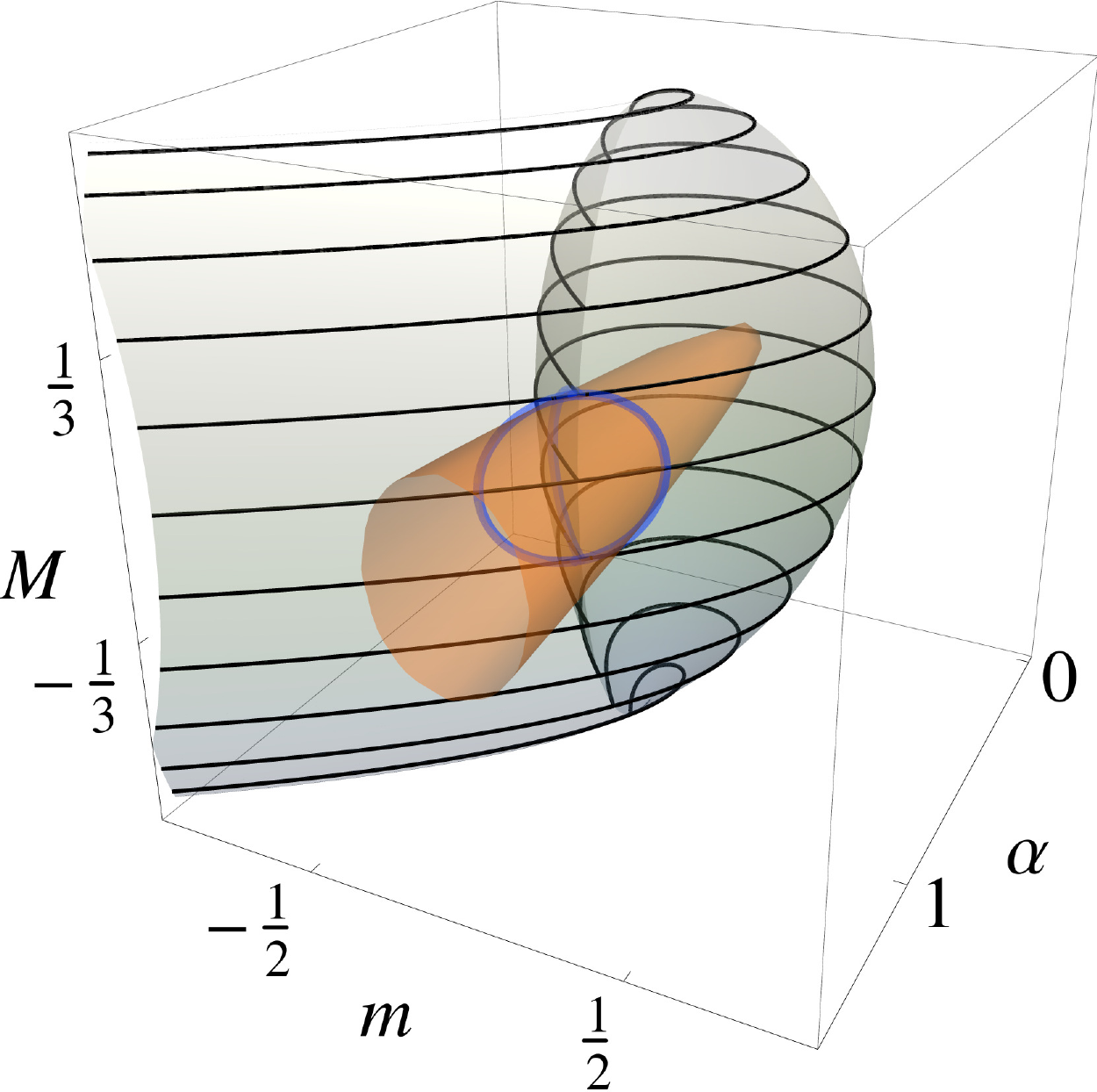}\hfill
\includegraphics[width=0.27\textwidth]{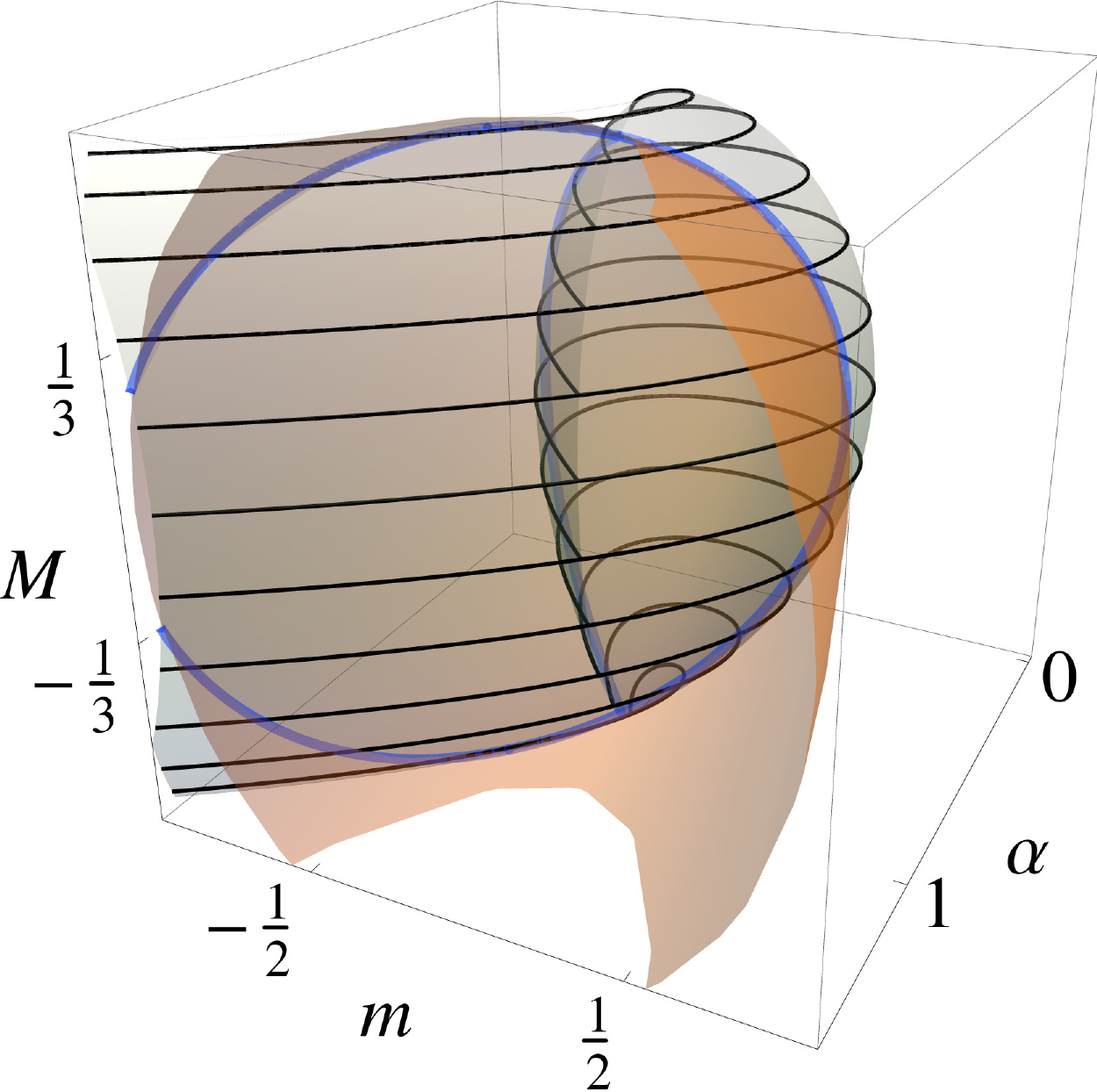}\hfill
\includegraphics[width=0.27\textwidth]{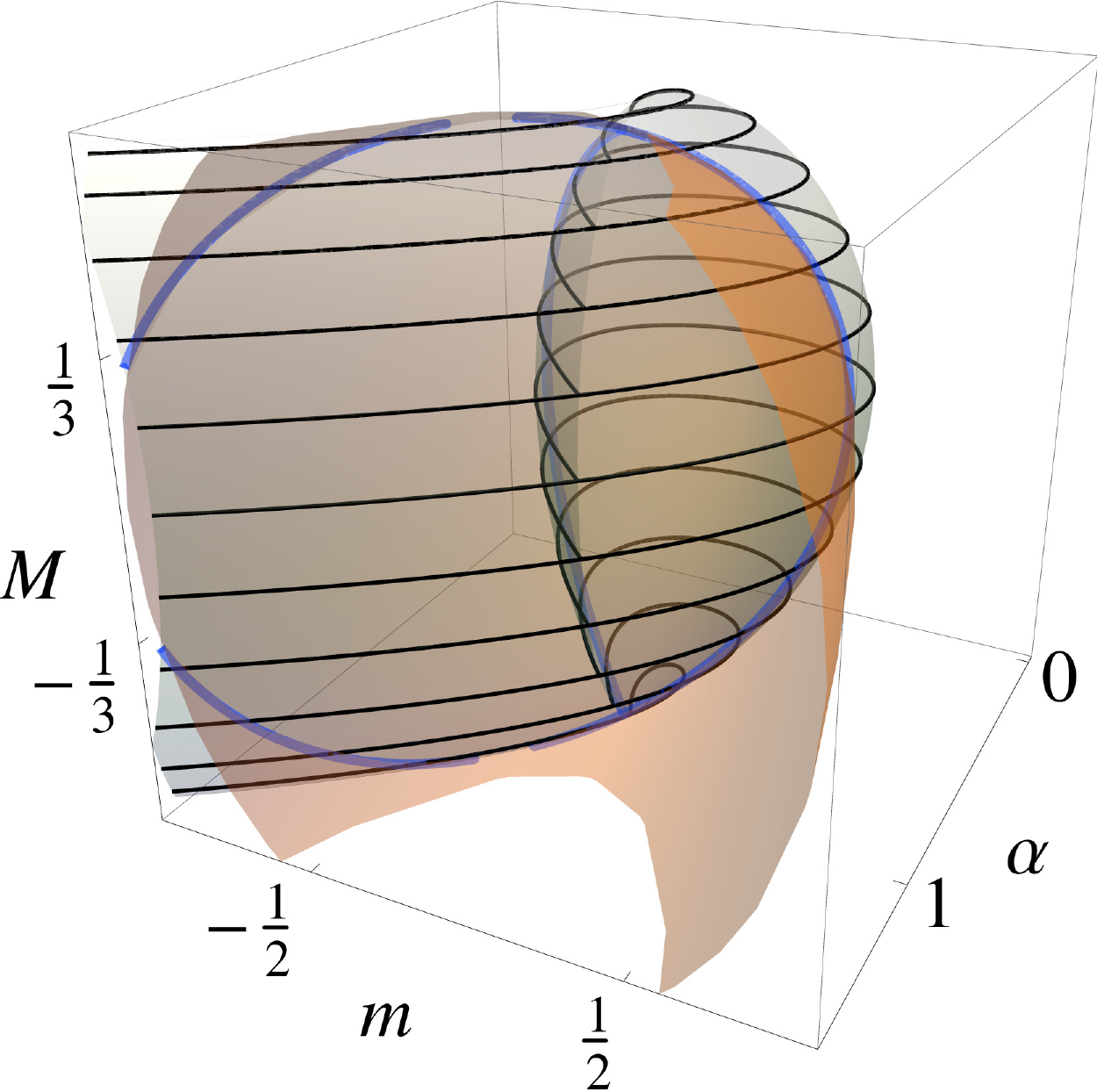}\\[0.2in]
\includegraphics[width=0.27\textwidth]{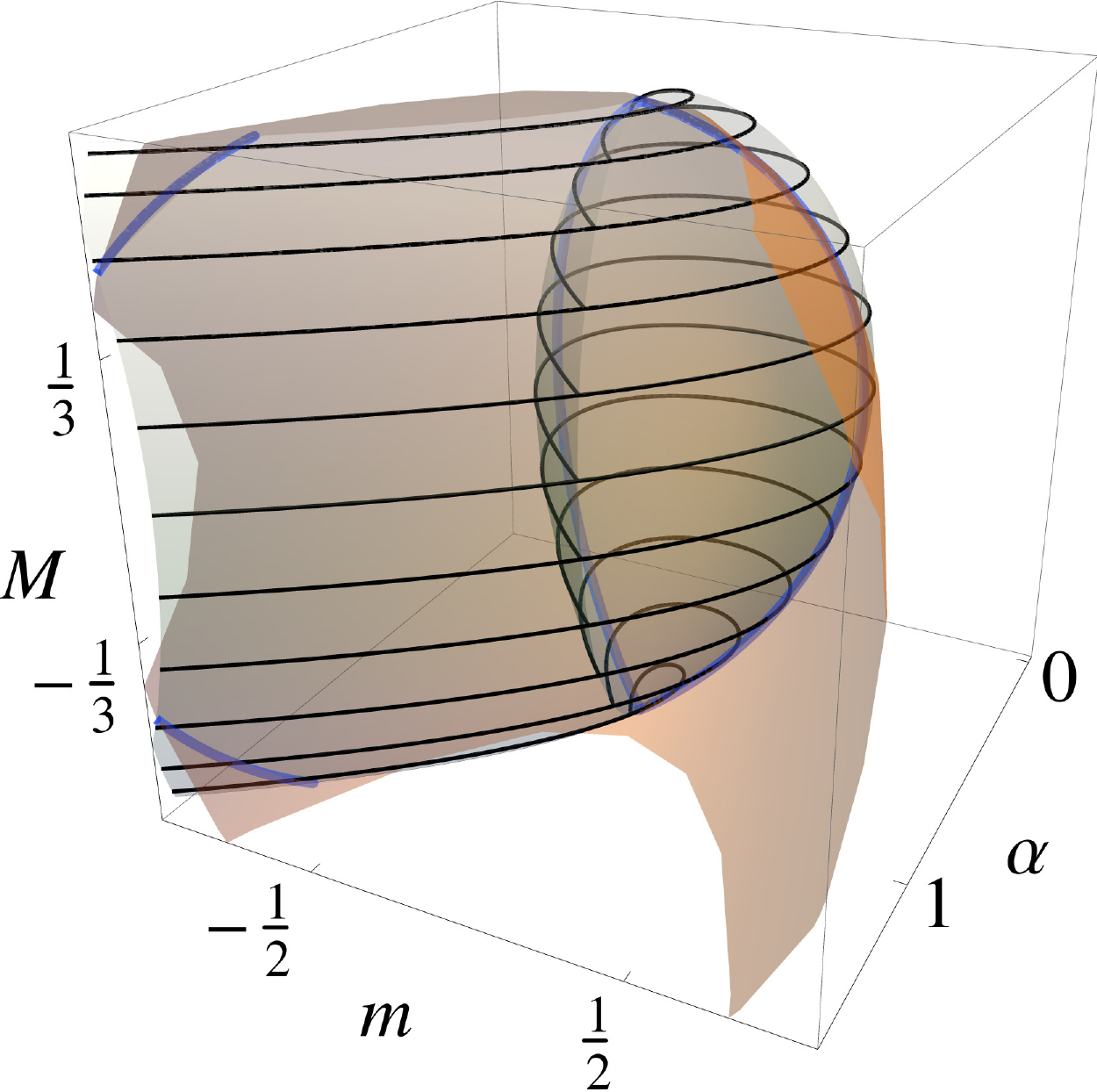}\hfill
\includegraphics[width=0.27\textwidth]{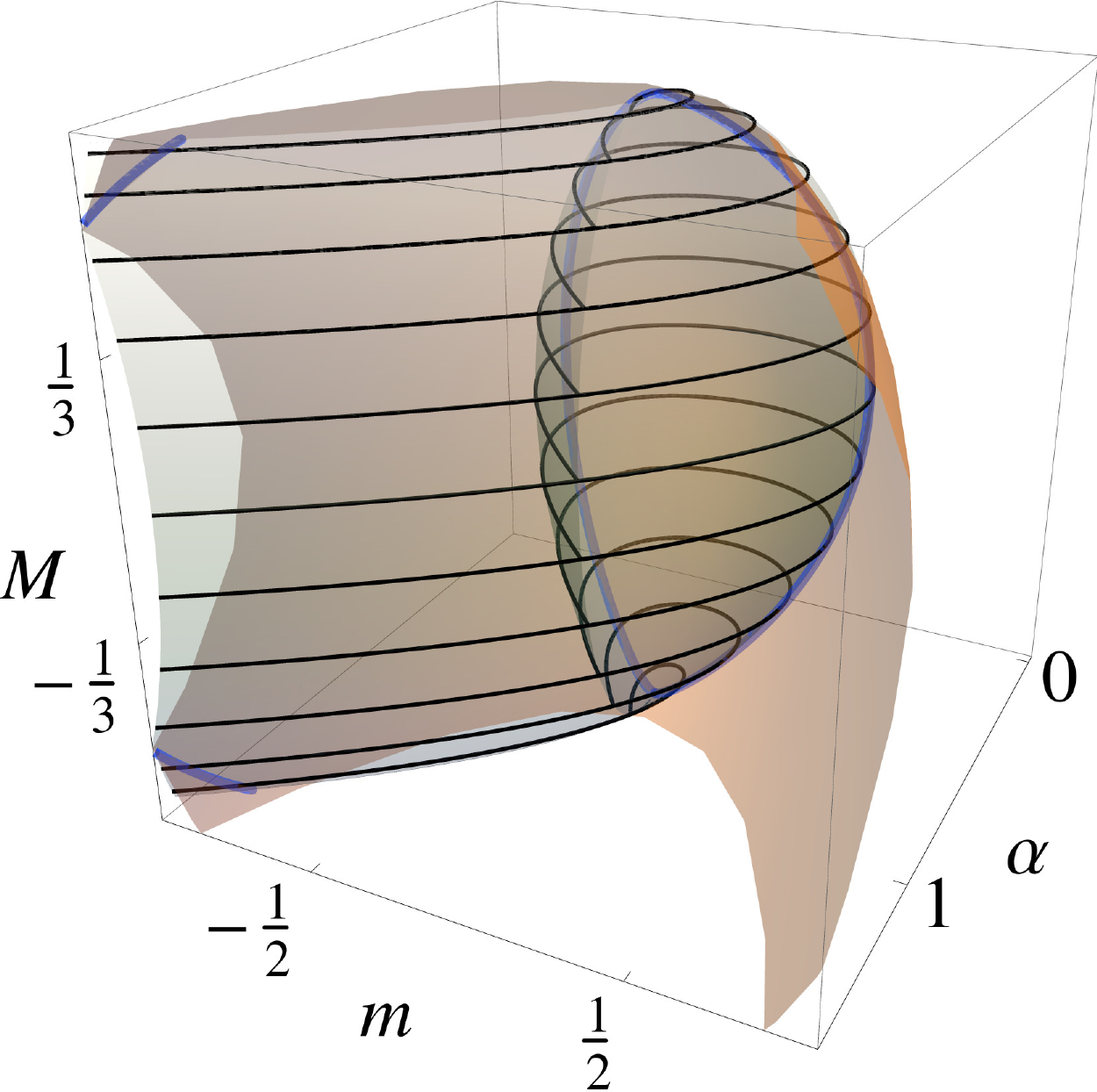}\hfill
\includegraphics[width=0.27\textwidth]{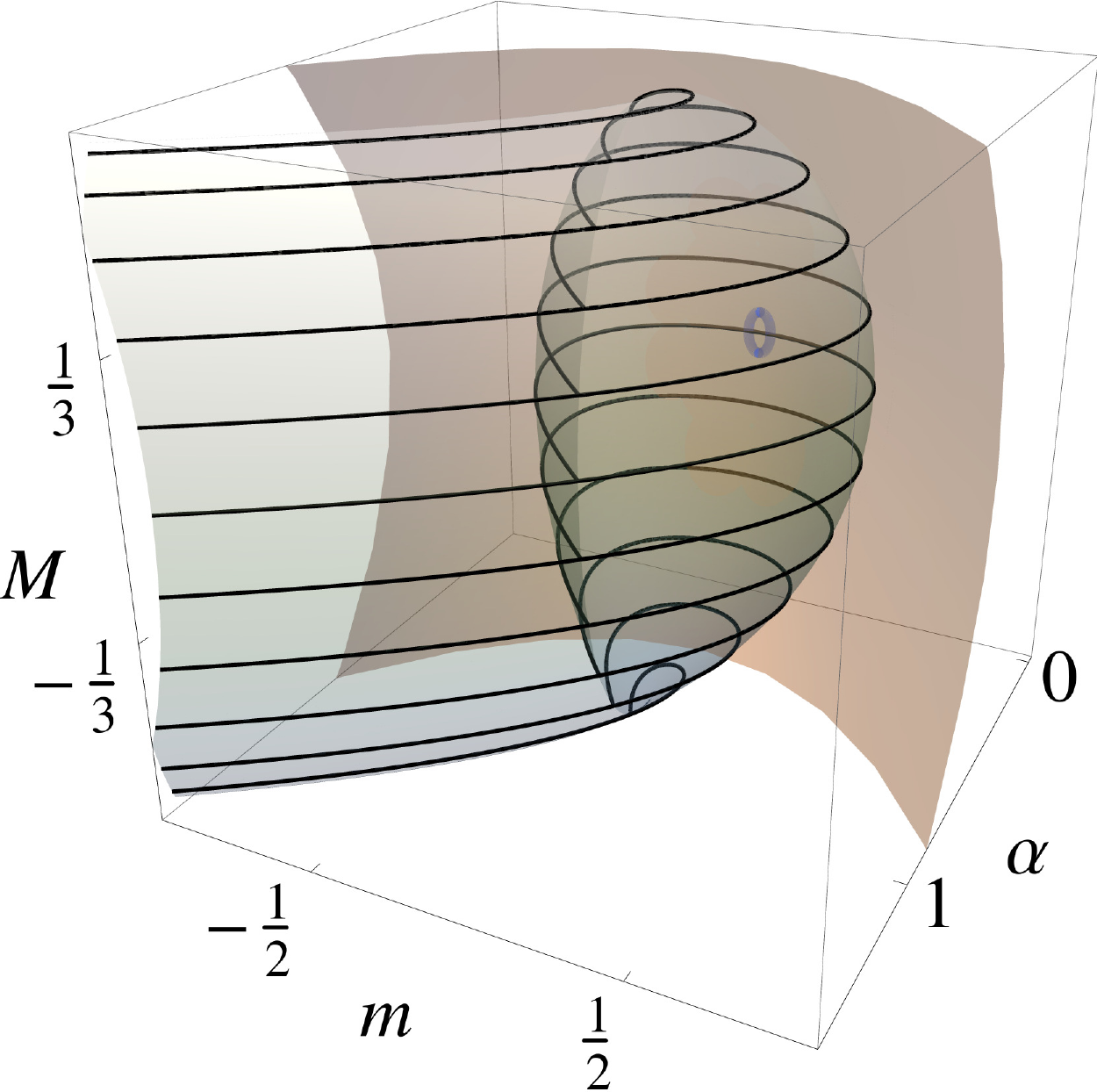}
\caption{Stable manifold in $(\alpha,m,M)$--space, together with the Neumann surface and the intersection; parameters are $\mu=0.1,0.4,0.41,0.46,0.49,0.7065$, left to right, top to bottom.}\label{f:neumann-fish}
\end{figure}
\begin{figure}[h!]
\includegraphics[width=0.27\textwidth]{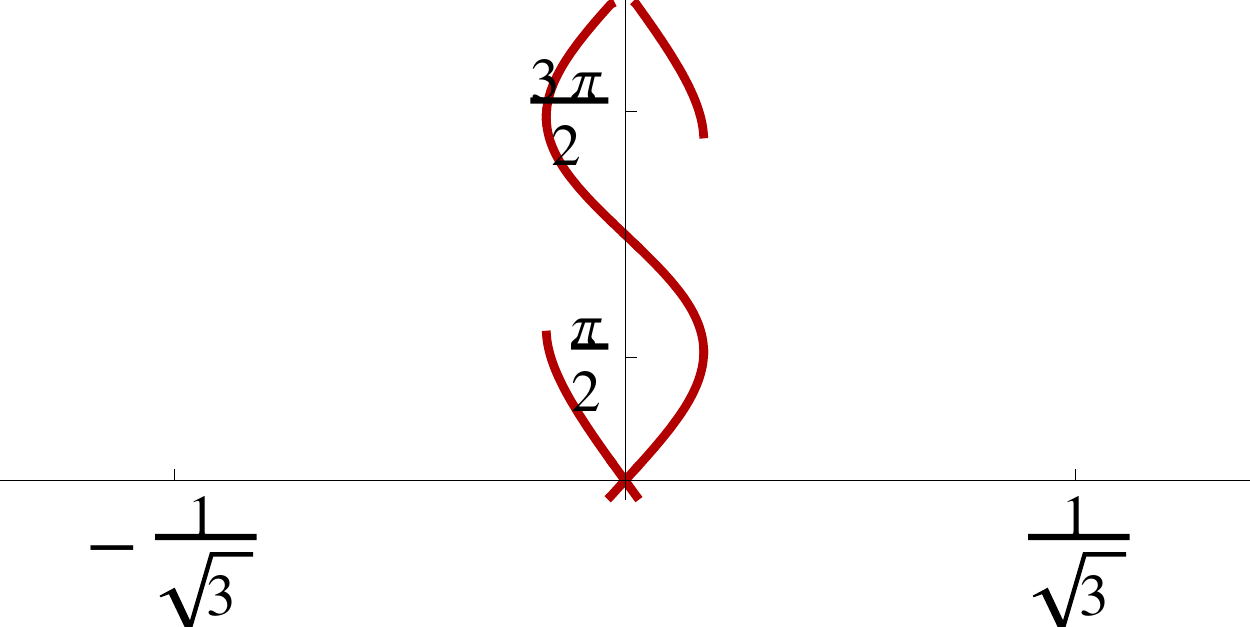}\hfill
\includegraphics[width=0.27\textwidth]{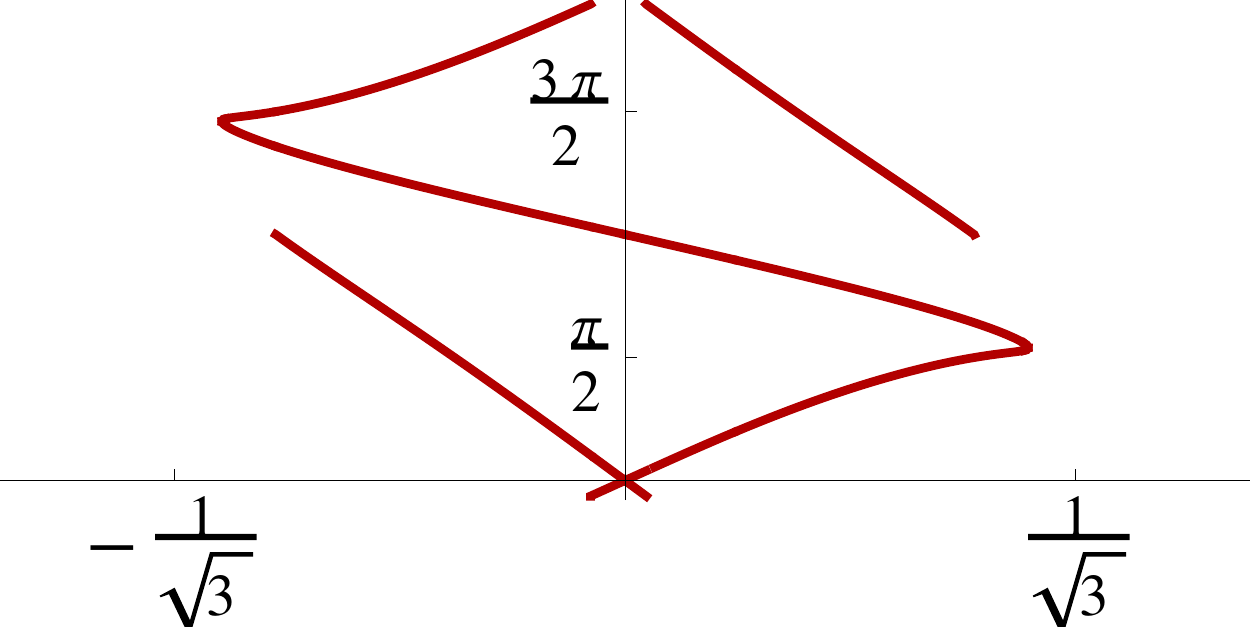}\hfill
\includegraphics[width=0.27\textwidth]{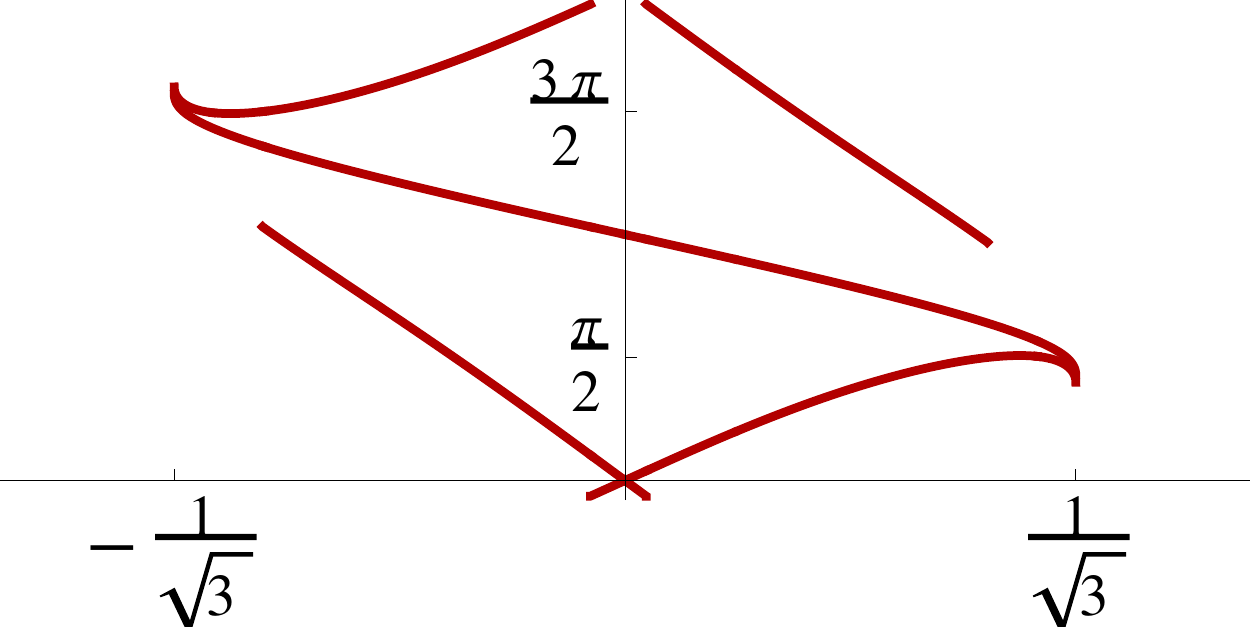}\\[0.2in]
\includegraphics[width=0.27\textwidth]{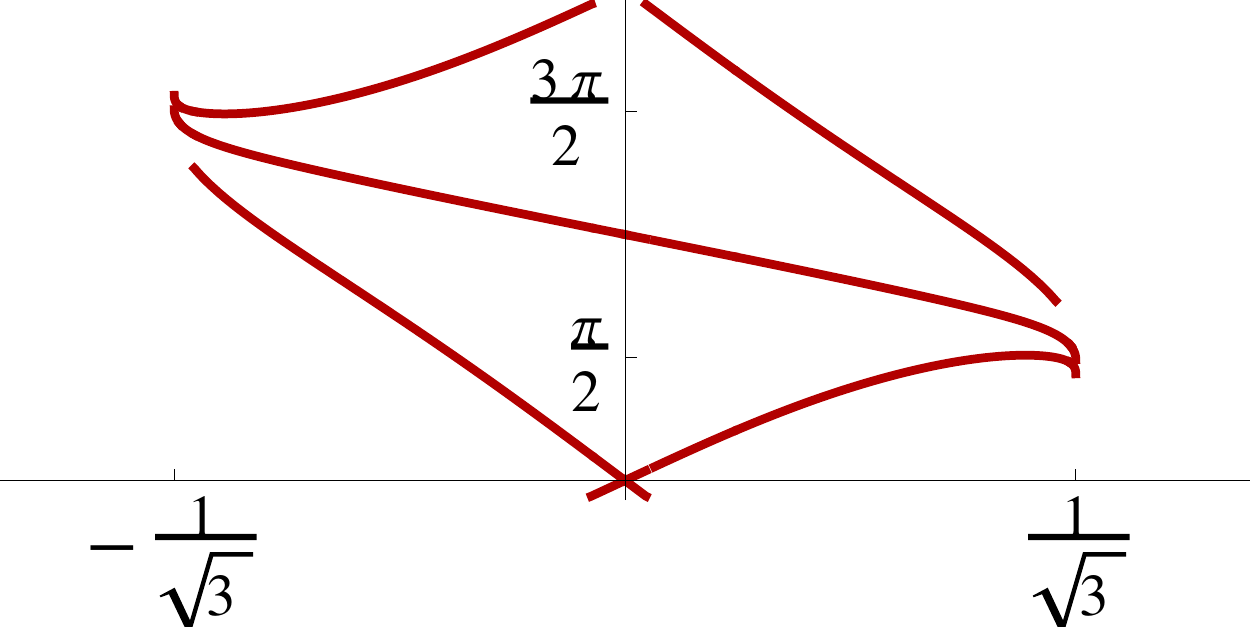}\hfill
\includegraphics[width=0.27\textwidth]{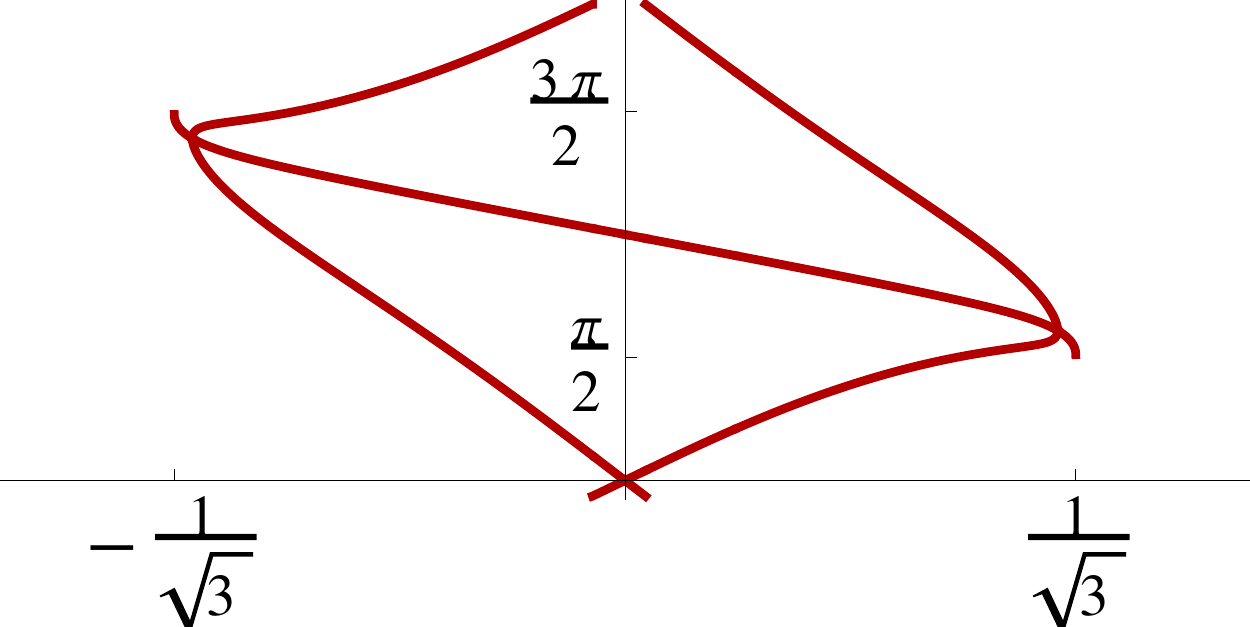}\hfill
\includegraphics[width=0.27\textwidth]{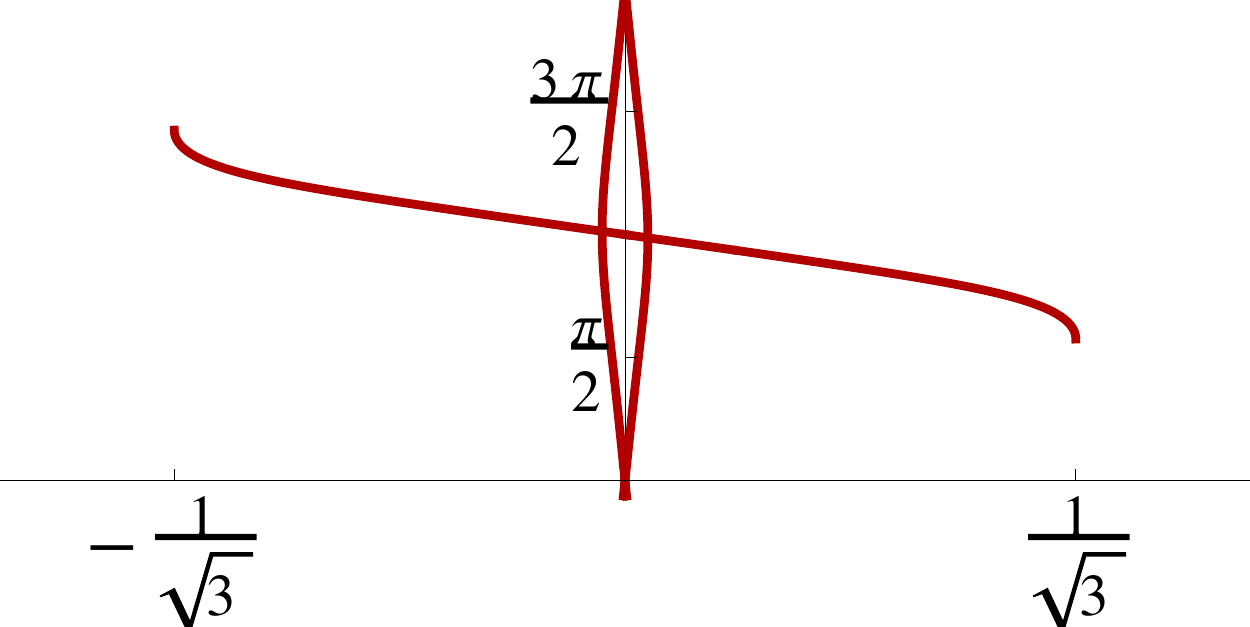}
\caption{Plots of displacement-strain curves in the $(\varphi,k)$-plane corresponding to Figure\ref{f:neumann-fish}. Note that  for $\mu>1/\sqrt{2}$, one of the displacement-strain curves has winding number 2. Open ended curves for small $\mu$ correspond to curves terminating on $\tau=\tau_*$; parameters $\mu=0.1,0.4,0.41,0.46,0.49,0.7065$, left to right, top to bottom.}\label{f:k-phi}
\end{figure}
A second crisis happens at $\mu=\sqrt{2}/3\sim 0.471$ when two curves join within the Eckhaus boundary. A third and last crisis occurs when $\mu=1/\sqrt{2}\sim 0.707$, when the origin becomes a critical point of (\ref{e:ktn}) and a curve vanishes at the origin. After this crisis, only one curve remains which covers the entire range of allowed wavenumbers. 

We can now reconstruct the phase using the original boundary condition $A=\xi$ and (\ref{e:phaseshift}).

Figure \ref{f:k-phi} shows the resulting displacement-strain curves. For small $\mu$, there are two families of boundary layers. One selects wavenumbers, with phase-winding number 1, the other family selects both wavenumber and phase, terminating on $\tau=\tau_*$. Along the curve a defect has nucleated from the boundary and the defect location diverges to infinity as the end point of the curve is approached. At the end point, the phase difference between this boundary layer and the other boundary layer is precisely $\pi$, the phase shift along the defect. For increasing $\mu$, the wavenumber-selecting boundary layer covers an increasing interval of $k$-values, until it reaches the Eckhaus boundary at the extremum, at which point it splits into two families of boundary layers, with opposite signs of displacement-strain derivative (and opposite stability properties). One of those families eventually reconnects with the third family, creating a wavenumber selecting family with winding number 2. 

\paragraph{Pure wavenumber-selection.}
A special case are boundary conditions of the form 
\begin{equation}\label{e:wns}
B=\mu A, \quad \mu\in\C,
\end{equation}
fixed. Since the boundary conditions are invariant under the gauge symmetry, the phase of boundary layers is arbitrary once we find a boundary layer. In invariant coordinates, we find the equivalent of (\ref{e:bcinv}),
\[
m-\rmi M=2\bar{\mathcal{M}}=2{\mu}\alpha,
\]
a curve rather than the surface we had seen in the previous example.  Together with (\ref{e:bcws}), we find 6 equations in 6 variables, which reduce to the system of polynomial equations in $\tau$ and $k$
\begin{align*}
2\Im(\mu) \left(2k^2+\frac{1}{2}\tau^2\right)&=-2k(1-k^2)\\
2\Re(\mu) \left(2k^2+\frac{1}{2}\tau^2\right)&=\tau\left(\frac{1}{2}\tau^2-\left(1-3k^2\right)\right).
\end{align*}
One readily notices that for $\Re(\mu)=0$, either $\tau=\pm\tau_*$ and $k=-\Im(\mu)$, the compatible equilibrium, or $\tau=0$, $\Im(\mu)\pm\sqrt{1+(\Im(\mu))^2}$.
For $|\Im(\mu)|<1/\sqrt{3}$, the first $k$-solution lies within the Eckhaus boundary, for $|\Im(\mu)|>1/\sqrt{3}$, the second $k$-solution 
does. 

\begin{figure}
\includegraphics[width=.27\textwidth]{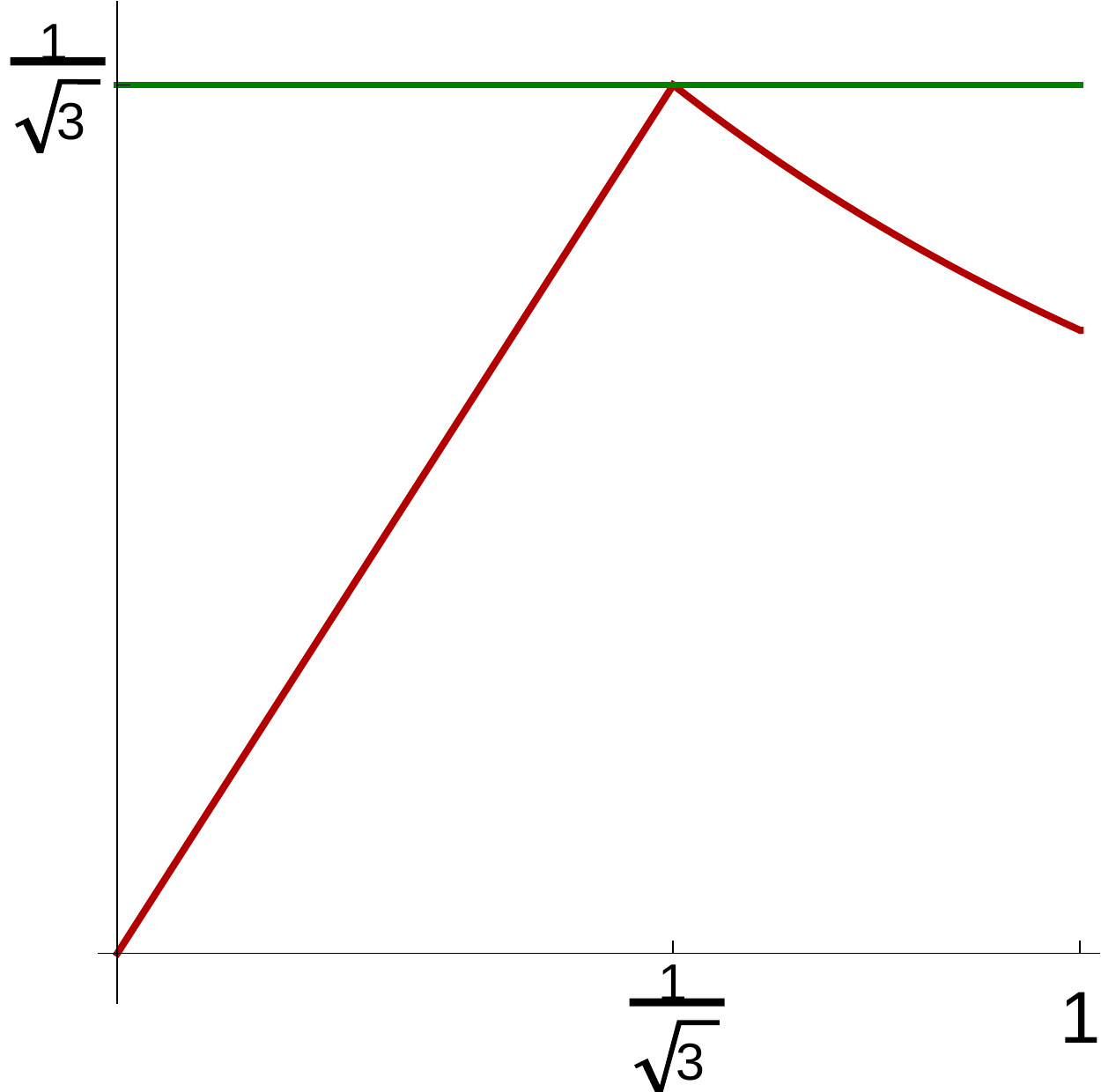}\hfill
\includegraphics[width=.27\textwidth]{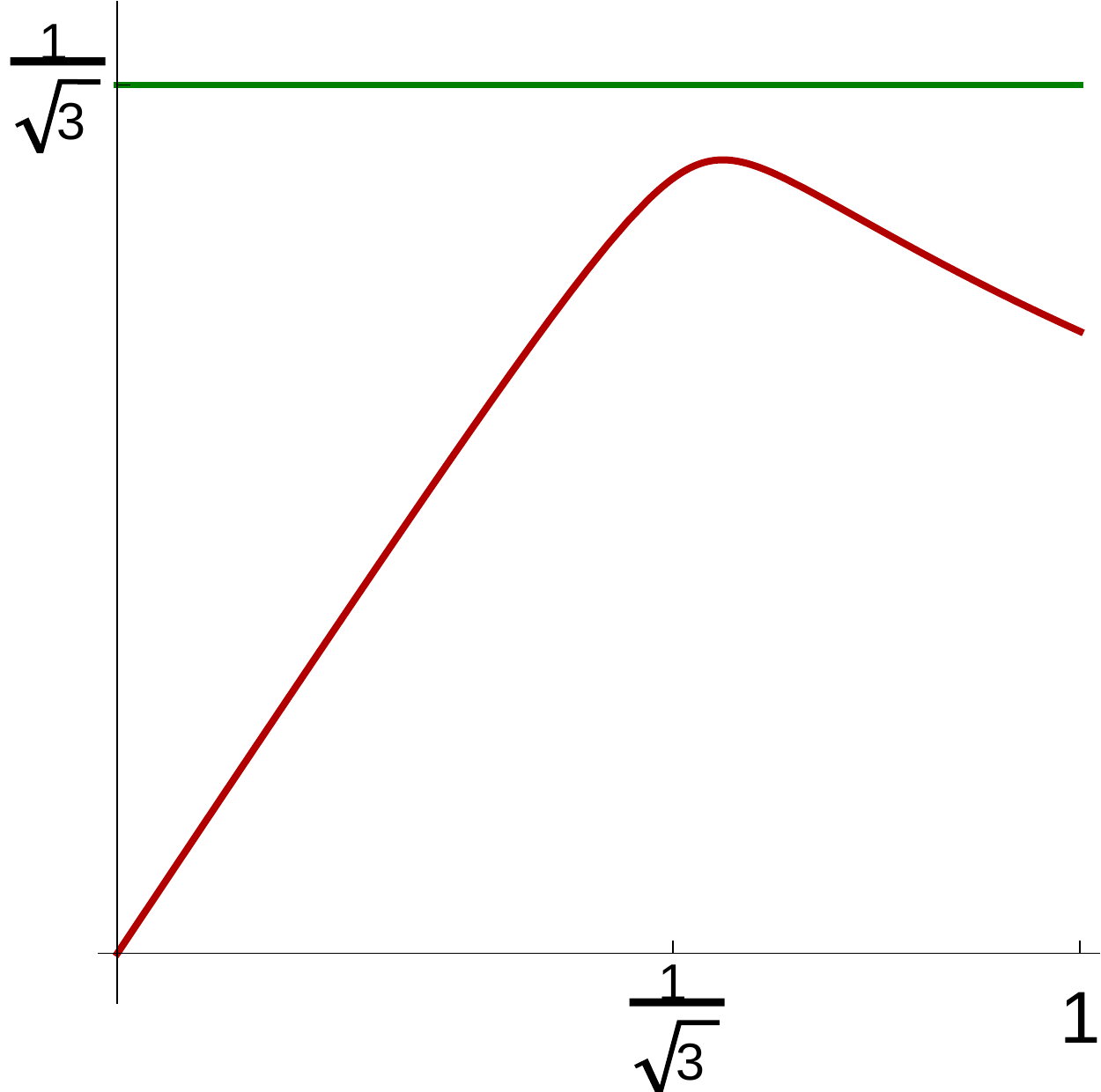}\hfill
\includegraphics[width=.27\textwidth]{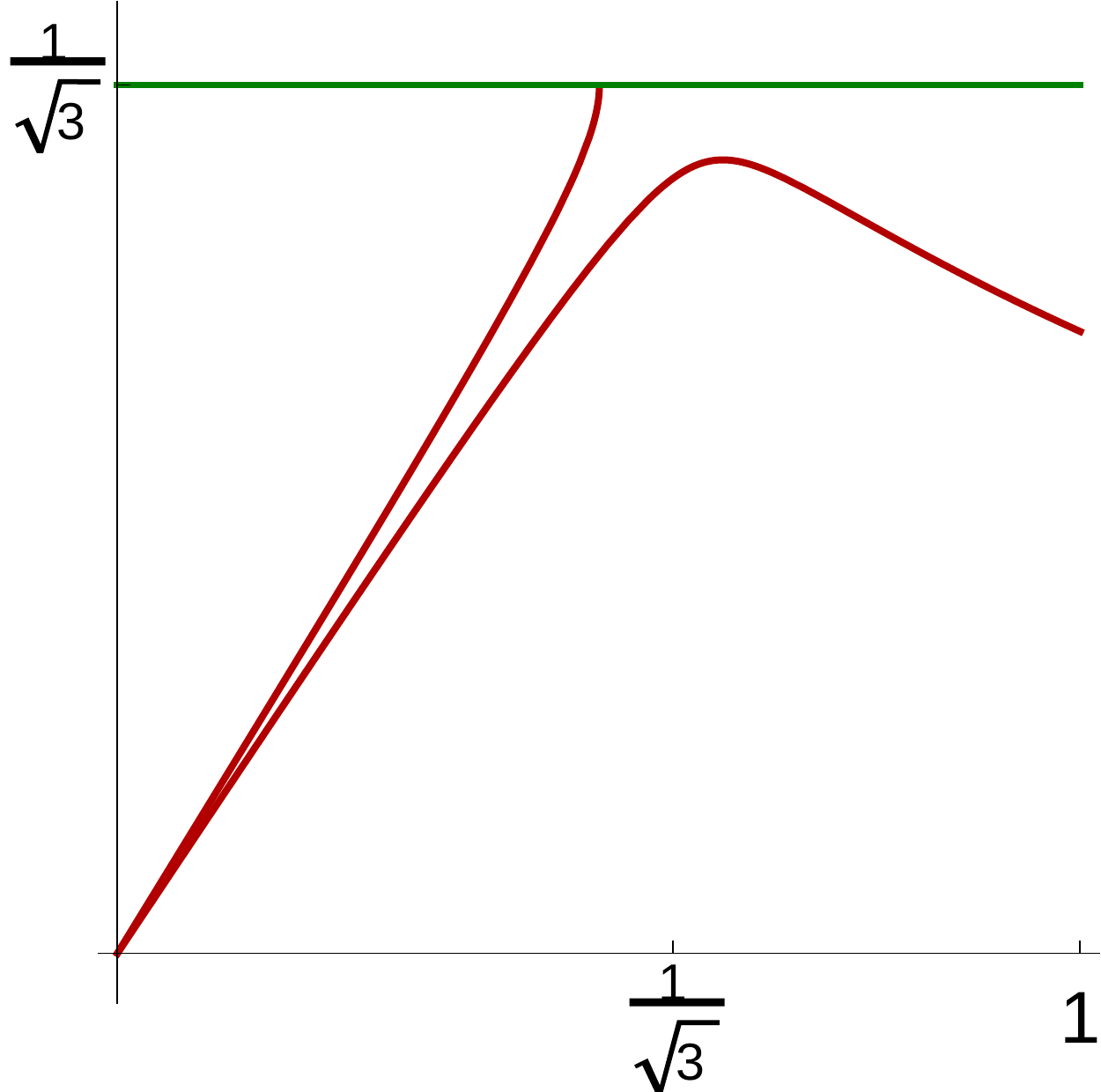}\hfill
\caption{Plot of selected wavenumber with boundary conditions (\ref{e:wns}) as a function of $\Im(\mu)$ for fixed $\Re(\mu)=-.155,0,.155$ from left to right.
}\label{f:purews}
\end{figure}

\begin{figure}
\includegraphics[width=.27\textwidth]{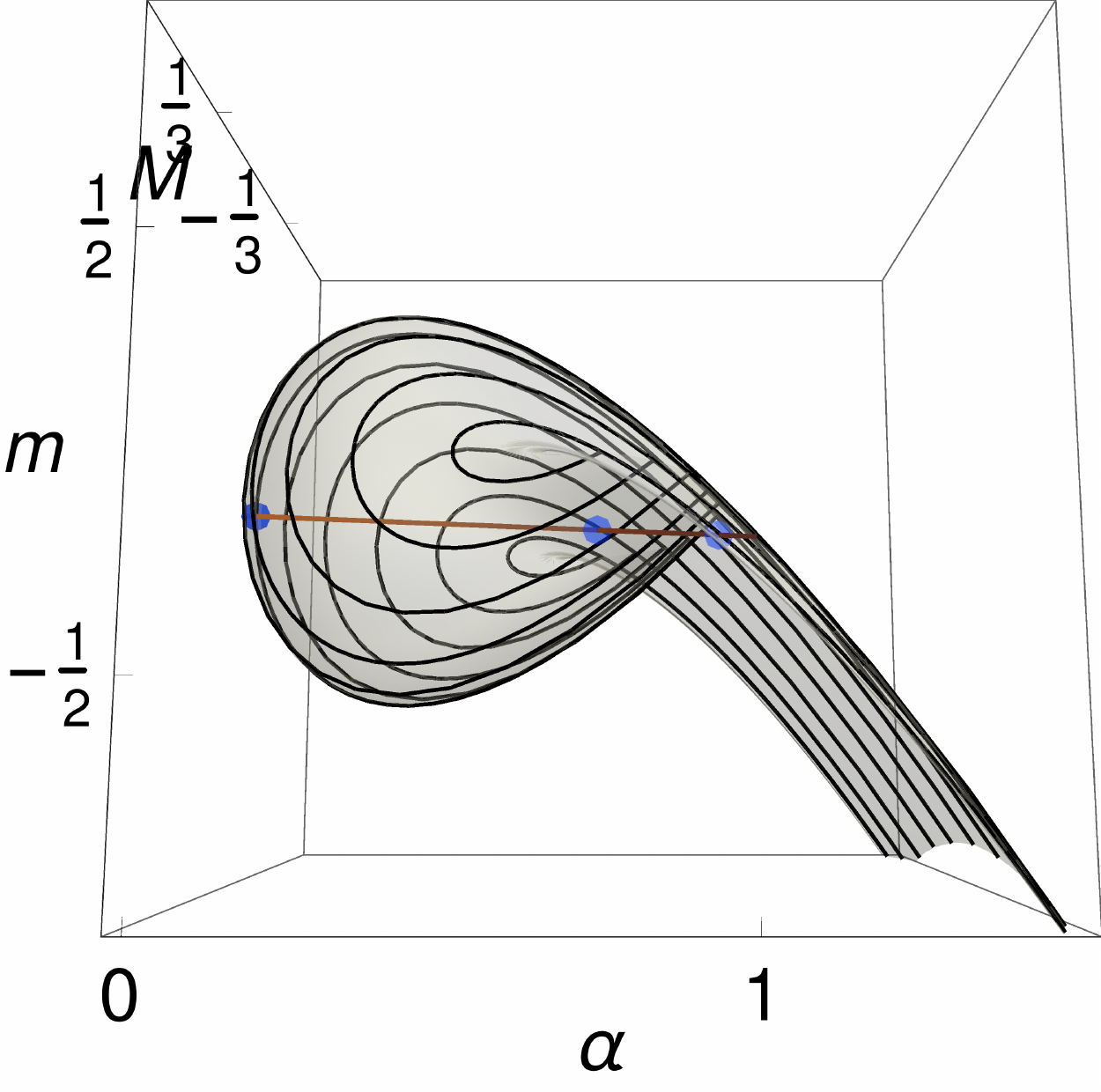}\hfill
\includegraphics[width=.27\textwidth]{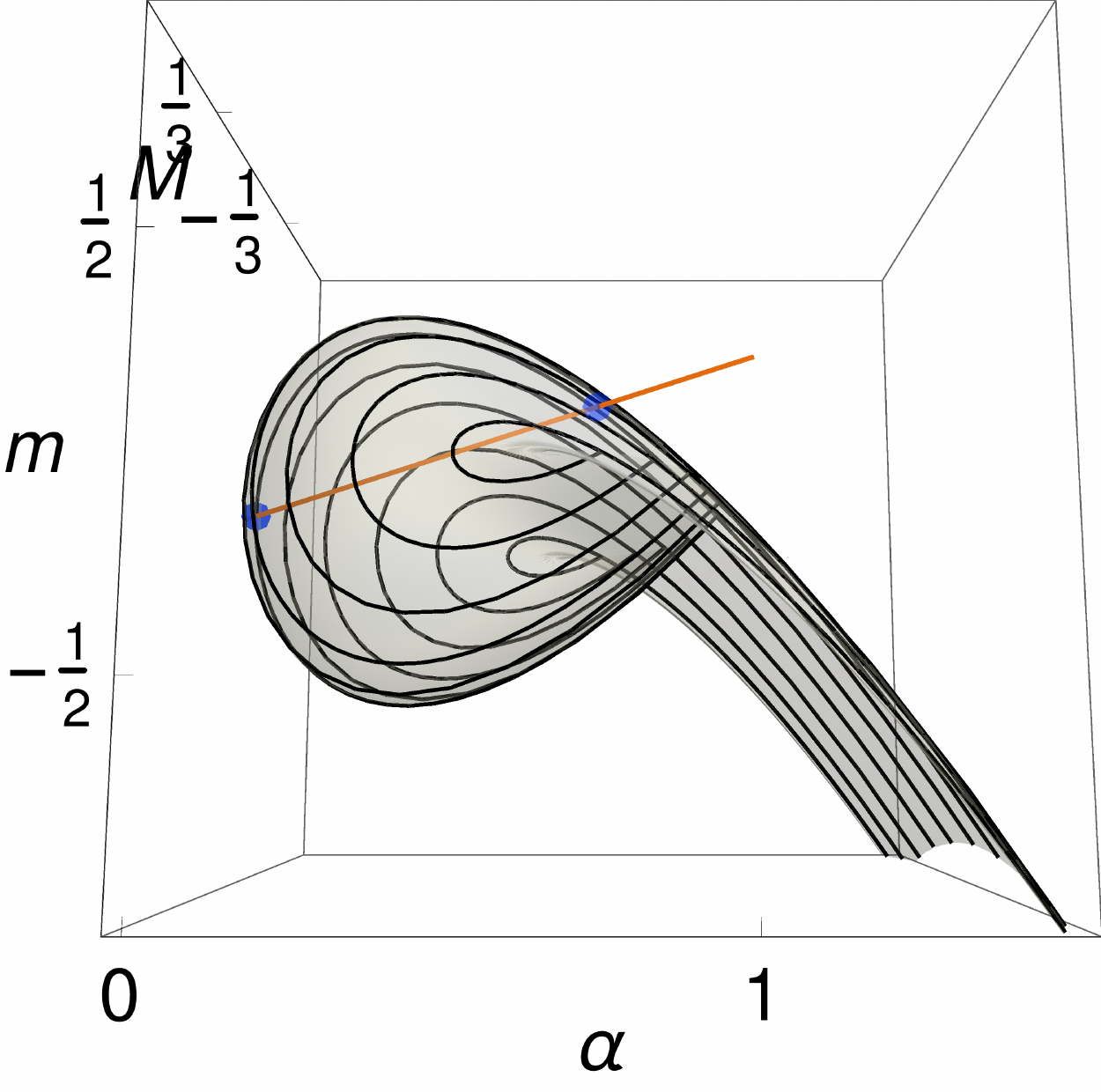}\hfill
\includegraphics[width=.27\textwidth]{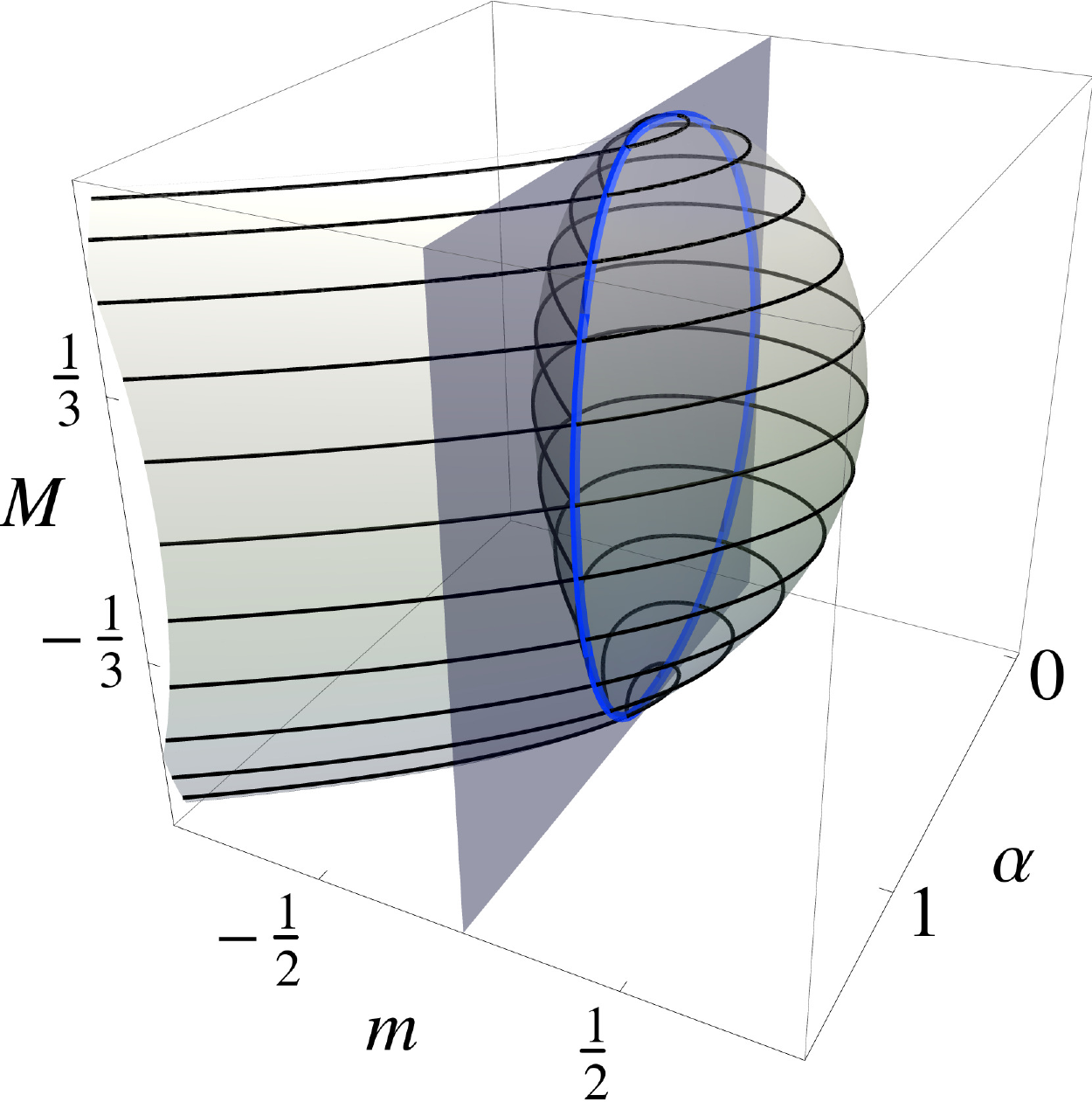}
\caption{Stable manifold in $(\alpha,m,M)$-space, together with the boundary conditions and intersections. Pure wavenumber selection ($\Re(\mu)=-0.1,0.1$, $\Im\mu=0.3$, left and middle) and pure phase selection (right).  }\label{f:purews-fish}
\end{figure}

For small $\Re(\mu)$, the intersection with the equilibrium splits into a pair of intersections close to $\pm\tau_*$, respectively. For $\Re(\mu)>0$, $|\tau|>|\tau_*|$, so that only one intersection with the level set of the Hamiltonian yields an intersection with the stable manifold. For $\Re(\mu)<0$, both intersections yield boundary layers. The dependence of $k$ on $\Im(\mu)$ is shown in Figure \ref{f:purews}. The associated pictures in $(\alpha,m,M)$-space are shown in Figure \ref{f:purews-fish}.

\paragraph{Pure phase selection.}

Boundary conditions $\Im(A)=\Re(B)=0$ encode the reflection symmetry $\mathcal{S}$. In fact, all periodic solutions satisfy this boundary condition for an appropriate phase shift. In order to study these boundary conditions more systematically, we write for the group orbit of the boundary condition $A=\rme^{\rmi\phi}(\xi+\bar{\xi}),B=\rme^{\rmi\phi}(\xi-\bar{\xi})$, which gives 
\[
m=0,\quad \mathcal{M}=-2\rmi(\bar{\xi}^2-\xi^2).
\]
Solving the remaining equations (\ref{e:bcws}), we find 4 solutions corresponding to the exact periodic pattern, $\tau=\pm\tau_*$, $\xi=\pm(1+\rmi k)\sqrt{1-k^2}$, and two solutions $\tau=0, \xi=(\rmi(1-k^2)/2 -k)/\sqrt{2}$. The first four solutions correspond to only two actual ``boundary layers'' ($\tau=\pm\tau_*$ yield the same boundary layer). They are in fact the exact periodic solution and the solution with a phase shift of $\pi$.  The latter two solutions correspond to an appropriate shift of the defect and yield a phase shift of $\pi/4$. 

A completely equivalent analysis also yields the boundary layers for the third reflection symmetry  $\Re A=0, \Im B=0$, with the same wavenumber selection, and phases shifted by $\pi/4$. 

\begin{figure}[h!]
\includegraphics[width=0.27\textwidth]{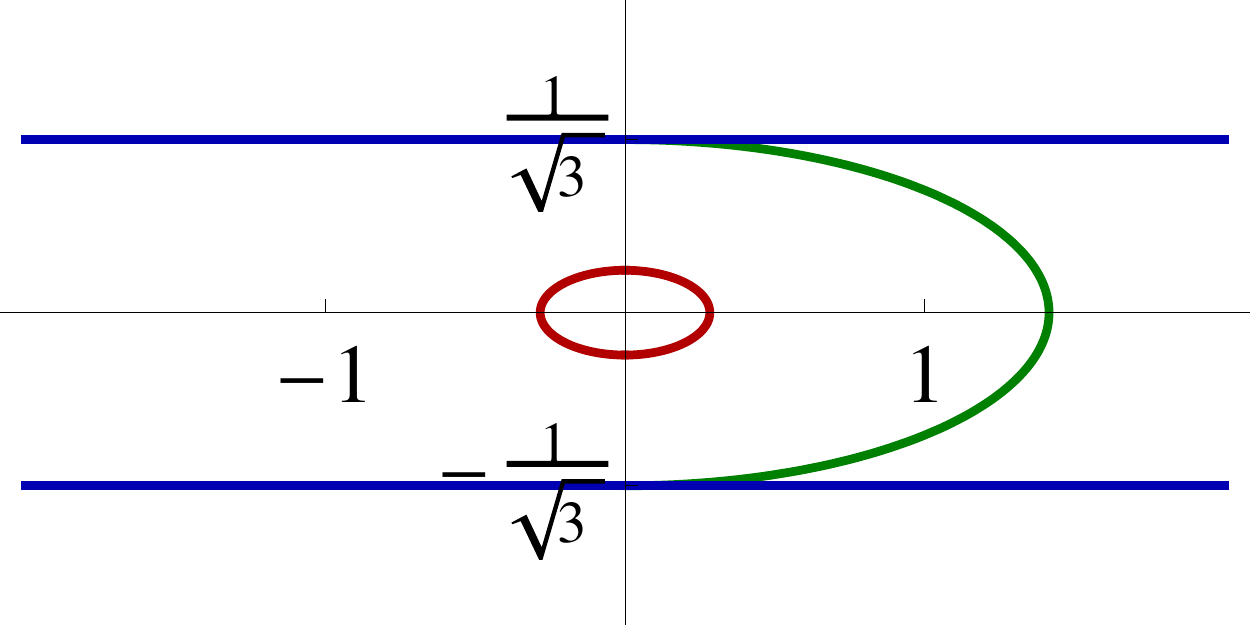}\hfill
\includegraphics[width=0.27\textwidth]{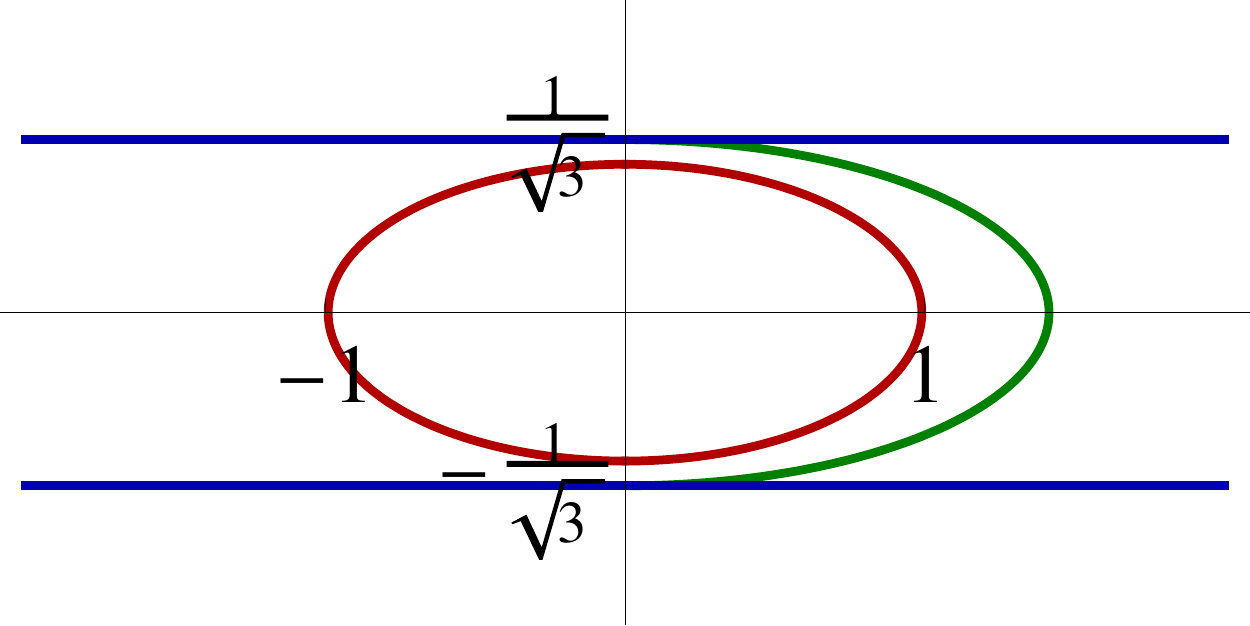}\hfill
\includegraphics[width=0.27\textwidth]{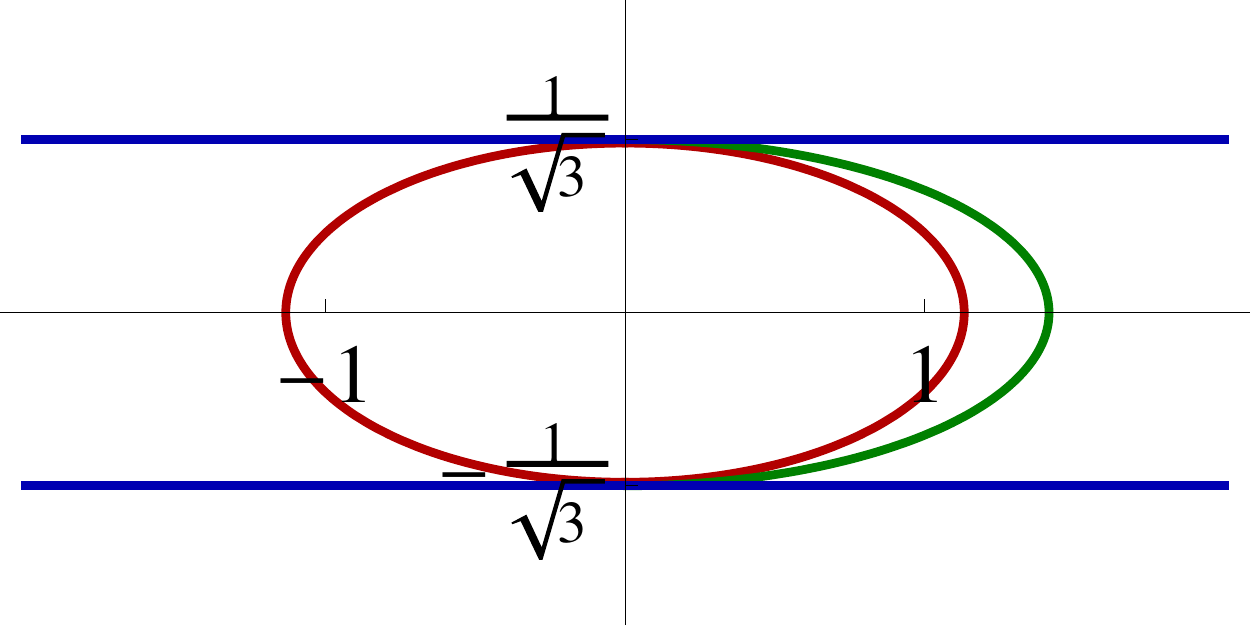}\\[0.2in]
\includegraphics[width=0.27\textwidth]{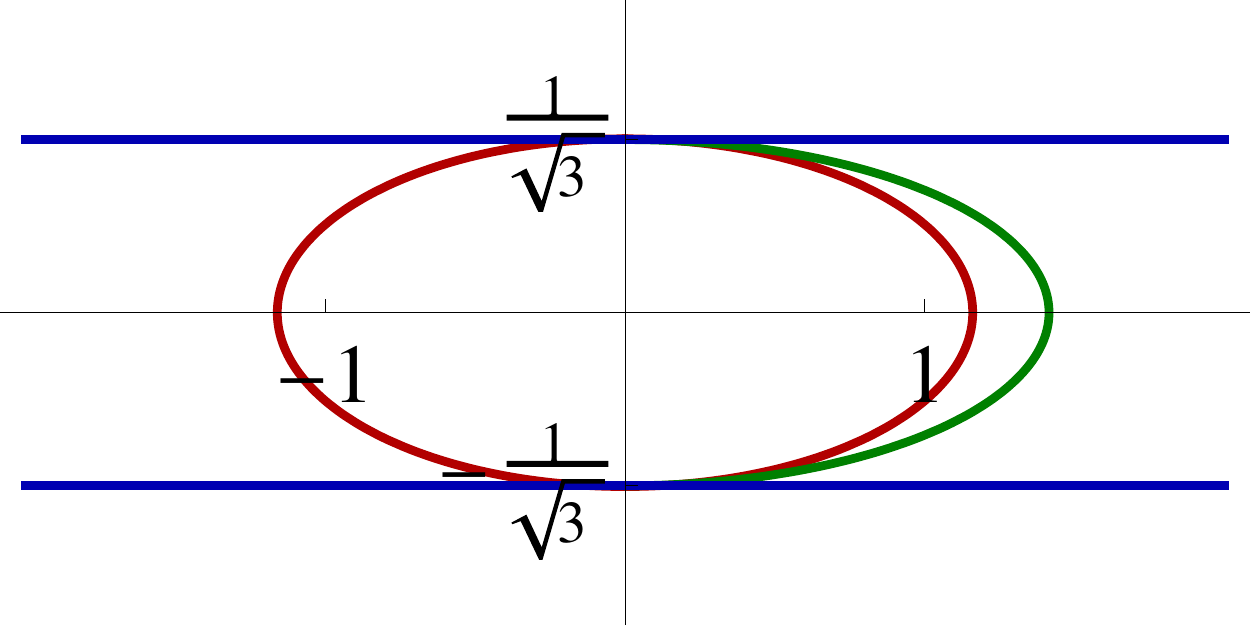}\hfill
\includegraphics[width=0.27\textwidth]{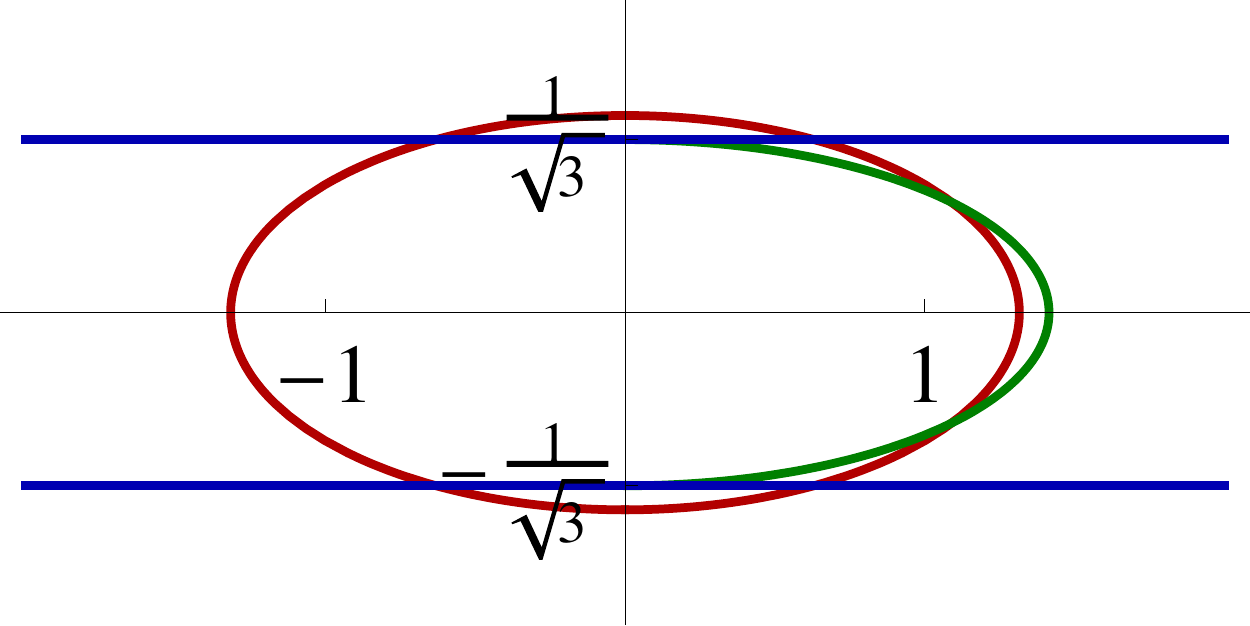}\hfill
\includegraphics[width=0.27\textwidth]{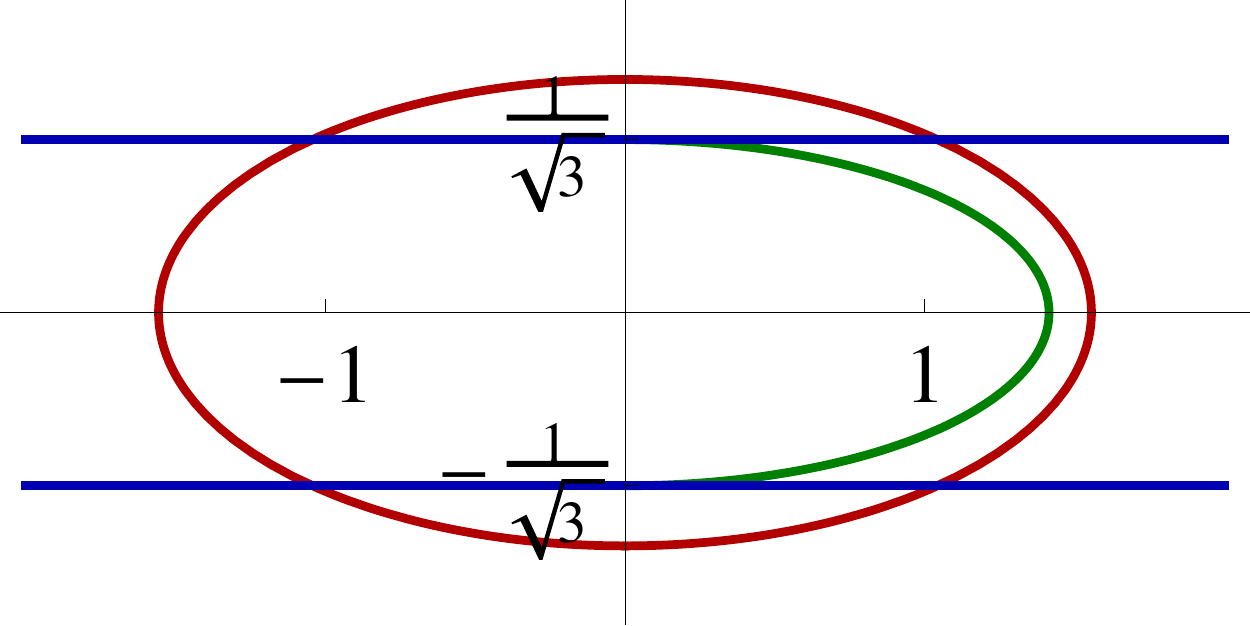}
\caption{Plots of solutions (\ref{e:kptdir}) in the $(\tau,k)$-plane, with  reconnection crises. Also shown the restriction $\tau\leq \tau_*$ as a half-ark in the right half-plane.
Parameter values are $\mu=0.2,0.7,0.93,1.1$.}\label{f:dirtau-k}
\end{figure}

\paragraph{Dirichlet boundary conditions.} 
For Dirichlet boundary conditions, $A=\mu>0$, the group invariant form gives 
\[
\{(A,B);B=\rme^{\rmi\phi}\xi,A=\rme^{\rmi\phi}\mu,\xi\in\C,\phi\in[0,2\pi)\},
\]
which gives the equivalent of (\ref{e:bcinv}),
\begin{equation}\label{e:bcinvd}
\alpha=\mu^2, \quad \mathcal{M}=\mu\bar{\xi}.
\end{equation}
Again, we can use the second equation in (\ref{e:bcinvd}) and the second and third equation in (\ref{e:bcws}) to express $\xi$ in terms of $\tau,k$. The first equations of  (\ref{e:bcinv}) and (\ref{e:bcws}) give 
\begin{equation}\label{e:kptdir}
\frac{1}{2}\tau^2+2k^2=\mu^2;
\end{equation}the solution curves, together with the restriction $\tau\leq \tau_*$ are shown in Figure \ref{f:dirtau-k}. The ellipse intersects the boundary $\tau=\tau_*(k)$ for $\mu>\sqrt{2/3}$. 

As $\mu$ is increased, we see displacement-strain curves changing from wavenumber selection (winding number 1) to phase-selection. The second curve of boundary layers terminates at defects and, as $\mu$ increases further, shrinks and vanishes on the defect associated with $k=0$. 

Note that for $\mu=1.1$, $\varphi(k)$ is not monotone. In fact, $\varphi(k)$ develops 3 critical points for $\mu<1$ in a pitchfork bifurcation. This phenomenon cannot occur for well-posed equations in the phase-diffusion systems, where an eigenvalue or resonance pole disappears at infinity at $\varphi'=0$.

\begin{figure}[h!]
\includegraphics[width=0.27\textwidth]{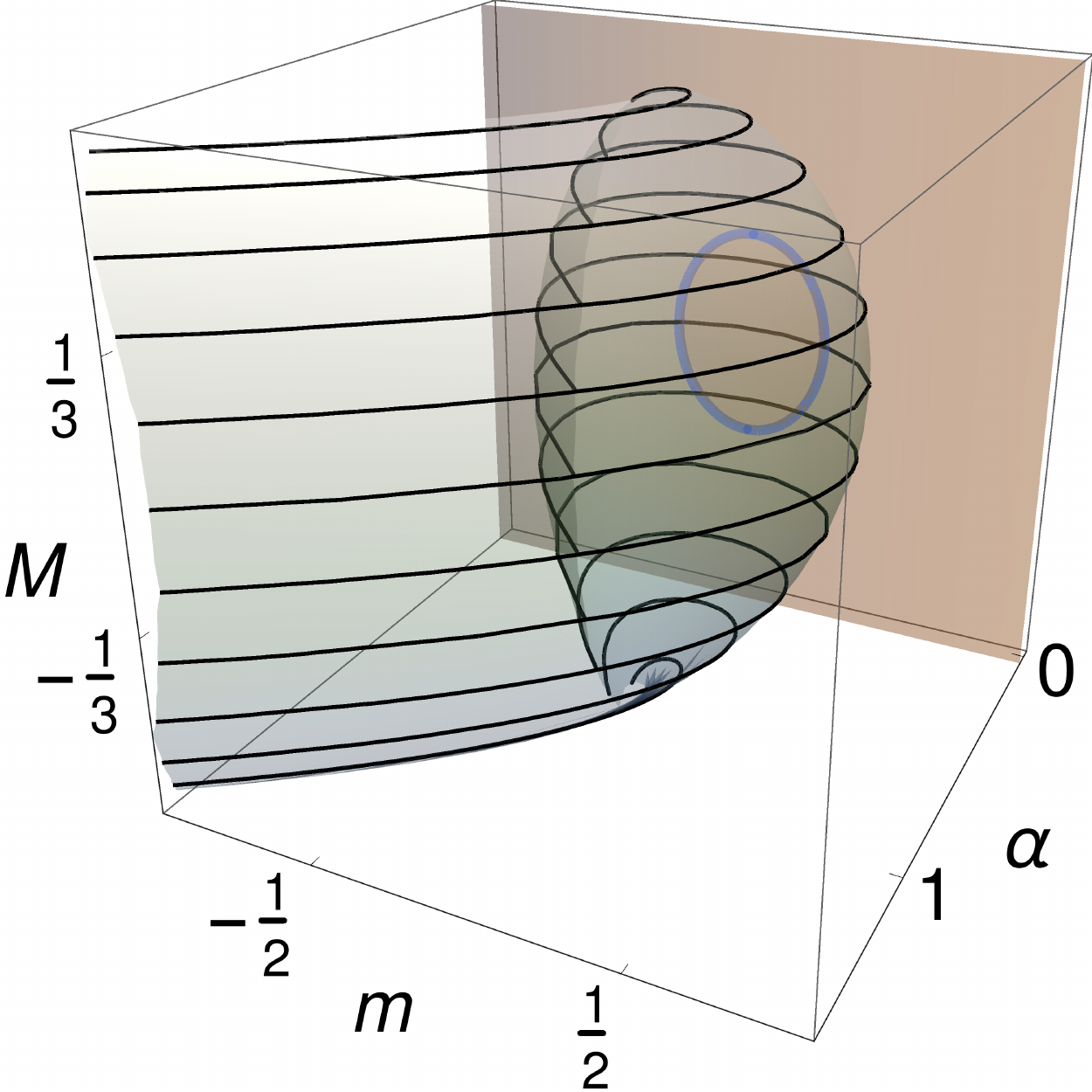}\hfill
\includegraphics[width=0.27\textwidth]{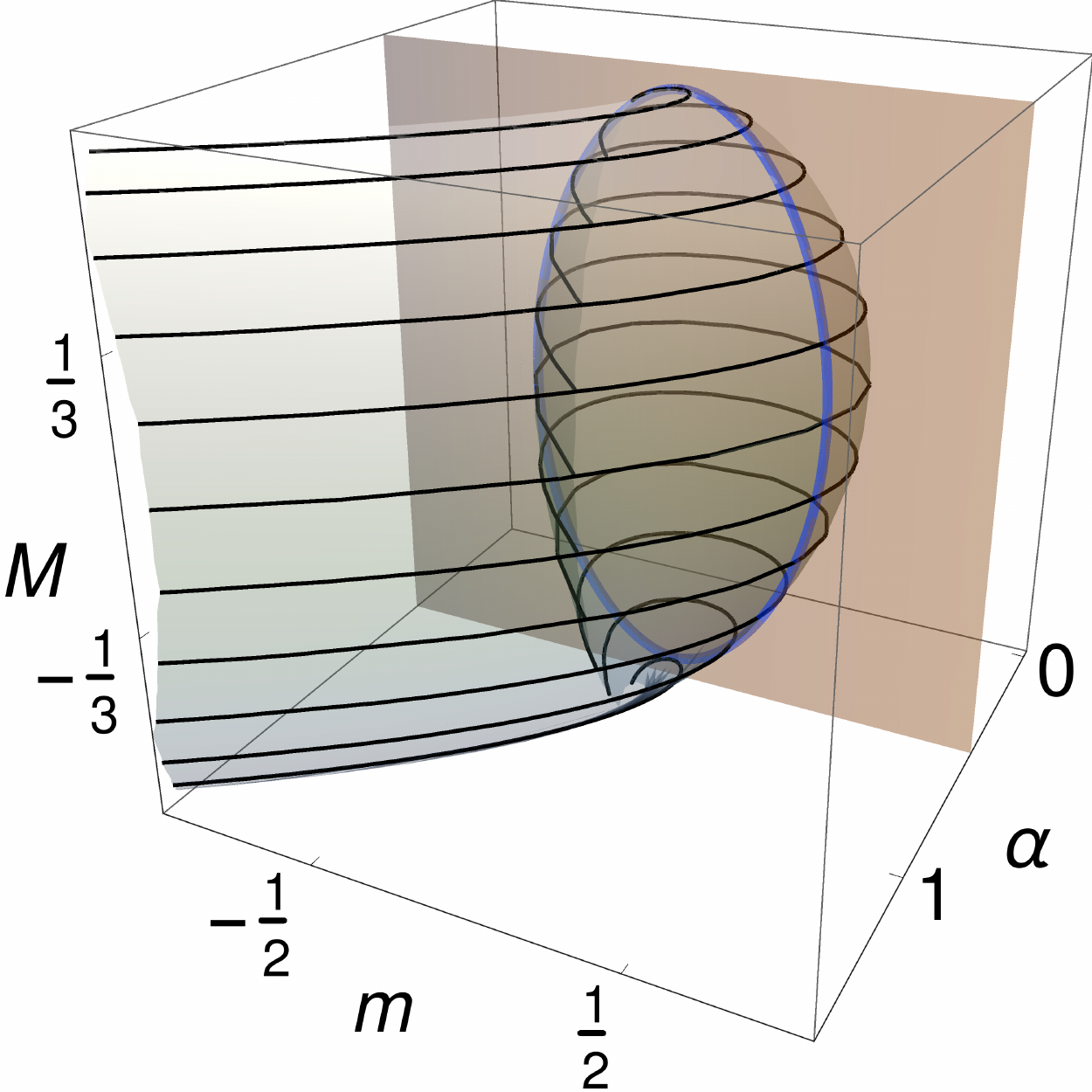}\hfill
\includegraphics[width=0.27\textwidth]{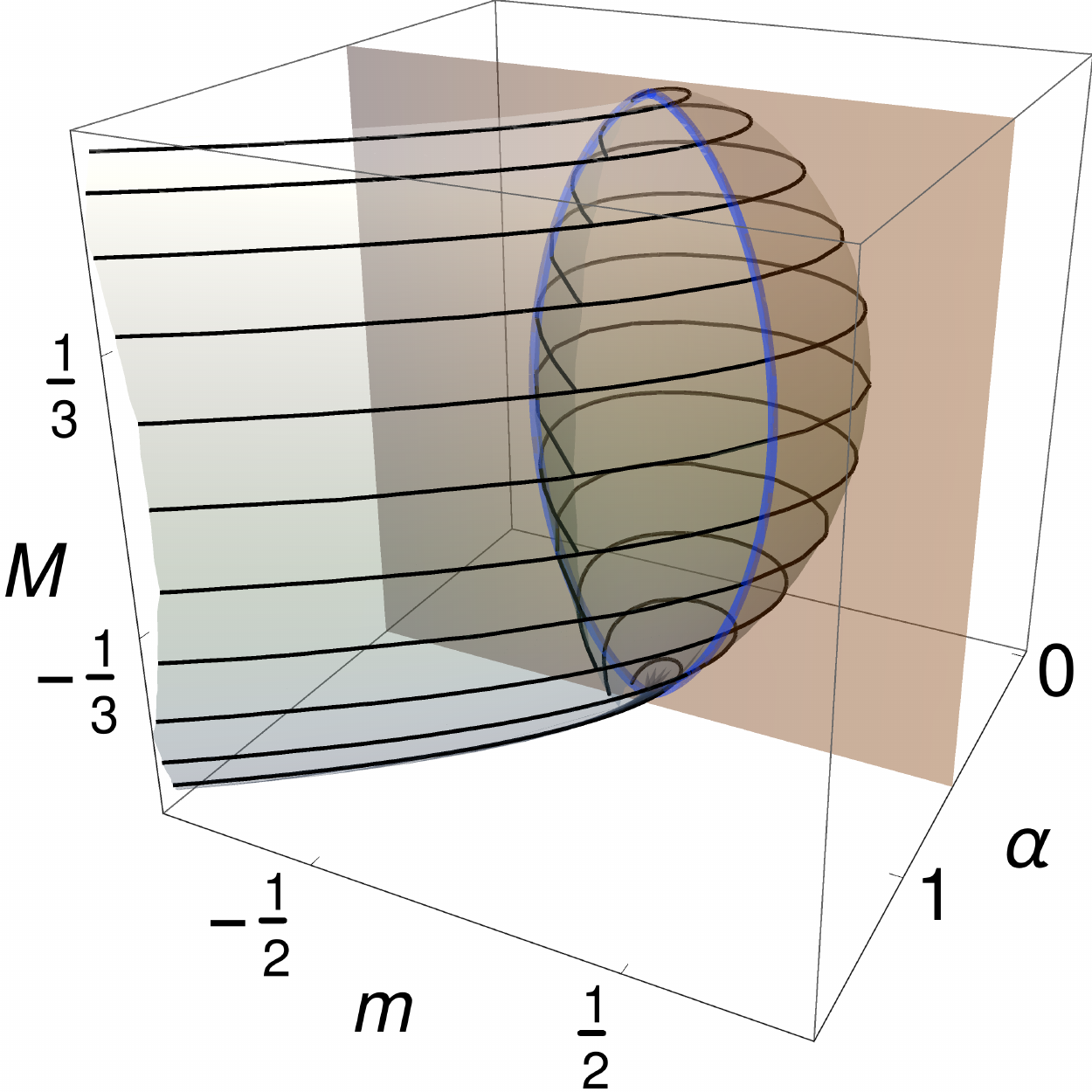}\\[0.2in]
\includegraphics[width=0.27\textwidth]{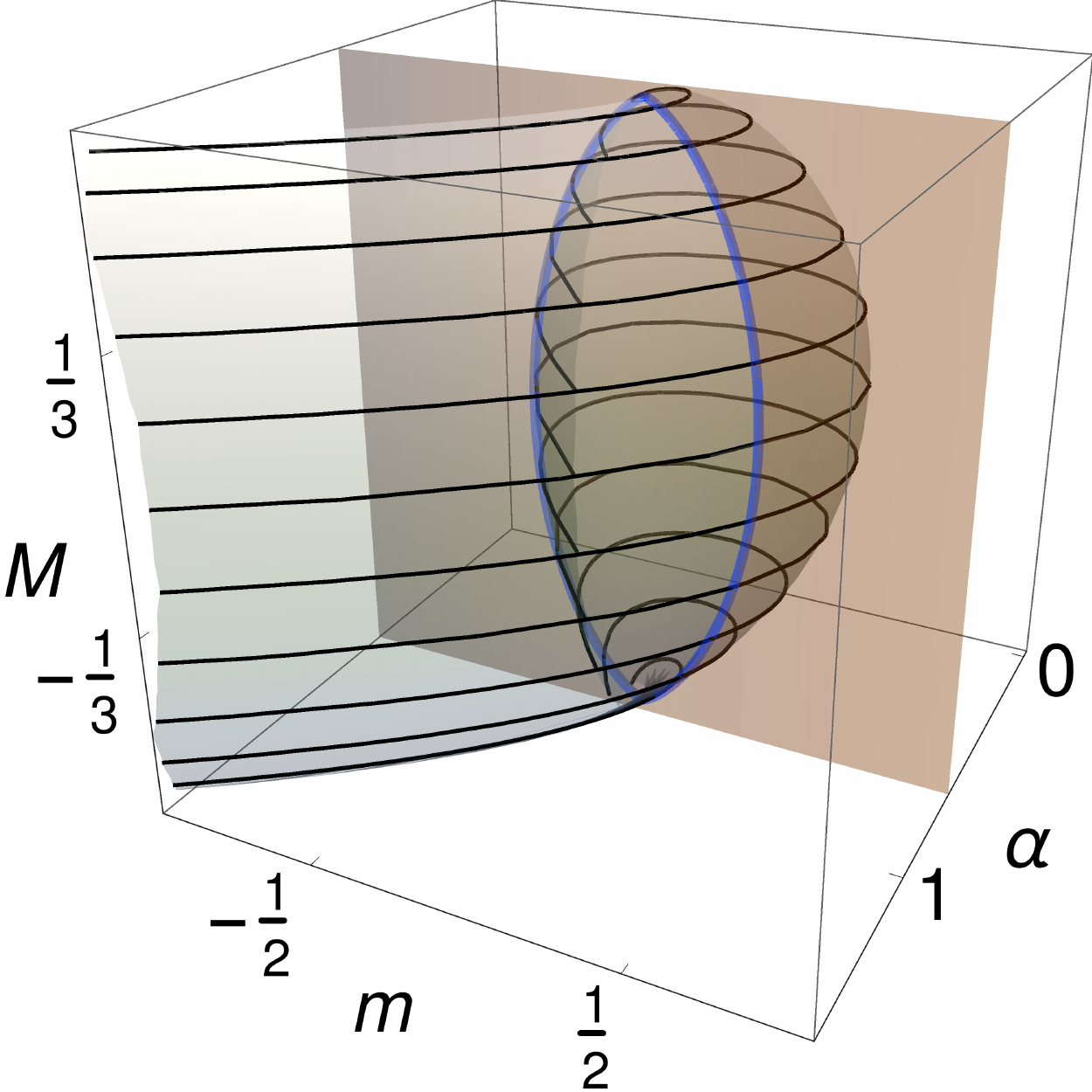}\hfill
\includegraphics[width=0.27\textwidth]{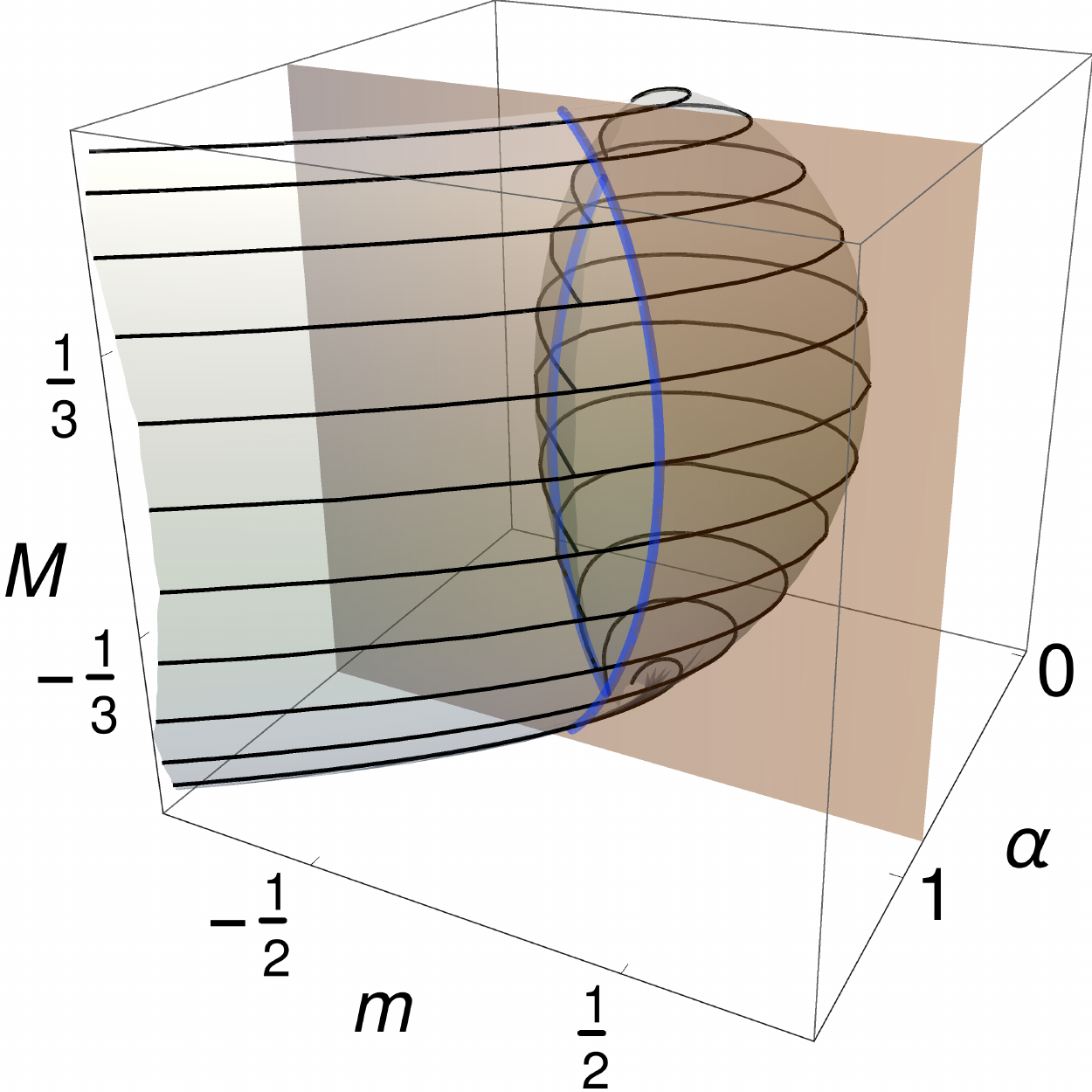}\hfill
\includegraphics[width=0.27\textwidth]{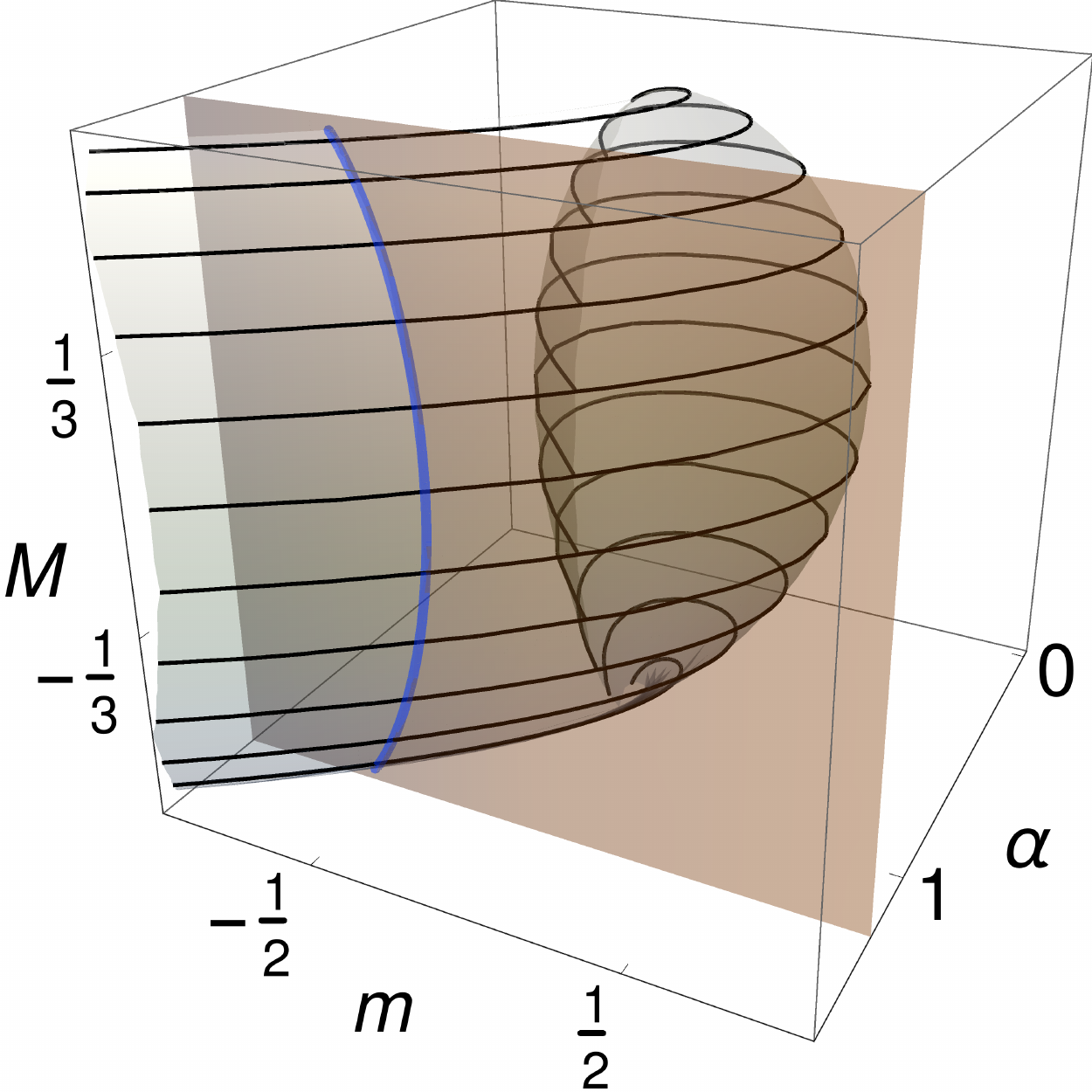}
\caption{Stable manifold in $(\alpha,m,M)$--space, together with the Dirichlet plane and the intersection. 
Parameter values are $\mu=0.2,0.7,0.8,0.82,0.93,1.1$.}\label{f:dirtau-fish}
\end{figure}
\begin{figure}[h!]
\includegraphics[width=0.27\textwidth]{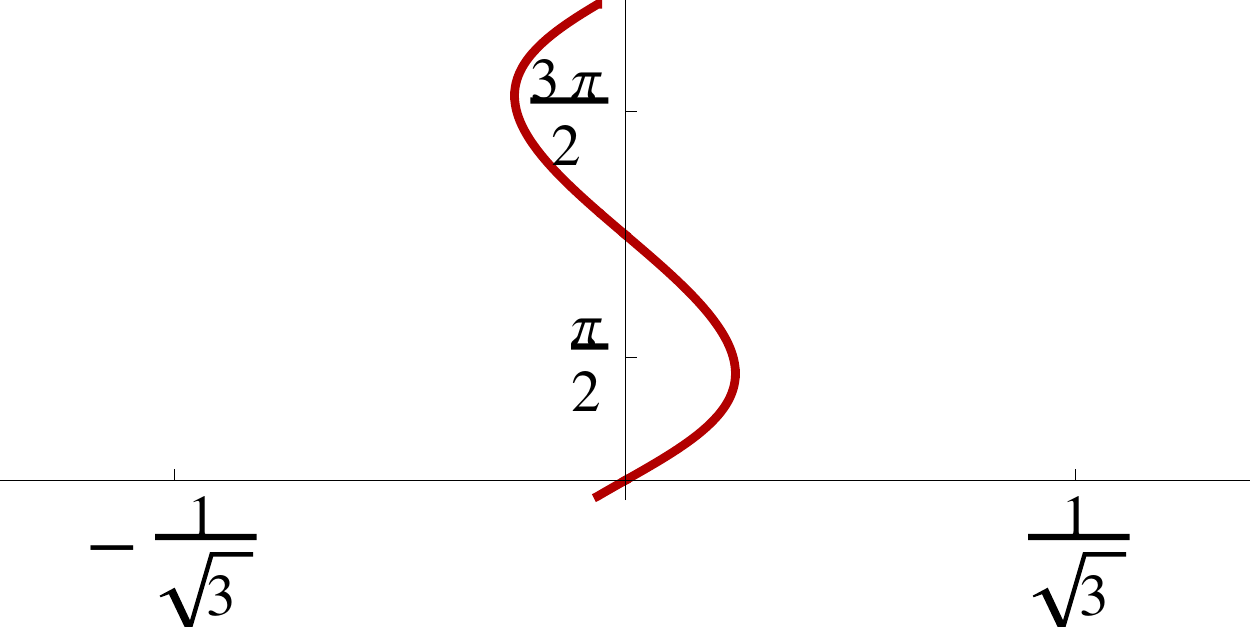}\hfill
\includegraphics[width=0.27\textwidth]{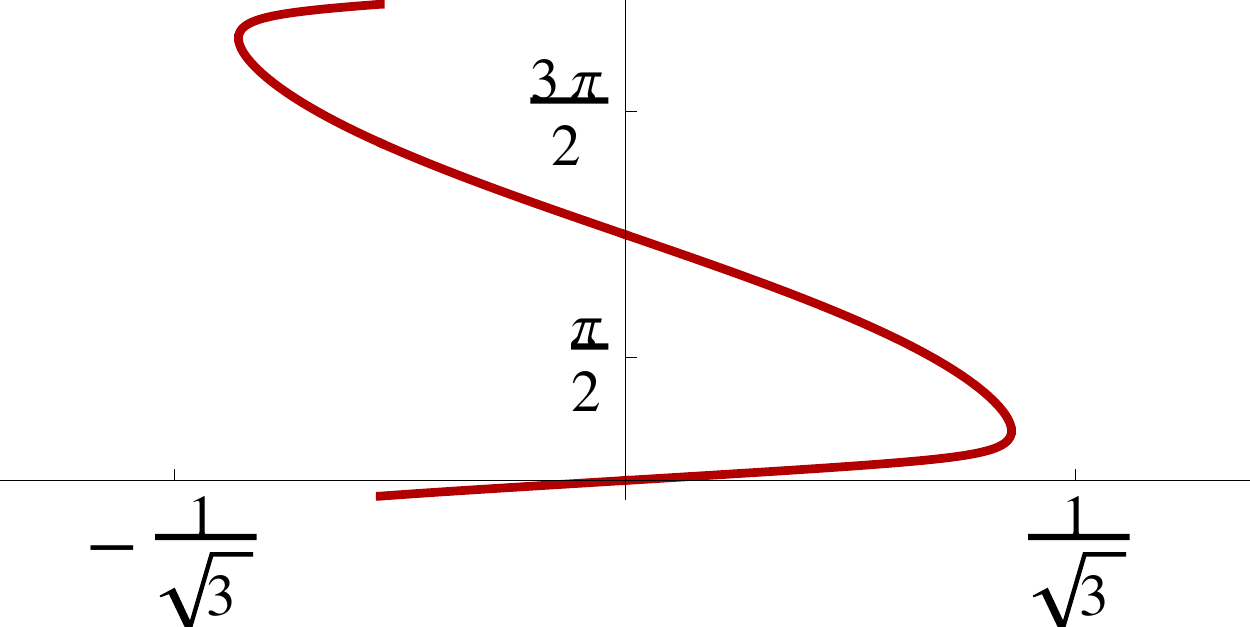}\hfill
\includegraphics[width=0.27\textwidth]{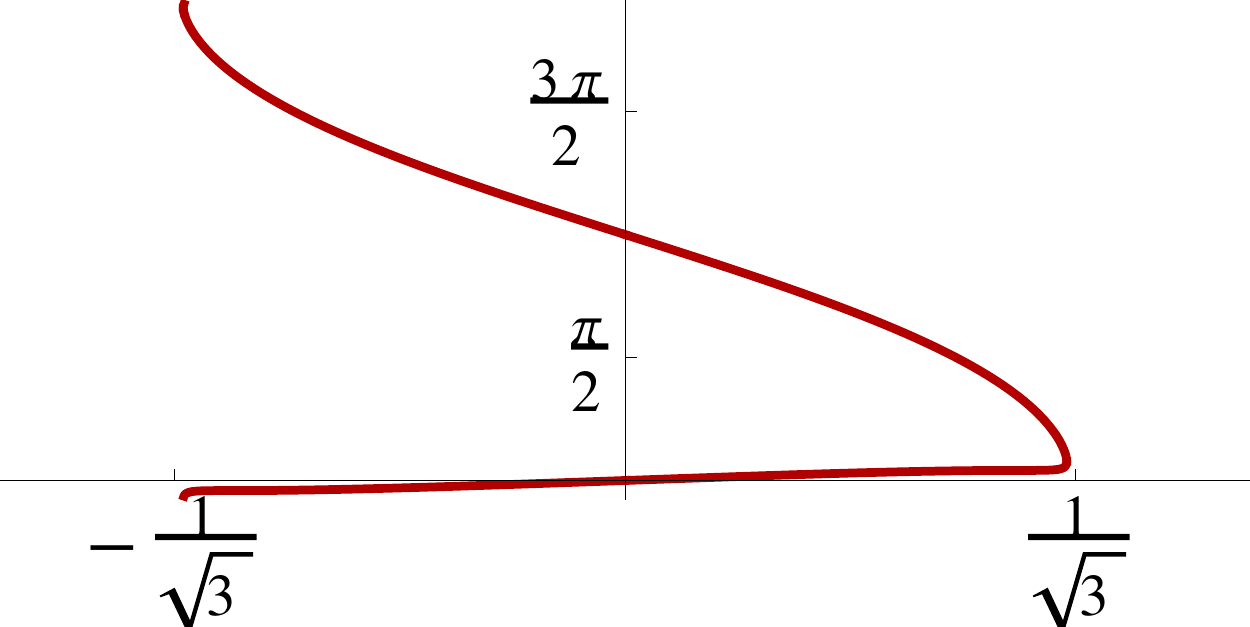}\\[0.2in]
\includegraphics[width=0.27\textwidth]{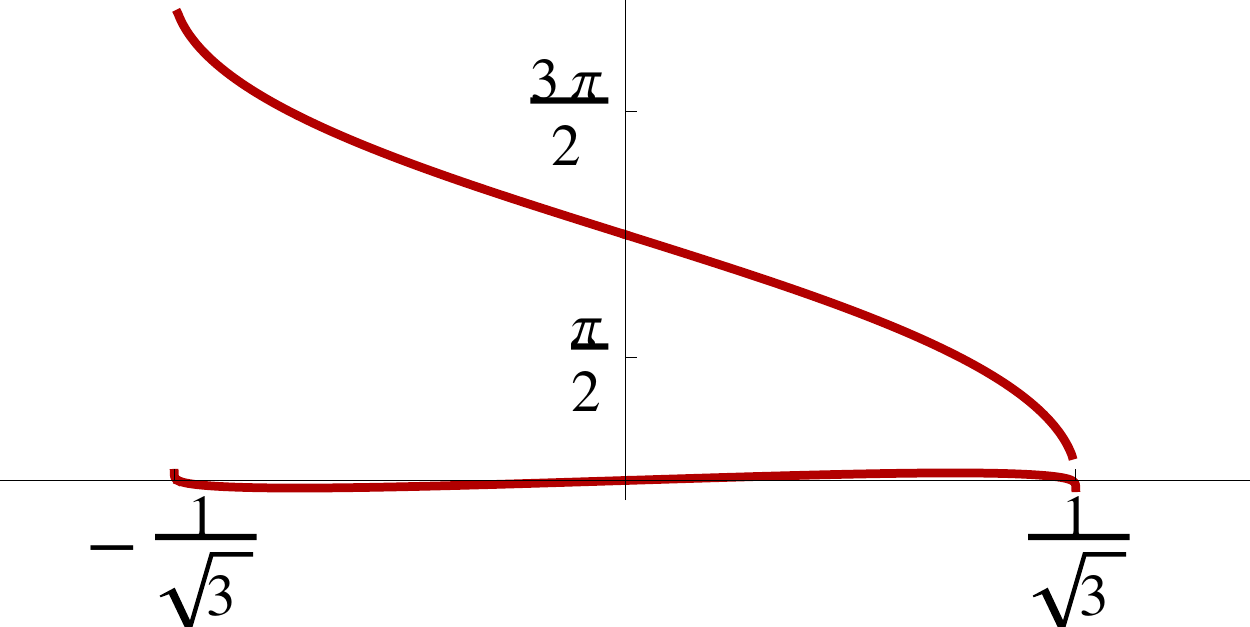}\hfill
\includegraphics[width=0.27\textwidth]{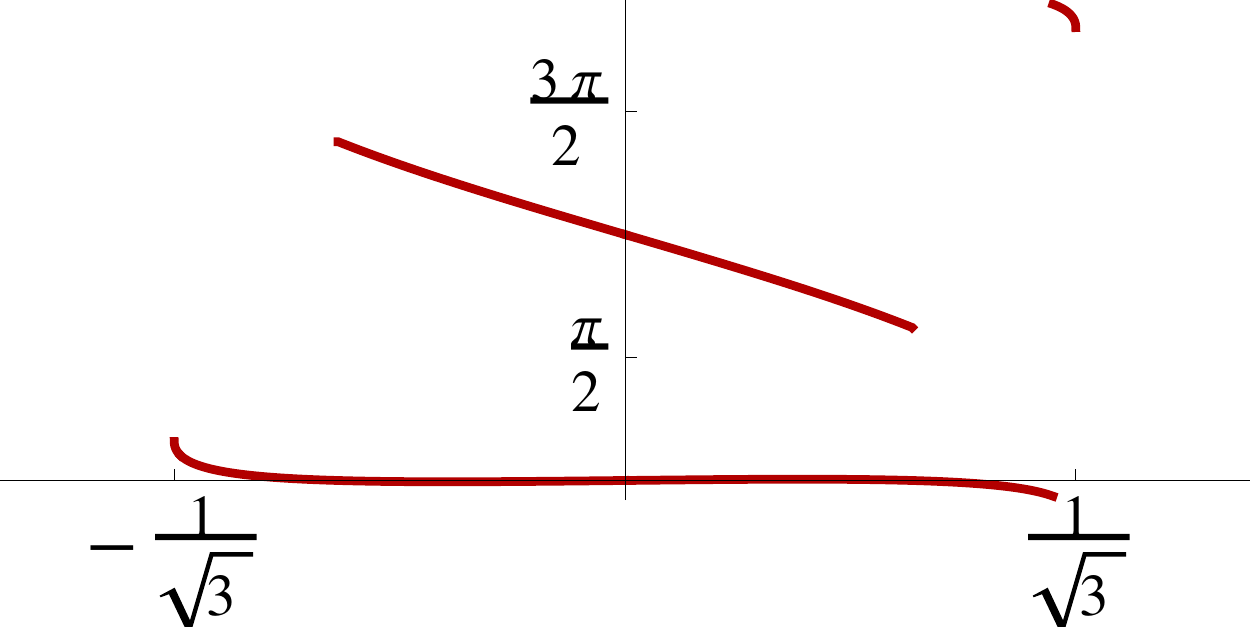}\hfill
\includegraphics[width=0.27\textwidth]{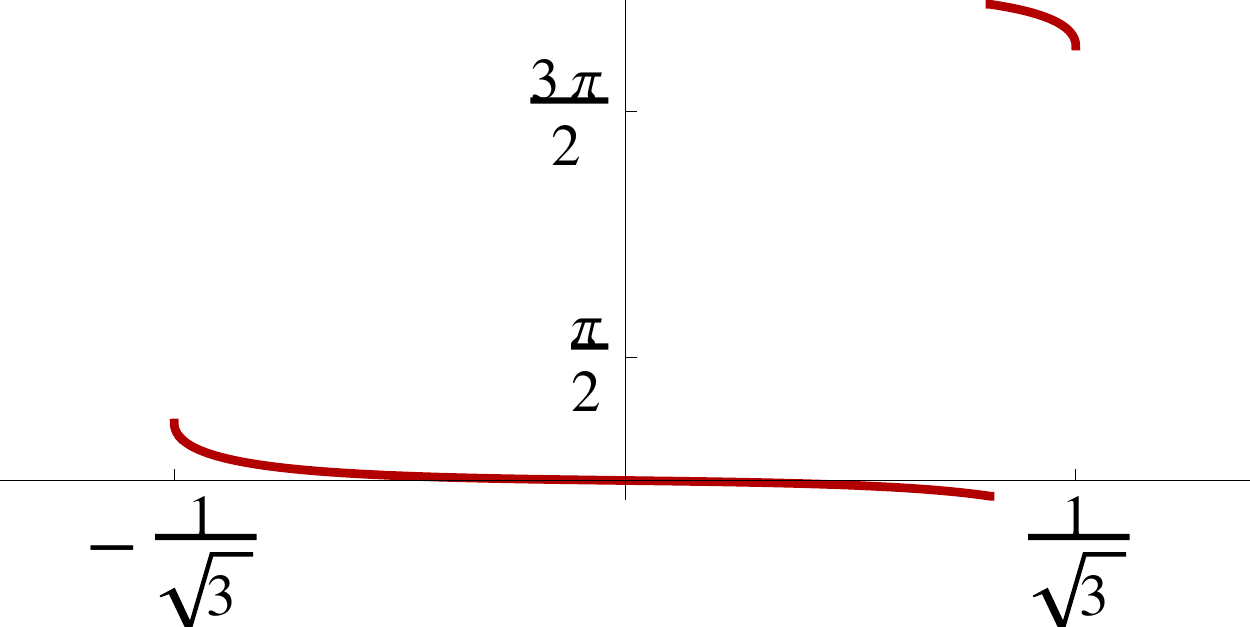}
\caption{Plots of displacement-strain curves in the  $(\varphi,k)$-plane, cporresponding to Figure \ref{f:dirtau-fish}. We first see the separation crisis at the Eckhaus boundary, followed by the vanishing of two defect limits. Eventually, $\varphi$ is not monotone.}\label{f:dirk-phi}
\end{figure}

\begin{Remark}\label{r:trans}
All boundary layers constructed in this section are in fact transverse in the sense of Definition \ref{d:bl}, provided that the intersection in the reduced orbit space is transverse. An example of a 
non-transverse boundary layer is the half-defect $A_\mathrm{d}=\rmi\tanh(x/\sqrt{2})q
$, $x>0$, with boundary conditions, $A_x=\rmi/\sqrt{2}$. 
Inspecting the formulas for the $k$-dependent defects (\ref{e:tanh}) and differentiating at the origin $x=0$, we find tangent vectors $(A,B)=(0,-1)$ from the phase, $(A,B)=(\rmi,0)$ from $x$-translation, and $(A,B)=(2,0)$ from the $k$-derivative, so that the real span of $(0,\rmi)$ is not contained in $TW^\mathrm{cs}+T\mathcal{B}$. The reason is of course that these boundary conditions mimic infinitesimal translation of the defect at the origin. 
\end{Remark}

\subsection{Existence --- linear and affine boundary conditions}

Linear boundary conditions translate into a cone in $\alpha,m,M$-space, by scaling invariance. On the other hand, the expansion of $W^\mathrm{s}$ at $\alpha=m=M=0$ is quadratic, so that any cone is contained in the interior of the body of the fish, which guarantees intersections between $W^\mathrm{s}$ and the cone, hence existence of boundary layers. 

\begin{Proposition}\label{p:ex}
For all affine boundary conditions
\[
T_1 A + T_2 B = \mu,\quad, 
\]
with $T_j:\R^2\to\R^2$ real linear, $\mu\in\R^2$, $\Rg T_1 + \Rg T_2=\R^2$, there exists a boundary layer. 
\end{Proposition}
\begin{Proof}
The proof exploits the topology of $(\alpha,m,M)$--space. We can distinguish three different cases depending on the rank of $T_2$:
\begin{itemize}
\item[(1)] $\mathrm{rk}\,T_2=2$;
\item[(2)] $\mathrm{rk}\,T_2=1$;
\item[(3)] $\mathrm{rk}\,T_2=0$.
\end{itemize}
In the following, we will treat various scenarios on a case-by-case basis. 
\paragraph{Case (3).} The map $T_1$ is invertible, so that the boundary conditions are equivalent to Dirichlet boundary conditions, and we already established existence in the previous section.

\paragraph{Case (1).} We can invert $T_2$ and obtain a boundary condition of the form 
\[ 
B=T A + \nu,\qquad T=T_2^{-1}T_1, \ \nu=T_2^{-1}\mu,  
\]
or $B=\kappa_1 A+\kappa_2 \bar{A}+\nu$. In case $\kappa_1=\kappa_2=0$, $T_1=0$, we recover the Neumann boundary conditions studied in the previous section and we can conclude existence. We consider the boundary manifold in the space of invariants $\alpha,m,M$, thinking of $m$ as the vertical coordinate. We set $A=\xi=r\rme^{\rmi\psi}$, which gives 
\[
\alpha=r^2,\qquad \mathcal{M}=r^2 \left(\bar{\kappa}_1+\bar{\kappa}_2\rme^{2\rmi\psi}\right) +r\nu\rme^{\rmi\psi}.
\]
We view this as a family in $r$ of closed curves, each parameterized by $\psi$, and compare the location of the curve with the position of the stable manifold. For large $r$, the stable manifold is located near $m\sim -\alpha^2$, so that the curves $\mathcal{M}\sim r^2 \left(\bar{\kappa}_1+\bar{\kappa}_2\rme^{2\rmi\psi}\right)$ 
lie above the stable manifold. For small $r$, suppose first that $\nu>0$. The curve $\mathcal{M}\sim r\nu\rme^{\rmi\psi}$ then is located either inside the compact part of the stable manifold, or it winds around or intersects it. In either case, homotoping between small and large $r$, we see that there necessarily is an intersection. 
In the case $\mu=0$, linear boundary conditions, the small-$r$ asymptotics are equal to the large-$r$ case and one readily sees that the curve is located inside the compact part of the stable manifold, which has a quadratic tangency with the plane $\alpha=0$. Again, a homotopy to large $r$ gives intersections between our family of curves and the stable manifold.

\paragraph{Case (2).} Without loss of generally, possibly rotating with the gauge symmetry, the boundary conditions can be written in the form
\[
\Re B = \kappa_1 \Re A + \kappa_2 \Im A + \mu_1,
\quad \Im A = \sigma_1 \Re A + \mu_2,
\]
with real parameters $\mu_1,\mu_2, \kappa_1,\kappa_2,\sigma_1$. Parameterizing $\Im B = y$, $\Re A = x$, we find 
\[
\alpha = x^2+ (\sigma_1 x + \mu_2)^2,\quad \mathcal{M}=\left(x+\rmi(\sigma_1 x + \mu_2)\right)\left(\kappa_1 x + \kappa_2 \sigma_1 x + \mu_1 + \kappa_2 \mu_2 + \rmi y \right).
\]
We view this surface as a family in $x$ of straight lines parameterized by $y$. We follow the family of lines from $x=+\infty$ to $x=-\infty$, and note that  near the endpoints of the homotopy, the lines are close, with tangent close to $(\rmi - \sigma_1) x^2$. Since $\mathcal{M}=\rmO(\alpha)$ when $y=\rmO(x)$, the line is located above the stable manifold for large $|x|$ (and large $\alpha$), with non-real tangent.  As $x$ is increased from $-\infty$ to $+\infty$, the tangent vector rotates by $\pi$. As a consequence, the family of lines  necessarily intersects the stable manifold. 
\end{Proof}

In the case where the equation comes from a variational problem, 
\[
\mathcal{E}[A]=\int_0^\infty \left(|A_x|^2+\frac{1}{2}(|A|^2-1)^2\right)\rmd x + \left.\left(\mu|A|^2+(\nu \bar{A}^2+c.c.)+(\rho\bar{A}+c.c.)\right)\right|_{x=0},\qquad \mu\in\R,\nu,\rho\in\C,
\]
one obtains boundary conditions 
\begin{equation}\label{e:symplbc}
B=\mu A + \nu \bar{A}+\rho,\quad \mu\in\R,\nu,\rho\in\C,
\end{equation}
which form a Lagrangian (affine) subspace of $\C^2$, with respect to the symplectic form \[\omega((A_1,B_1),(A_2,B_2))=\Re(A_1\bar{B}_2-A_2\bar{B_1}).\]
Since in the interior of the domain, $|A|=1$ and the associated solutions $A=\rme^{\rmi\varphi}$ are the only minimizers, one expects boundary layers with finite energy, thus connecting to $k=0$. While such a result can presumably be established using variational techniques, possibly for more general systems, we restrict here to the Ginzburg-Landau equation and the geometric viewpoint from Proposition \ref{p:ex}. 

\begin{Proposition}\label{p:exsympl}
For all Lagrangian boundary conditions (\ref{e:symplbc}), there exists a boundary layer with $k=0$.
\end{Proposition}
\begin{Proof}
We describe the boundary surface in the invariant space setting 
$A=\xi=r\rme^{\rmi\psi}$, so that
\[
\alpha=r^2,\qquad \mathcal{M}=(\mu+\bar{\nu}\rme^{2\rmi\psi})r^2+\bar{\rho}\rme^{\rmi\psi}r.\]
For $\bar{\rho}=0$, these boundary conditions describe a cone with elliptical cross-section, which necessarily intersect the compact part of the stable manifold in the plane $M=0$. For $\bar{\rho}\neq 0$, we can view the boundary manifold as a family in $r$ of curves parameterized by $\psi$. One can choose two intersection of these curves with the real axis (maximum and minimum) depending on $r$ in a continuous fashion (the curves are Lima\c{c}ons). For $r$ large, this intersection point lies above $W^\mathrm{s}$. For $r$ small, the intersection points lie either above and below, or inside, or on the compact part of the stable manifold. For large $r$, the intersection points lie above the stable manifold. As a consequence, there is an intersection, which concludes the proof. 
\end{Proof}

\begin{Remark}[Non-existence and isolas]\label{r:ne}
Boundary layers may not exist even with well-posed boundary conditions. A simple and instructive example is the cubic-quintic Ginzburg-Landau equation 
\[
A_t=A_{xx}-A+\gamma A|A|^2-A|A|^5.
\]
Periodic orbits $r\rme^{\rmi kx}$ exist for $k^2=-1+\gamma r^2-r^4$, provided $\gamma\geq 2$. The steady-state equation preserves the Hamiltonian $H=|A_x|^2-|A|^2+\frac{\gamma}{2}|A|^4-\frac{1}{3}|A|^6$. For $\gamma\sim 2$, all periodic orbits are close to $|A|=1$ and have $H\sim -1/3$. Choosing ``asymptotic boundary conditions'', $\mathcal{B}= W^\mathrm{u}(A=0,B=0)$, we conclude that all boundary layers would have Hamiltonian $H=0$, which is incompatible with the Hamiltonian of periodic orbits. The unstable manifold $W^\mathrm{u}$ can in fact be expressed as a smooth graph $B=A\sqrt{1-\frac{\gamma}{2}|A|^2+\frac{1}{3}|A|^4}$, so that the equation with these nonlinear boundary conditions is in fact locally well-posed. Alternatively, imposing Dirichlet boundary conditions $A=0$ gives $H\geq 0$ on the boundary, which also precludes intersections with the stable manifold of periodic patterns, where $H\sim -1/3$.

Considering inhomogeneous Dirichlet boundary conditions $A=\mu$ with increasing $\mu$, one finds a tangency between boundary conditions and stable manifold (much like the tangency in the cubic case at $\mu=0$, just at some finite $\alpha>0$, $k=0$), which then leads to a small closed curve of intersections and a small isola in the $(\varphi,k)$-plane.
\end{Remark}

\subsection{Scaling boundary conditions}\label{s:3.4}

In general, Ginzburg-Landau spatial dynamics are recovered on a 4-dimensional center manifold for spatial dynamics \cite{radial}.  Boundary conditions will typically intersect the center-stable manifold along a two-dimensional manifold,which can be collapsed along the stable fibration to obtain effective 2-dimensional boundary manifolds in the center manifold. In the center manifold, the  dynamics are, up to  changes of coordinates and higher-order terms, 
\[
A_x=\rmi k A + B, \quad B_x=\rmi k B - \mu A +A|A|^2,
\]
where $\mu$ is the bifurcation parameter. The usual scaling $A=\sqrt{\mu}\tilde{A}$, $B=\mu\tilde{B}$, $k=\sqrt{\mu}\tilde{k}$, $y=\sqrt{\mu}x$, eliminates the parameter $\mu$ and gives, up to the fast rotation, the Ginzburg-Landau equation. Under the scaling, most linear boundary conditions will reduce to homogeneous Dirichlet boundary conditions at leading order since derivatives are small in the expansion. 

\section{Computation of boundary layers in the Swift-Hohenberg equation}\label{s:4}

Computing boundary layers is related to the computation of heteroclinic orbits. Numerically one needs to approximate the infinite domain $[0,\infty)$  by a finite domain $[0,L]$ and discretize in space. The first step, truncation, involves choosing appropriate boundary conditions at $x=L$. The correct boundary condition, $\underline{u}(L)\in W^\mathrm{s}$ is usually approximated by its linearization. A difficulty in our case is the presence of neutral directions, phase and wavenumber, and the fact that the asymptotic state is not explicitly known. 

In the following we present the general strategy of our approach, Section \ref{s:4.1} and some results for a specific set of boundary conditions, Section \ref{s:4.3}. 

\subsection{Decomposing periodic and heteroclinic orbit}\label{s:4.1}
Let $\chi(x)$ be a smooth function with $\chi(x)=0$ on $x<\ell_-$, $\chi(x)=1$ on $x>\ell_+$, for some $0<\ell_-<\ell_+<L$. Also, denote by $\mathcal{L}_\mu=-(\partial_{xx}+1)^2+\mu $ the linear part of the Swift-Hohenberg equation. Lastly, we write $[\mathcal{L}_\mu,\chi]u:=\mathcal{L}_\mu(\chi u)-\chi \mathcal{L}_\mu u$. 

With these preparations, we substitute the ansatz $u(x)=\chi(x)u_\mathrm{st}(kx-\varphi;k)+w(x)$ into the Swift-Hohenberg equation, subtract $\chi(Lu_\mathrm{st}-u_\mathrm{st}^3)$, and obtain
\begin{align}
\mathcal{\tilde{F}}_0(w,\varphi,k)=\mathcal{L}_\mu w - \left((\chi u_\mathrm{st} +w)^3-(\chi u_\mathrm{st})^3\right)+ [\mathcal{L}_\mu,\chi]u_\mathrm{st} -(\chi-\chi^3)u_\mathrm{st}^3&=0,\nonumber \\
(w,w_x,w_{xx},w_{xxx})(0)&\in \mathcal{B}.\label{e:tF}
\end{align}
Since $\chi$ vanishes near $x=0$, $w$ needs to satisfy the same boundary conditions that we imposed on $u$. The first two terms vanish for $w=0$, reflecting the fact that $u_\mathrm{st}$ solves the Swift-Hohenberg equation.  The last two terms, which we refer to as residuals,  are compactly supported on $[\ell_-,\ell_+]$. 

We are interested in solutions $w$ that are exponentially localized. Define therefore $L^2_\eta=\{u; u(x)\rme^{\eta x}\in L^2\}$, with induced norm, and, analogously $H^k_\eta=\{u;\partial_x^\ell u\in L^2_\eta, \ell\leq k\}$.  It turns out that $\tilde{\mathcal{F}}_0$ is not differentiable in with respect to $k$ as a map: the derivative with respect to $k$ of the term $3(\chi u_\mathrm{st})^2 w$ does in general not belong to $L^2_\eta$ due to the linear multiplier $6\chi^2 u_\mathrm{st}x\partial_x u_\mathrm{st} $. A simple rescaling can remedy this difficulty as we shall see below. Consider the scaled differential operator $\mathcal{L}_{\mu,k}=-(k^2\partial_{yy}+1)^2+\mu$, with appropriately scaled boundary conditions, and define 
\begin{align}
\mathcal{F}_0(w,\varphi,k)=\mathcal{L}_{\mu,k} w - \left((\chi u_\mathrm{st} +w)^3-(\chi u_\mathrm{st})^3\right)+ [\mathcal{L}_{\mu,k},\chi]u_\mathrm{st} -(\chi-\chi^3)u_\mathrm{st}^3&=0, \nonumber\\
(w,kw_y,k^2w_{yy},k^3w_{yyy})(0)&\in \mathcal{B},\label{e:F}
\end{align}
where $u_\mathrm{st}=u_\mathrm{st}(y-\varphi;k)$ is the $2\pi$-periodic, scaled family of stripes, $\chi=\chi(y)$ can be left unscaled. One readily finds that solutions $w,k,\varphi$ to (\ref{e:F}) are in one-to-one correspondence with solutions $\tilde{w},k,\varphi$ to (\ref{e:tF}) via the rescaling $\tilde{w}(x)=w(kx)$ and the appropriate modification of $\chi$.

\begin{Lemma}[Well-Posedness]\label{l:well}
The map $\mathcal{F}_0:H^4_\eta\cap \mathcal{B}\times \R^2\to L^2_\eta$ is well-defined and smooth near a given boundary layer $w=u_\mathrm{bl}-\chi u_\mathrm{st}$, for $\eta>0$, sufficiently small. Moreover, the linearization at a boundary layer is Fredholm with index 1. If in addition the boundary layer is transverse, according to Definition \ref{d:bl}, then $\mathcal{F}_0'$ is onto and the $k$- and $\varphi$-components of the one-dimensional kernel do not both vanish. 
\end{Lemma}
\begin{Proof}
The map is well-defined since residuals are compactly supported, thus belong to exponentially weighted spaces. The linearization with respect to $w$ is Fredholm with index -1 as a simple counting of Floquet exponents shows. Adding two parameters increases the Fredholm index to 1. The map and derivatives with respect to $w$ depend smoothly on $k$ and $\varphi$ as coefficients $u_\mathrm{st}$ depend smoothly on $k$ and $\varphi$ in $L^\infty$ and the dependence of the linear operator on $k$ is smooth. The last claim is an immediate consequence of the assumption of the intersection between boundary manifold and center-stable manifold in Definition \ref{d:bl}.
\end{Proof}
The Lemma can be viewed as a substitute for the smooth dependence of the strong-stable foliation  of the periodic on phase and wavenumber, parameterizing the base points of the fibration. As a direct consequence of the lemma, we obtain the existence of a one-parameter family of boundary layers, described by a smooth curve in the $(k,\varphi)$--plane. Using Lyapunov-Schmidt reduction.

Truncating to a finite domain would require $w(L)$ to lie in the (unknown) strong stable fiber of $u_\mathrm{st}$ (translated by $u_\mathrm{st}$). Since $w$ is exponentially small, this fiber is well approximated by its tangent space, so that a good asymptotic boundary condition would be $w(L)\in E^\mathrm{ss}(kL-\varphi;k)$, where $E^\mathrm{ss}(y)$ denotes the strong stable Floquet subspace of the linearization at $u_\mathrm{st}(y)$. Again, we find the computation of this Floquet subspace somewhat intricate. Alternatively, one can impose other boundary conditions at $x=L$, as long as the boundary subspace is transverse to $E^\mathrm{cu}$. \footnote{Such an approach is well understood in the context of heteroclinic or homoclinic orbits to hyperbolic equilibria. One finds that the convergence of solutions is $\rme^{-2\eta L}$ with correct linear Floquet boundary conditions and only $\rme^{-\eta L}$ with general transverse boundary conditions, where $\eta$ measures the spectral gap separating Floquet 
exponents on either side of the imaginary axis; see \cite[\S 6]{hettrunc} and references therein.} In practice, we choose Dirichlet boundary conditions and add a phase condition, to achieve the correct number of boundary conditions. As a phase condition, we choose 
\begin{equation}\label{e:phase}
\mathcal{F}_\mathrm{ph}(w,\varphi)=\int_K^{K+2\pi}u_\mathrm{st}'(y-\varphi)w(y)\rmd y,
\end{equation}
which ensures that $w$ cannot encode a phase shift in the periodic pattern (which is already encoded explicitly in the parameter $\varphi$). 

We obtain the following problem in a bounded domain:
\begin{align}
\mathcal{F}_0^L(w,\varphi,k)=\mathcal{L}_{\mu,k} w - \left((\chi u_\mathrm{st} +w)^3-(\chi u_\mathrm{st})^3\right)+ [\mathcal{L}_{\mu,k},\chi]u_\mathrm{st} -(\chi-\chi^3)u_\mathrm{st}^3&=0, \nonumber\\
(w,kw_y,k^2w_{yy},k^3w_{yyy})(0)&\in \mathcal{B}\nonumber\\
(w,w_{yy})(L)&=0\nonumber\\
\mathcal{F}_\mathrm{ph}(w,\varphi)=\int_K^{K+2\pi}u_\mathrm{st}'(y-\varphi)w(y)&=0.\label{e:FL}
\end{align}
One readily verifies that this problem is Fredholm index 1 and we therefore expect to find solution curves in the  $(k,\varphi)$-plane. We expect exponential convergence of the solutions to (\ref{e:FL}) to solutions of (\ref{e:F}), exponentially as $L\to\infty$. 

%
%

%

%

\subsection{Wavenumber selection in Swift-Hohenberg --- numerical results}\label{s:4.3}

We use second order finite differences to discretize the fourth-order  differential operator, using a mesh of width $h$ on a domain of length $L$. As part of the problem, we compute $u_\mathrm{st}(y;k)$ using Neumann boundary conditions on a domain of size $2\pi$. We then construct $u_\mathrm{st}(y-\varphi;k)$ using cubic splines and periodic extrapolation. We then solve (\ref{e:FL}) using arc length continuation within \textsc{Matlab}. We found the expected exponential convergence in the domain length, which breaks down near the Eckhaus boundary. We verified that results are independent of the cut-off function. 

\begin{figure}
\centering{
 \includegraphics[height=1.4in]{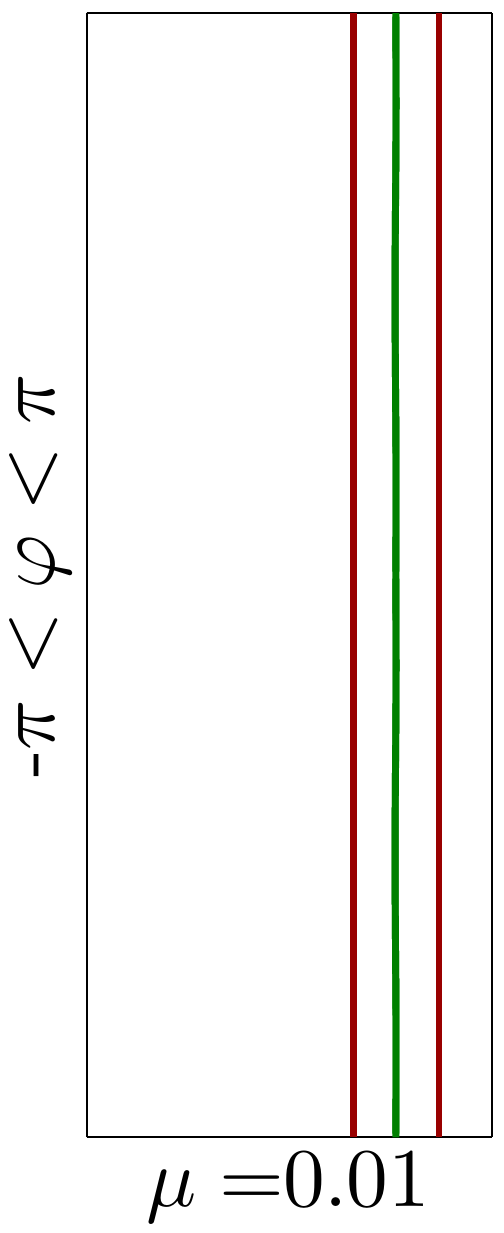}
\includegraphics[height=1.4in]{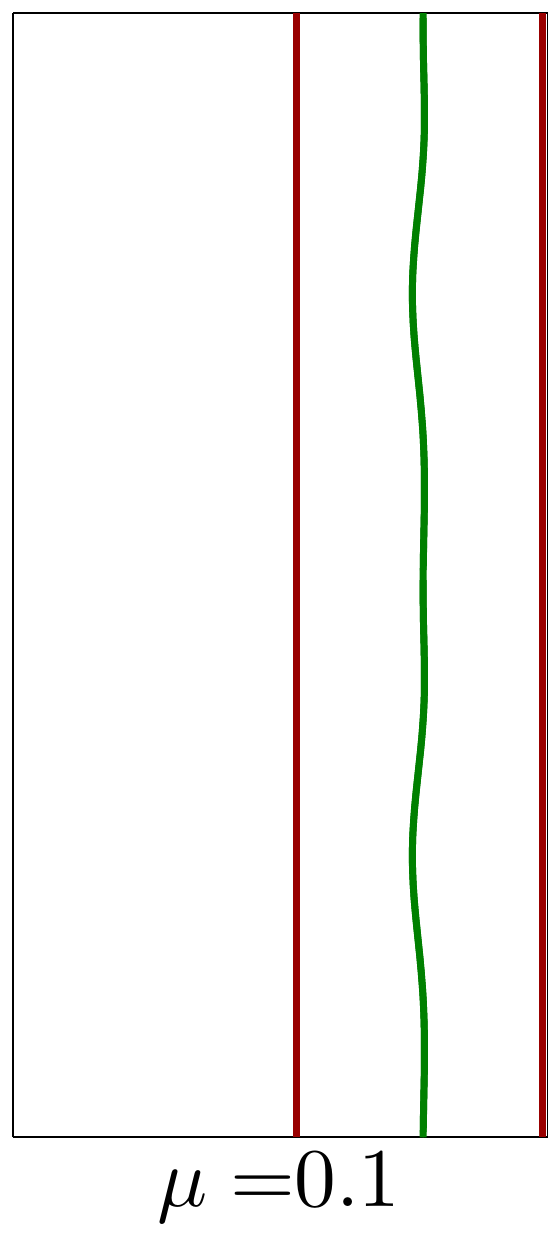}
\includegraphics[height=1.4in]{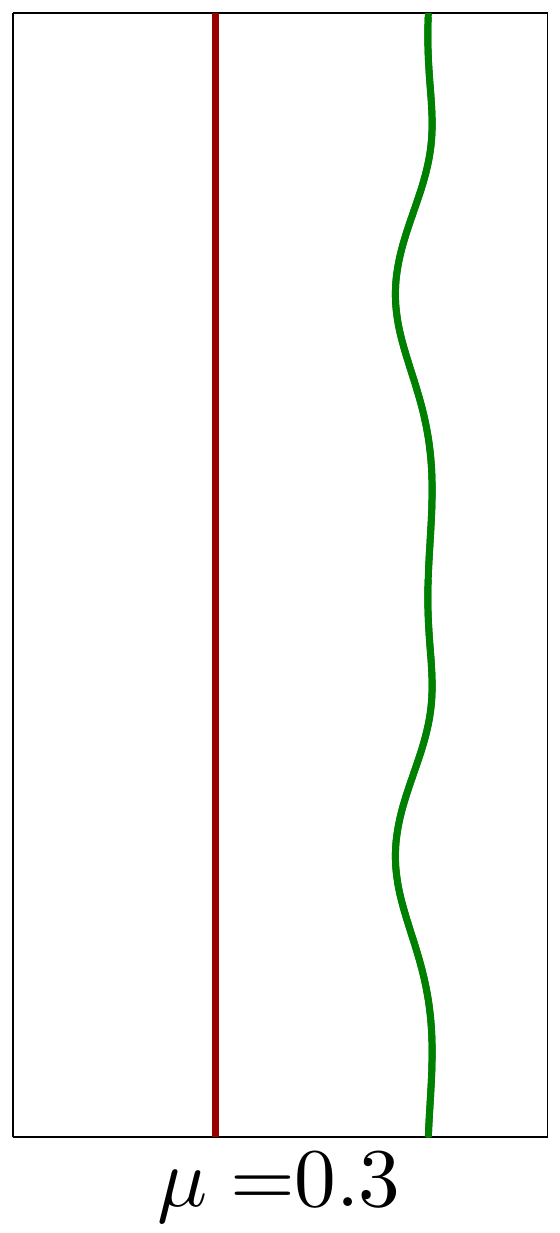}
\includegraphics[height=1.4in]{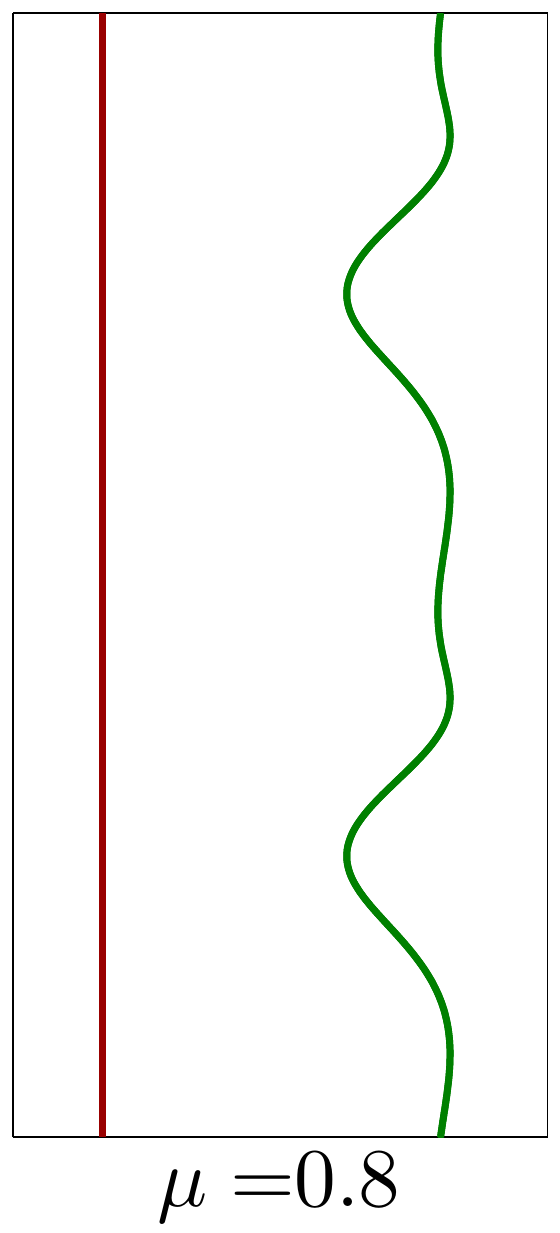}
\includegraphics[height=1.4in]{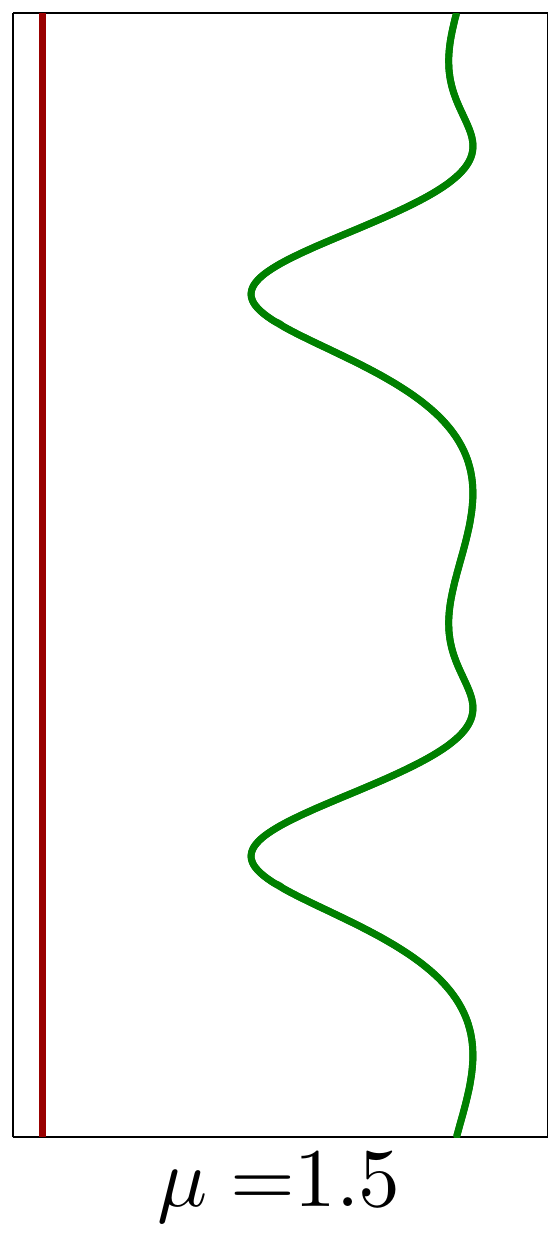}
\includegraphics[height=1.4in]{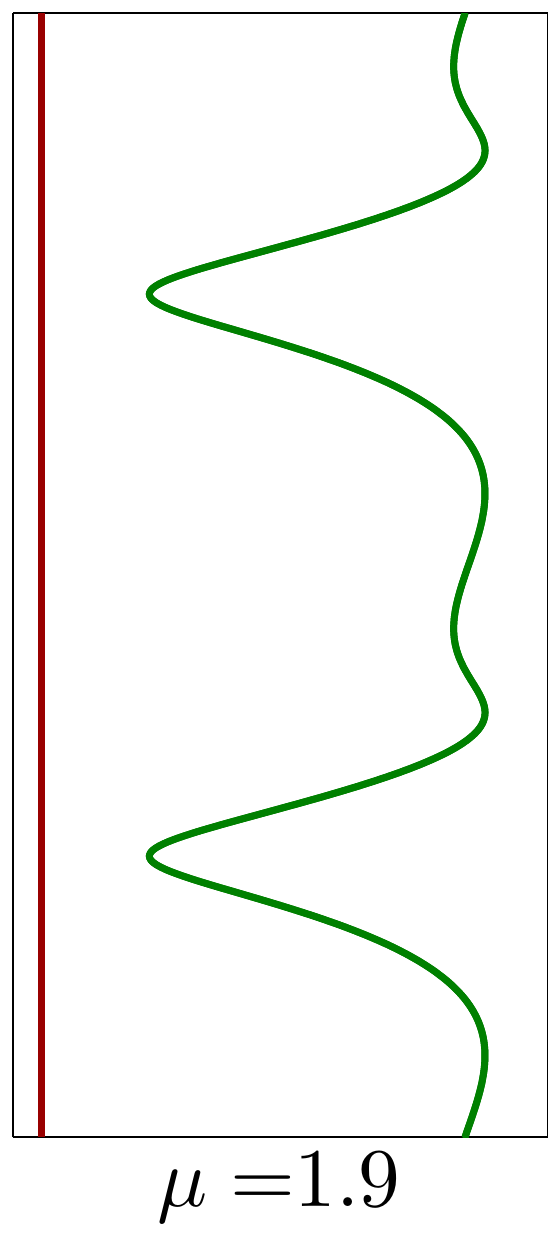}
\includegraphics[height=1.4in]{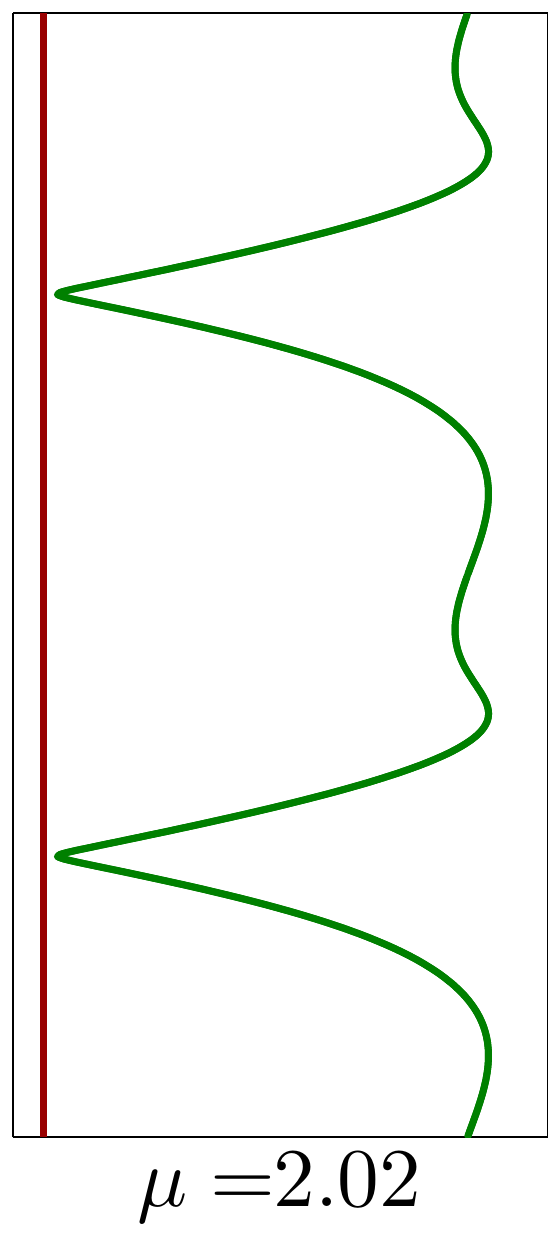}
\includegraphics[height=1.4in]{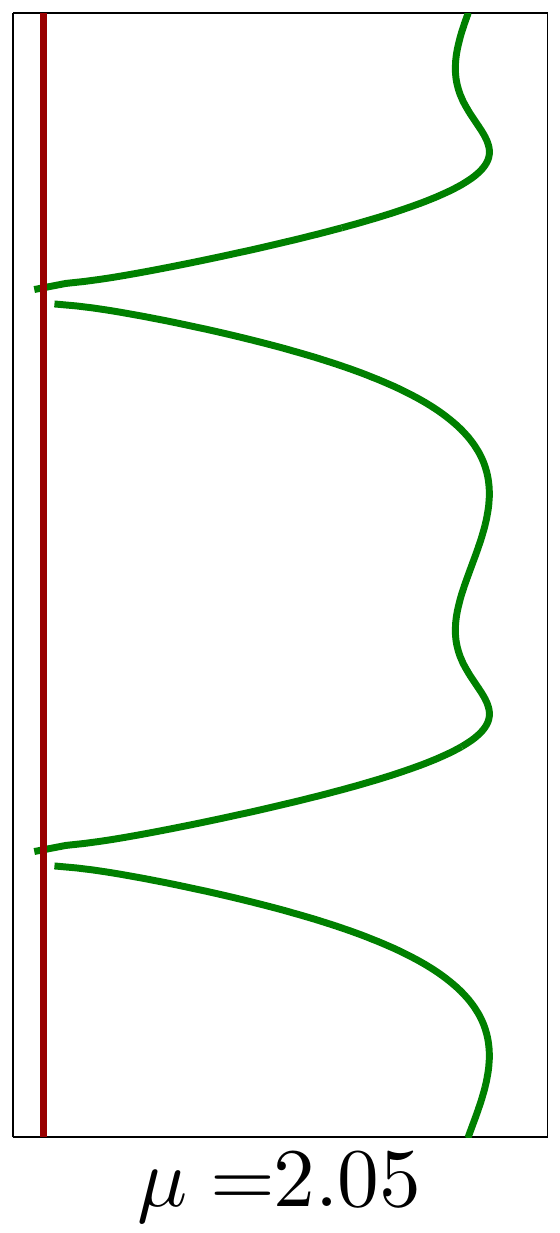}
\includegraphics[height=1.4in]{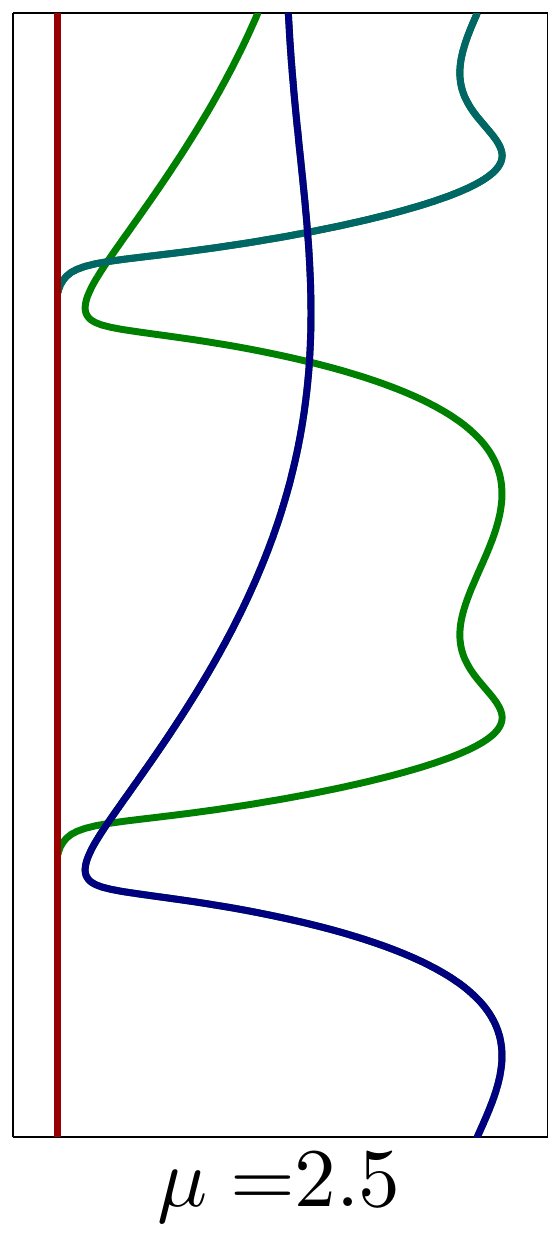}
\includegraphics[height=1.4in]{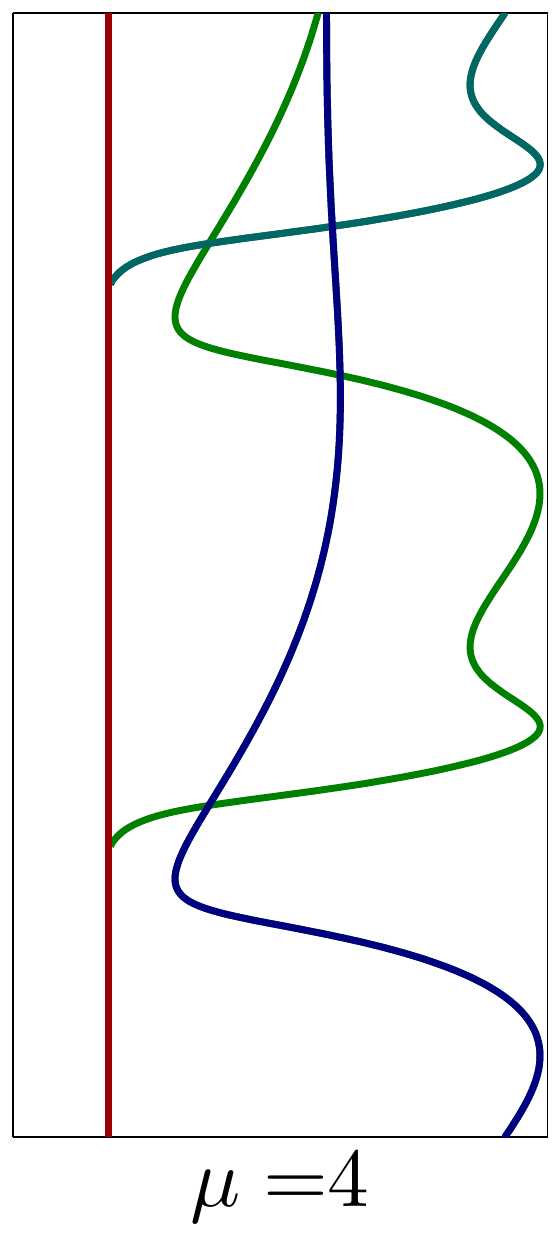}
}
\caption{Displacement-strain relation in the $(k,\varphi)$-plane for in the Swift-Hohenberg equation with boundary condition (\ref{e:free}), for various values of $\mu$, $k\in (0.6,1.1)$. The vertical lines indicate Eckhaus boundaries. Note the separation crisis at the Eckhaus boundary near $\mu=2.02$ followed by the reconnection with a different branch (near $\mu=2.15$).}
\label{f:shpk}
\end{figure}

In order to illustrate the phenomena discussed here, we computed the displacement-strain curves for free boundary conditions, \begin{equation}\label{e:free}
u_{xx}+u=0,\qquad u_{xxx}+u_x=0.
\end{equation}
Free boundary conditions arise naturally when considering the operator generated by the bilinear form $\int_0^\infty(\partial_{xx} u+u)(\partial_{xx}v+v)\rmd x$ on $H^2\times H^2$. Figure \ref{f:shpk} shows results of these computations for various values of $\mu$.  Numerical parameters were $dx=0.04, L=14*\pi$; for these parameters, we found errors in the displacement-strain relations below $10^{-4}$.  

Not surprisingly, we find phase selection for small $\mu$, reflecting the fact that the energy density along periodic patterns is almost constant in $x$. For moderate sizes of $\mu$, we still find wavenumber selection, with a characteristic double-well shape for $k>1$ that persists throughout the entire range of $\mu$. Note also that the displacement-strain curves are symmetric by $\varphi\mapsto \varphi+\pi$ but break the $\varphi\mapsto -\varphi$ symmetry. For $\mu\sim 2.02$, the displacement-strain curve touches the Eckhaus boundary and separates, the typical transition between wavenumber and phase selection. Shortly after, $\mu\sim 2.15$, the displacement-strain relation reconnects with a different branch, which is not periodic in $\varphi$, but rather ``snakes'' as $\varphi\to\infty$. The profiles indicate that this snaking behavior is mediated by the emergence of a stationary interface between the stable state $u=\sqrt{\mu-1}$ and the periodic patterns. This stationary interface interacts with the 
boundary condition in an oscillatory, exponentially weak fashion, yielding steady states with increasing separation of the interface from the boundary and exponentially small, oscillatory oscillations in the wavenumber.

\begin{figure}
\tiny
\begin{minipage}{0.26\textwidth}\includegraphics[height=1.9in]{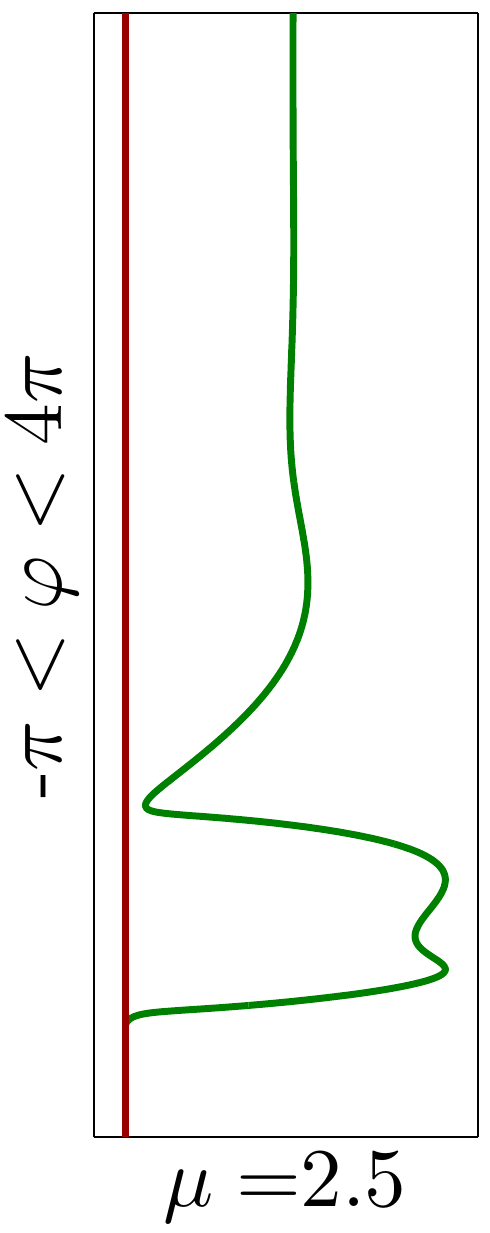}
\includegraphics[height=1.9in]{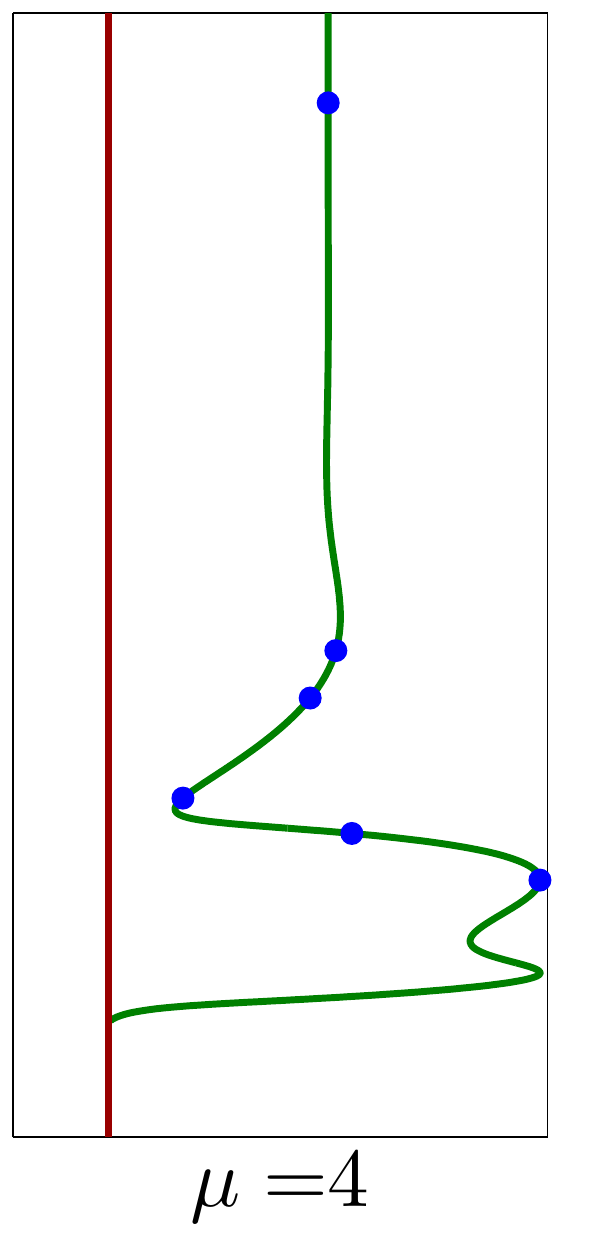}
\end{minipage}\hfill
\begin{minipage}{0.73\textwidth}
\raisebox{0.25in}{$\begin{array}{c}k=1.09\\\varphi=0.45\end{array}$}\includegraphics[width=0.365\textwidth]{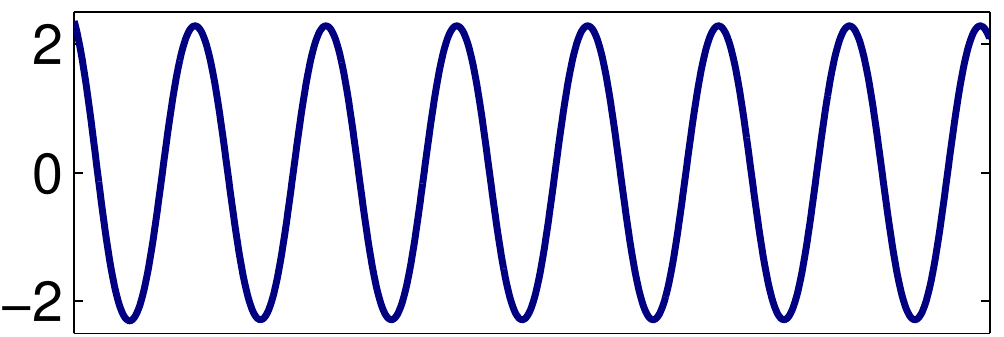}\qquad 
\includegraphics[width=0.365\textwidth]{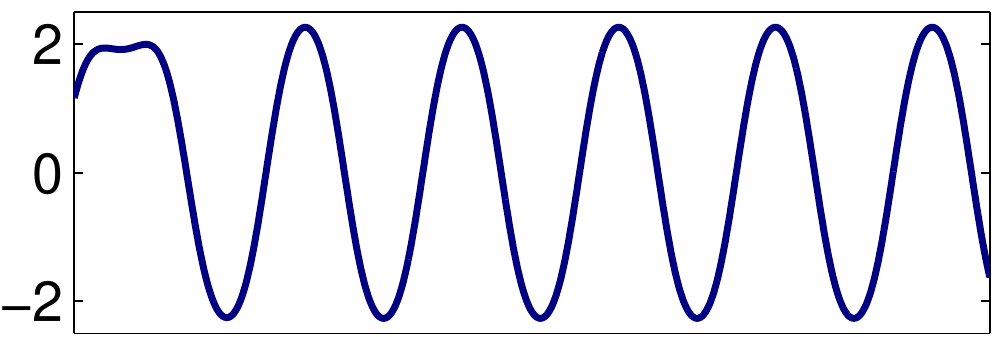}
\raisebox{0.25in}{$\begin{array}{c}k=0.91\\\varphi=3.00\end{array}$}\\
\raisebox{0.25in}{$\begin{array}{c}k=0.95\\\varphi=1.10\end{array}$}\includegraphics[width=0.365\textwidth]{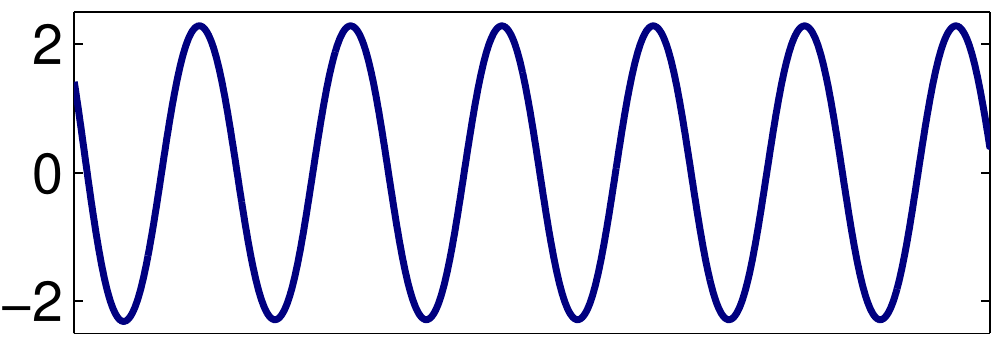}\qquad 
\includegraphics[width=0.365\textwidth]{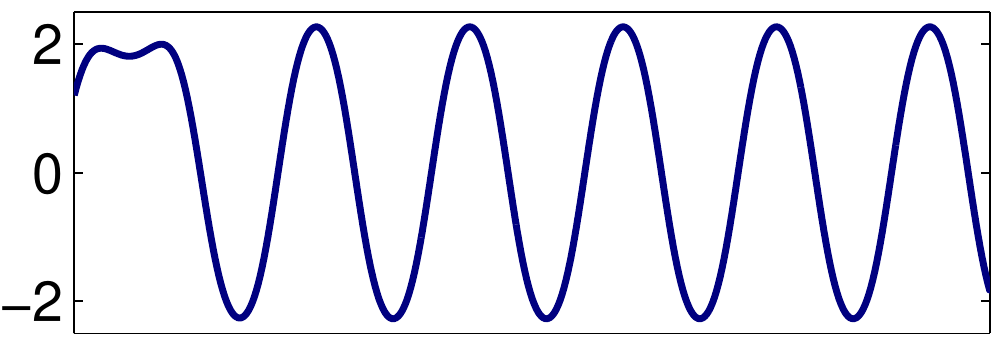}\raisebox{0.25in}{$\begin{array}{c}k=0.93\\\varphi=3.66\end{array}$}\\
\raisebox{0.25in}{$\begin{array}{c}k=0.83\\\varphi=1.60\end{array}$}\includegraphics[width=0.365\textwidth]{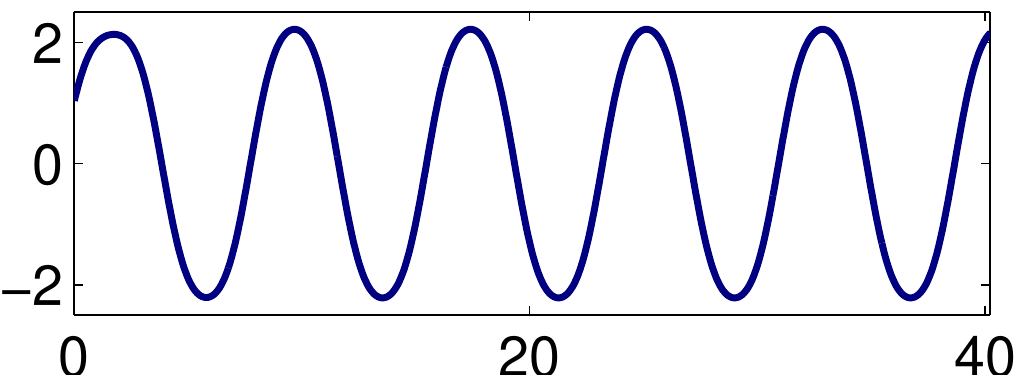}\qquad 
\includegraphics[width=0.365\textwidth]{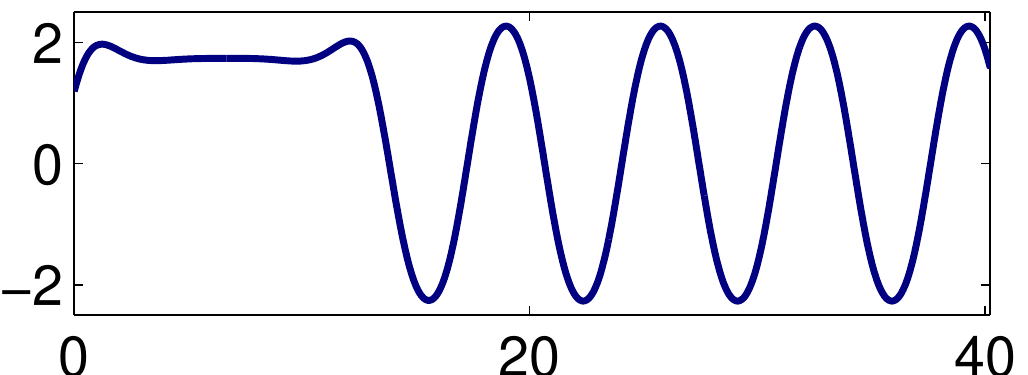}\raisebox{0.25in}{$\begin{array}{c}k=0.9276\\\varphi=11.31\end{array}$}
\end{minipage}
%
%
%
%
%
%

\caption{Displacement-strain relations as in Figure \ref{f:shpk} for $\mu=2.5,4$, extended for large phases (left). Plots of associated profiles for $\mu=4$ associated with points marked on the displacement-strain curve for $\mu=4$; $\varphi$ increasing top to bottom, left to right. One clearly notices the separation of an interface from the boundary along the snaking curve (right). }
\label{f:snake}
\end{figure}

\section{Discussion}\label{s:5}

Summarizing, we have introduced a systematic description of boundary layers in terms of curves that can be interpreted as displacement-strain relations. Generically, in terms of the boundary conditions, these are smooth curves which may terminate either near instability boundaries, or when the boundary layer itself disappears. We showed how such curves can be used to describe systematically the set of equilibria in bounded domains, in particular as the domain size is varied. More specifically, we defined an operation on pairs of displacement-strain relations corresponding to left and right boundary layers in large but finite domains, that consists of ``vertical subtraction'' and subsequent ``quantization''. The resulting families of intersection points are in one-to-one correspondence with equilibria in large bounded domains. 

Most of the paper is devoted to actually constructing or computing such displacement-strain curves. We saw that the common phase-diffusion approximation only gives very limited classes of examples! On the other hand, the Ginzburg-Landau equation, equipped with a variety of inhomogeneous Dirichlet, Neumann, or mixed boundary conditions gives a plethora of interesting phenomena, including phase- and wavenumber selection, non-monotone displacement-strain, pinch-off of defects leading to end points of displacement-strain curves, absence of boundary layers and isolas, and displacement-strain curves of higher winding number. We also outline how to systematically explore such displacement-strain curves using numerical continuation and illustrate some preliminary results in the Swift-Hohenberg equation. 

Our point of view raises numerous questions. First, we noted throughout that changes in the type of displacement-strain curve are accompanied by reconnection or separation events at the Eckhaus boundary. A local analysis of the stable manifold there should reveal a universal description of such events. One may also wonder if a variational structure in the boundary condition may lead to additional restrictions on boundary layers and bifurcations.We showed existence in the Ginzburg-Landau equation, but we suspect that, more generally, variational problems where periodic structures minimize the energy always accommodate boundary layers. 

Much of our analysis is reminiscent of the analysis of snaking diagrams in, for instance, the weakly subcritical Swift-Hohenberg \cite{snakes}. There, the trivial stable state coexists with stable periodic patterns and interfaces between trivial and periodic states provide ``effective boundary conditions''. Due to translation invariance, such boundary conditions are always pure wavenumber selecting, $k(\varphi)\equiv const$, so that interesting questions are concerned with changes of the selected wavenumber as parameters $\mu$ are varied. Fold points  of $k(\mu)$ can be interpreted as unpinning transition. The analysis in \cite{snakes} is concerned with phase matching of such ``effective boundary layers'', which yields spatially localized structures. Due to the singularities of displacement-strain curves, equilibria are not locally unique, generating various bifurcations (ladders, etc.).  It would clearly be very interesting to analyze phase matching as described here near singularities of displacement-
strain curves, and possibly elucidate the connection with snaking diagrams further. 

In higher dimensions, rolls allow for rotation as an additional degree of freedom. Setting up boundary layer problems in $x>0$, and $ y\in \R/L\Z$ periodic, we do encounter however a translational symmetry in $y$ which implies wavenumber selection due to arbitrary phases, for rolls that are not parallel to the boundary. The case of Neumann boundary conditions has been studied extensively, since boundary layers can then be interpreted as symmetric grain boundaries; see \cite{gb1,gb2} and references therein. 

Recognizing the analogy with defects such as grain boundaries, one can ask if there is a similar description for defects and their interactions. This appears to be an interesting avenue of research. Some results in two space-dimensions \cite{JS} indicate that properties such as phase selection can be attributed to inhomogeneities in higher-dimensional striped phases.

Generalizations of our setup also include discrete atomic lattices, possibly multi-atom lattices, which possess natural periodic equilibria. In a different direction, boundary conditions (and defects) for oscillatory media, where striped patterns propagate with non-zero group velocities, have been classified in \cite{SSdef}. It would be interesting to find a more unified description of all these interesting phenomena, and we hope that the present work is a contribution in this direction.

\paragraph{Acknowledgments.} The authors acknowledge discussions with David Lloyd on the implementation of the computational strategy outlined in Section \ref{s:4.1}. The decomposition described there was first implemented in a joint forthcoming work on grain boundaries \cite{lloyd} and the computational results in Section \ref{s:4.3} mimic the approach in \cite{lloyd}.

\end{document}